\documentclass[preprint,3p,12pt]{elsarticle}

\usepackage{amssymb}
\usepackage{latexsym}

\usepackage{url}
\usepackage{xcolor}
\definecolor{newcolor}{rgb}{.8,.349,.1}

\usepackage{bm}
\usepackage[ruled,vlined]{algorithm2e}
\usepackage[caption = false]{subfig}
\usepackage{amsmath,amsfonts}
\usepackage{hyperref}
\hypersetup{colorlinks=true}
\usepackage{capt-of}
\usepackage{array}
\usepackage{stackengine}

\newcommand{\real}{\mathbb{R}}
\newcommand{\complx}{\mathbb{C}}

\newcommand{\dv}[2]{\frac{\partial #1}{\partial #2}}
\newcommand{\ip}[2]{\langle #1, #2 \rangle}

\renewcommand{\vec}[1]{\bm{#1}} %comment this out to revert to classic 'over-arrow' notation
\newcommand{\mat}[1]{\mathbf{#1}}
\newcommand{\tens}[1]{\mathfrak{#1}}
\newcommand{\defn}{\stackrel{\Delta}{=}}
\newcommand{\grad}{\vec{\nabla}}

\begin{document}

\begin{frontmatter}

% TITLE
\title{Reduced order models for nonlinear radiative transfer based on moment equations and POD/DMD of Eddington tensor}%

% AUTHORS
\author{Joseph M. Coale}

\author{Dmitriy Y. Anistratov}

\address{Department of Nuclear Engineering, North Carolina State University, Raleigh, NC}
\address{jmcoale@ncsu.edu, anistratov@ncsu.edu}

\begin{abstract}
%%%
A new group of reduced-order models  (ROMs) for nonlinear thermal radiative transfer (TRT) problems is presented.
They are formulated by means of the nonlinear projective approach  and
data compression techniques.
The nonlinear projection is applied to the Boltzmann transport equation (BTE)
 to derive a hierarchy of low-order moment equations. The Eddington (quasidiffusion) tensor that
 provides exact closure for the system of  moment equations is
 approximated via one of several data-based methods of model-order reduction.
 These methods are the (i) proper orthogonal decomposition, (ii) dynamic mode decomposition (DMD),
 (iii) an equilibrium-subtracted DMD variant. Numerical results
  are presented to demonstrate the performance of these ROMs
for the simulation of evolving radiation and heat waves.
 Results show these models to be  accurate even with very low-rank representations of the Eddington tensor.
 As the rank of the approximation is increased, the errors  of solutions generated by the ROMs gradually decreases.
\end{abstract}


\end{frontmatter}

%\linenumbers

\section{Introduction}
Radiative transfer is the process in which energy is  transported  through the mechanisms of propagation, absorption and emission of photon radiation, and plays an essential role in many different physical phenomena.
These phenomena are present in a wide range of fields including plasma physics, astrophysics, atmospheric and ocean sciences, and high-energy-density physics \cite{drake-hedp,shu-astro,thomas-stamnes-atm}.
The  multiphysics models for each of those corresponding  phenomena involving radiation transport  (e.g. radiation-hydrodynamics problems) are described by complex systems of differential equations.
Solving this class of problems is associated with an array of fundamental challenges.
The systems of governing equations are generally (i) tightly coupled, (ii) strongly nonlinear, (iii) characterized by multiple scales in space-time, and (iv) high-dimensional.

\newpage
The Boltzmann transport equation (BTE) describes the propagation of particles in matter.
It is an essential part of multiphysics models describing physical systems in which energy redistribution is affected by radiation transport.
The BTE drives the dimensionality of these problems.
Its solution depends on 7 independent variables in 3D geometry and typically resides in a higher dimensional space than the other multiphysics equations it becomes coupled to.
Employing a reduced-order model (ROM) for the BTE can be an effective means to decrease computational costs associated with multiphysics problems involving radiative transfer.

Of the ROMs that have been previously developed for radiation transport, some of the most well-known include the flux-limited diffusion, $P_1$ and $P_{1/3}$ models \cite{olson-auer-hall-2000,morel-2000,simmons-mihalas-2000}.  The capabilities of these ROMs have been extensively studied and they remain as useful and computationally cheap methods for many applications. Even so, the accuracy of these models is limited. As such the development of ROMs for
radiation transport with the goal of achieving high accuracy while remaining computationally
efficient continues to be an active area of research.
Recently this research has been developing a new class of ROMs with the potential to meet this goal.
These methods are founded on leveraging the vast amounts of data available from experiments and simulations that have been amassed over the years, with the idea to take advantage of general model-order reduction techniques combined with some given databases to achieve a reduction in dimensionality.  Many such techniques are available to choose from \cite{benner-gugercin-willcox-2015,brunton-kutz-2019,lucia-beran-silva-2004}, some notable examples including:
(i) the proper orthogonal decomposition (POD) (a.k.a. principle component analysis (PCA) or the Karhunen-Lo\`eve expansion) \cite{aubry-1991,berkooz-holmes-lumley-1993,holmes-1996},
(ii) the dynamic mode decomposition (DMD) \cite{rowley-2009,schmid-2010,tu-rowley-2014,williams-kevrekidis-rowley-2015},
(iii) the proper generalized decomposition \cite{chinesta-2011} and (iv) balanced truncation \cite{moore-1981}. These techniques have seen extensive use in the fluid dynamics community for the modeling of general nonlinear flows \cite{kunisch-volkwein-2002,rowley-dawson-2017}, linearized flows \cite{rowley-2004}, compressible flows \cite{rowley-colonius-2004}, turbulence \cite{smith-moehlis-holmes-2005,berkooz-holmes-lumley-1993} and other applications \cite{taira-2017,bui-damodaran-willcox-2003}. Naturally the same techniques also have a wide range of applicability in the development of ROMs for
particle  transport, and have been used to model linear
particle  transport  problems \cite{choi-2021,behne-ragusa-morel-2019,hardy-morel-ahrens-2019,peng-mcclarren-frank-2020,huang-I-2021,huang-II-2021},
neutron transport in reactor-physics problems \cite{elzohery-roberts-2021,alberti-palmer-2018,cherezov-sanchez-joo-2018},
and nonlinear radiative transfer \cite{peng-chen-cheng-li-2021,soucasse-2019,jc-dya-m&c2021,jc-dya-m&c2019,jc-dya-tans2019}.

In this paper, we consider the nonlinear thermal radiative transfer (TRT) problem.
It is defined by  the BTE coupled with the material energy balance (MEB)  equation that describes energy exchange between  radiation and matter. This TRT problem  models  a supersonic radiative flow \cite{Moore-2015}.
It also serves as a useful platform for the development and testing of computational methods for the more general class of
radiation hydrodynamics problems
and retains all of the associated fundamental challenges as discussed earlier.
We present a novel class of ROMs  based on a combined approach incorporating nonlinear projective and data-based techniques of model order reduction.
These ROMs are constructed from a set  of  low-order equations for moments of the specific intensity
with data-driven approximate closures.

The ROMs are based on the multilevel quasidiffusion (QD) method \cite{gol'din-1964},
also known as the variable Eddington factor (VEF) method  \cite{auer-mihalas-1970}. This method is in essence a nonlinear method of moments that takes on a multigrid algorithm over the variables describing
particle frequency (energy) and direction of motion \cite{gol'din-1972,Goldin-sbornik-82,mihalas-FRH-1984,winkler-norman-mihalas-85,PASE-1986,dya-aristova-vya-mm1996,aristova-vya-avk-m&c1999,dya-vyag-nse-2011,at-dya-nse-2014,lrc-dya-pne-2017,dya-jcp-2019}.
 It is formulated by (i) the high-order BTE
 and (ii) a hierarchy of low-order QD (LOQD) equations
 for moments of the radiation intensity. This system is exactly closed through the use of
the Eddington (QD) tensor and other linear-fractional factors
 that are weakly dependent on the BTE solution.
Multiphysics equations, e.g. the MEB equation, are coupled to these low-order moment equations.
  This constitutes our initial reduction in dimensionality (with no approximation) and has been shown to give significant advantage compared to other methods in solving multiscale, multiphysical problems \cite{adams-larsen-2002}.

The hierarchy of LOQD equations coupled with the MEB equation
can be applied as a basis for model reduction
with the use of approximate closures.
The  LOQD equations are  coupled with the BTE through the  Eddington (QD) tensor and boundary factors
that define the exact closures.
A spectrum of ROMs can  thus be derived
by means of
  various approximations to the  Eddington tensor.
The classical approach is to apply a linear approximation of the radiation intensity in angle.
This reduces the   LOQD equations to the  $P_1$ equations.
The model based on the $P_{1/3}$ approximation is derived from the
$P_1$ equations
by modification of the time derivative term of the flux in the momentum equation using a factor of $\frac{1}{3}$ \cite{olson-auer-hall-2000}.
Another group of models are based on the variable Eddington factor method
that uses an approximation of the Eddington tensor by means of the first two  moments
of the intensity.
The Minerbo model is derived by means of a maximum entropy closure for the Eddington tensor.
The $M_N$ method applies  the maximum entropy closure for a system of $N$ moment equations \cite{levermore-1996,hauck-2011,hauck-2012}.
Hence, the Minerbo model is  the $M_1$ method.
Other commonly used models apply Kershaw, Wilson, Livermore closures
 \cite{kershaw-1976,wilson-1970,levermore-1984,m1-2017}.

The novel class of ROMs developed here find closure for the LOQD system by means of a data-informed approximation of the Eddington (QD) tensor. There exist many unique methods to create this approximation as discussed earlier, and as such the use of each particular technique defines a variant model under this class. We consider two such methods of approximation: the POD and DMD. Both of these methods seek an optimal reduced basis to project the dynamics of a given system onto a low-dimensional subspace, although the optimality condition differs.

The remainder of this paper is organized as follows. The TRT problem is defined in Section \ref{sec:trt}. In Section \ref{sec:quasidiffusion} we formulate the multilevel QD (MLQD) method. The developed class of ROMs is formulated in Section \ref{sec:rom}, followed by an overview of the POD and DMD in Section \ref{sec:data_methods}. Section \ref{sec:results} presents numerical results and analysis of performance of the ROMs on the well-known Fleck-Cummings test problem. In Section \ref{sec:conclusion} a brief discussion is given to close out the paper.

\section{Thermal Radiative Transfer} \label{sec:trt}
We consider the TRT problem  given by the multigroup BTE  \cite{mihalas-FRH-1984,zel-1966}
\begin{gather}
	\frac{1}{c}\dv{I_g}{t} + \vec{\Omega}\cdot\grad I_g + \varkappa_g(T)I_g = \varkappa_g(T)B_g(T), \label{bte_mg}\\
	\vec{r}\in \Gamma,\quad \vec{\Omega}\in\mathcal{S},\quad  g=1,\dots,N_g,\quad t\geq t_0, \nonumber\\
	I_g\big|_{\vec{r} \in\partial\Gamma} = I_g^{\text{in}}, \quad  \vec{\Omega}\cdot\vec{e}_n<0,
	\quad I_g \big|_{t=t_0}=I_g^0,  \label{bte_bc-ic}
\end{gather}

\noindent and the material energy balance (MEB) equation,
which models energy exchange between radiation and matter
\begin{equation}
	\dv{\varepsilon(T)}{t} = \sum_{g=1}^{N_g} \int_{4\pi} \varkappa_g(T)\big( I_g - B_g(T) \big)\ d\Omega,\quad 	\label{meb}
\end{equation}
\begin{equation}
T|_{t=t_0}=T^0.
	\label{meb_ic}
\end{equation}

\noindent Here $I_g=I_g(\vec{r},\vec{\Omega},t)$ is the intensity of radiation, $T=T(\vec{r},t)$ is the material temperature, $\varepsilon$ is the material energy density, $\varkappa_g$ is the material opacity, and $B_g$ is the Planckian black-body radiation distribution function given by
\begin{equation}
	B_g(T) = \frac{2h}{c^2}\int_{\nu_{g-1}}^{\nu_{g}}\ \frac{\nu^3 \ d\nu}{e^{\frac{h\nu}{kT}}-1},
\end{equation}

\noindent where $\nu_g$ is the upper boundary of the $g^\text{th}$ frequency group and $\nu_0=0$. $c$ is the speed of light, $\vec{r}\in\real^3$ is spatial position,
$\vec{\Omega}$ is the unit vector in the direction of particle motion, $g$ is the frequency group index,
$N_g$ is the number of frequency groups,
$t$ is time.
We denote
$\mathcal{S}=\{ \vec{\Omega}\in\real^3\ :\ |\vec{\Omega}|=1 \}$,
$\Gamma\subset\real^3$ is the spatial domain,
$\partial\Gamma$ is the boundary surface of $\Gamma$ and $\vec{e}_n$ is the outward-facing unit normal vector to $\partial\Gamma$.
The TRT problem \eqref{bte_mg} and  \eqref{meb}   neglects photon scattering, material motion and heat conduction.

%=================================================================================
%
%=================================================================================
\section{The Multilevel Quasidiffusion Method}\label{sec:quasidiffusion}

%=================================================================================
\subsection{Method Formulation}
The multilevel quasidiffusion method is formulated by means of a nonlinear projective approach.
The BTE is projected  in several stages  onto a sequence of subspaces to reduce dimensionality of the transport problem
and derive a closed hierarchy of low-order equations \cite{Goldin-sbornik-82}.
Note that the BTE solution is a 7-dimensional function of space, angle, frequency group and time.
At the first stage, the BTE is projected onto  the 5-dimensional  subspace of functions of space, frequency group, and time.
The  projection  operators  are given by
$\mathcal{P}_0 \stackrel{\Delta}{=}\ip{1}{\cdot}_\Omega$ and  $\mathcal{P}_1 \stackrel{\Delta}{=}\ip{\vec{\Omega}}{\cdot}_\Omega$, where
$\ip{w}{u}_\Omega \stackrel{\Delta}{=} \int_{4\pi} w u \ d\Omega$.
Applying $\mathcal{P}_0$ and $\mathcal{P}_1$ to Eq. \eqref{bte_mg} leads to the zeroth and first angular moments of the BTE
\begin{subequations}
	\begin{gather}
		\frac{1}{c}\dv{}{t}\ip{1}{I_g}_\Omega + \grad\cdot\ip{\vec{\Omega}}{I_g}_\Omega + \varkappa_g(T)\ip{1}{I_g}_\Omega = 4\pi\varkappa_g(T)B_g(T),\\
		\frac{1}{c}\dv{}{t}\ip{\vec{\Omega}}{I_g}_\Omega + \grad\cdot\ip{\vec{\Omega}\otimes\vec{\Omega}}{I_g}_\Omega + \varkappa_g(T)\ip{\vec{\Omega}}{I_g}_\Omega = 0.
	\end{gather}
	\label{ang_moments}
\end{subequations}

\noindent From these moment equations, a set of low-order equations is formulated for the first two angular moments of the radiation intensity:
(i) the group radiation energy density \linebreak $E_g (\vec{r},t)\stackrel{\Delta}{=} \frac{1}{c}\ip{1}{I_g}_\Omega$  and
(ii)  group radiation flux
 $\vec{F}_g (\vec{r},t) \stackrel{\Delta}{=} \ip{\vec{\Omega}}{I_g}_\Omega$.
Closure is defined between the group moment equations \eqref{ang_moments} and the BTE \eqref{bte_mg} by casting the group radiation pressure tensor $\mathcal{H}_g\stackrel{\Delta}{=} \ip{\vec{\Omega}\otimes\vec{\Omega}}{I_g}_\Omega$ (i.e. the second moment of $I_g$) by means of the Eddington (QD) tensor given by
\begin{equation}
	\boldsymbol{\mathfrak{f}}_g  = \frac{\int_{4\pi}\vec{\Omega}\otimes\vec{\Omega}I_g\ d\Omega}{\int_{4\pi}I_g\ d\Omega},
	\label{mgqd_tensor}
\end{equation}
to get
\begin{equation}
\mathcal{H}_g=c \boldsymbol{\mathfrak{f}}_g E_g \, .
\end{equation}
This yields the system of the multigroup
LOQD equations given by
\cite{gol'din-1964,auer-mihalas-1970}
\begin{subequations} \label{gp-loqd}
	\begin{gather}
		\dv{E_g}{t} + \grad\cdot\vec{F}_g + c\varkappa_g(T)E_g = 4\pi\varkappa_g(T)B_g(T), \label{gp_ebal}\\
		\frac{1}{c}\dv{\vec{F}_g}{t} + c\grad\cdot (\boldsymbol{\mathfrak{f}}_g E_g) + \varkappa_g(T) \vec{F}_g = 0. \label{gp_mbal}
	\end{gather}
	\label{mgqd}
\end{subequations}

\noindent
Equation \eqref{gp_ebal} is
 the group radiation energy balance equation,
 and equation \eqref{gp_mbal} is
  the group radiation momentum  balance  equation.
The boundary and initial  conditions (BCs and ICs) for
equations \eqref{gp-loqd} have the following form \cite{gol'din-1964,gol'din-1972}:
\begin{equation} \label{g-bc-ic}
 	\vec{e}_n \cdot  \vec{F}_g\big|_{\vec{r}\in\partial\Gamma} = c C_g\big(E_g\big|_{\vec{r}\in\partial\Gamma} - E_g^{\text{in}}\big) + F_g^{\text{in}},  \quad E_g\big|_{t=t_0} = E_g^0,\quad \vec{F}_g\big|_{t=t_0} = \vec{F}_g^0,
\end{equation}
with the group boundary factors defined as
\begin{equation}
	C_g = \frac{\int_{\vec{\Omega}\cdot\vec{e}_n>0}  \vec{e}_n \cdot  \vec{\Omega}I_g \ d\Omega}{\int_{\vec{\Omega}\cdot\vec{e}_n>0} I_g \ d\Omega},
	\label{mgqd_bcfac}
\end{equation}
and
\begin{equation}
E_g^{\text{in}}  =   \frac{1}{c} \int_{4\pi}I_g^{\text{in}} \ d\Omega , \quad
F_g^{\text{in}} =\int_{4\pi} \vec{e}_n \cdot  \boldsymbol{\Omega} I_g^{\text{in}} \ d\Omega   , \quad
 E_g^0  =  \frac{1}{c} \int_{4\pi}I_g^0 \ d\Omega  , \quad
\vec{F}_g^0 = \int_{4\pi} \boldsymbol{\Omega} I_g^0 \ d\Omega .
\end{equation}

At the second stage, the multigroup LOQD equations \eqref{gp-loqd} are projected onto the 4-dimensional subspace
of functions of space and time by
applying  the projection operator  $\mathcal{P}_{\nu} \stackrel{\Delta}{=}\ip{1}{\cdot}_{\nu}$, where $\ip{w_g}{u_g}_{\nu}\stackrel{\Delta}{=}\sum_{g=1}^{N_g}w_g u_g$. This yields
\begin{subequations}
	\begin{gather}
		\dv{}{t}\ip{1}{E_g}_{\nu} + \grad\cdot\ip{1}{\vec{F}_g}_{\nu} + c\ip{\varkappa_g(T)}{E_g}_{\nu} = 4\pi\ip{\varkappa_g(T)}{B_g(T)}_{\nu},\\
		\frac{1}{c}\dv{}{t}\ip{1}{\vec{F}_g}_{\nu} + c\grad  \cdot \ip{\boldsymbol{\mathfrak{f}}_g}{E_g}_{\nu} + \ip{\varkappa_g}{\vec{F}_g}_{\nu} = 0.
	\end{gather}
	\label{nu_moments}
\end{subequations}

\noindent To derive the low-order equations for  the total radiation energy density
$E(\vec{r},t) \defn \ip{1}{E_g}_{\nu}$
and total radiation flux
$\vec{F}(\vec{r},t) \defn \ip{1}{\vec{F}_g}_{\nu}$,
the following set of spectrum averaged quantities are introduced:

\begin{subequations}
\begin{gather}
	\bar{\boldsymbol{\mathfrak{f}}} =
\frac{1}{\sum_{g=1}^{N_g} E_g}  \sum_{g=1}^{N_g} \boldsymbol{\mathfrak{f}}_g E_g , \quad
	\bar{\varkappa}_E = \frac{\sum_{g=1}^{N_g} \varkappa_gE_g }{\sum_{g=1}^{N_g} E_g},\quad
	%%%
	\bar{\varkappa}_B = \frac{\sum_{g=1}^{N_g} \varkappa_gB_g }{\sum_{g=1}^{N_g} B_g},
	%%%
	\\[5pt]   \bar{\mathbf{K}}_{R} = \text{diag}\big(\bar{\varkappa}_{R,x},\ \bar{\varkappa}_{R,y},\ \bar{\varkappa}_{R,z}\big) \, ,
 \quad
 	\bar{\varkappa}_{R,\alpha} = \frac{\sum_{g=1}^{N_g} \varkappa_g|F_{\alpha, g}| }{\sum_{g=1}^{N_g} |F_{\alpha, g}|},\quad \alpha=x,y,z, \\[5pt]
	%%%
	\bar{\bm{\eta}}=\frac{1}{\sum_{g=1}^{N_g} E_g} \sum_{g=1}^{N_g} (\varkappa_g-\bar{\mathbf{K}}_{R})\vec{F}_g \, .
	%%%
\end{gather}
\label{gr_closures}
\end{subequations}

\noindent
As a result
 the effective grey LOQD equations are obtained and written as \cite{PASE-1986}
\begin{subequations}	
	\begin{gather}
		\dv{E}{t} + \grad\cdot\vec{F} + c\bar{\varkappa}_E E = c\bar{\varkappa}_B a_RT^4, \label{tot_ebal}\\
		\frac{1}{c}\dv{\vec{F}}{t} + c\grad \cdot (\bar{\boldsymbol{\mathfrak{f}}}E) + \bar{\mathbf{K}}_{R}\vec{F} + \bar{\bm{\eta}}E = 0. \label{tot_mbal}
	\end{gather}
	\label{gqd}
\end{subequations}

\noindent
The BCs and ICs for Eqs. \eqref{gqd} are defined by
\begin{equation} \label{grey-bc-ic}
		  \vec{e}_n \cdot \vec{F}\big|_{\vec{r}\in\partial\Gamma} = c\bar{C}\big(E\big|_{\vec{r}\in\partial\Gamma} - E^{\text{in}}\big) + F^{\text{in}}, \quad E\big|_{t=t_0} = E^0,\quad \vec{F}\big|_{t=t_0} = \vec{F}^0 \, ,
\end{equation}
where
\begin{gather}
	 \bar{C} = \left. \frac{\sum_{g=1}^{N_g} C_g\big(E_g - E_g^{\text{in}}\big) }
	{\sum_{g=1}^{N_g} \big(E_g  - E_g^{\text{in}}\big)} \right|_{\vec{r}\in\partial\Gamma},\\[5pt]
	%%%
	E^{\text{in}}  = \sum_{g=1}^{N_g} E_g^{\text{in}} , \quad
	F^{\text{in}} = \sum_{g=1}^{N_g} F_g^{\text{in}}, \quad
	E^0  = \sum_{g=1}^{N_g} E_g^0 , \quad
	\vec{F}^0 = \sum_{g=1}^{N_g} \vec{F}_g^0 .
\end{gather}
Lastly, the material energy balance equation \eqref{meb}
is cast in grey form
\begin{equation}
	\dv{\varepsilon(T)}{t} = c\bar{\varkappa}_E E - c\bar{\varkappa}_Ba_RT^4
	\label{grey_meb}
\end{equation}
to couple with the grey LOQD equations \eqref{gqd}.

Finally, the hierarchy of equations of the MLQD method for TRT problems  consists of
\begin{enumerate}
%%%
\item the multigroup BTE  for the group intensity $I_g$
\begin{equation} \label{a-bte}
\frac{1}{c} \frac{\partial I_g}{\partial t} + \mathcal{L}_g I_g =   Q_g \, , \quad
g=1,\dots,N_g \, ,
\end{equation}
where  $\mathcal{L}_g= \mathcal{L}_g \big[T \big]$
and $ Q_g=  Q_g(T)$ are given by Eqs. \eqref{bte_mg}  and \eqref{bte_bc-ic},
%%%
\item the multigroup LOQD  equations   for $E_g$ and $\vec{F}_g$
  defined by Eqs. \eqref{mgqd} and \eqref{g-bc-ic} that have the following general form:
\begin{equation} \label{a-mloqd}
\frac{\partial \boldsymbol{\varphi}_g}{\partial t} +
\mathcal{K}_g \boldsymbol{\varphi}_g =  \vec{q}_g \, , \quad
 \boldsymbol{\varphi}_g =  \begin{pmatrix}E_g \\ \vec{F}_g \end{pmatrix} \, , \quad g=1,\dots,N_g \, ,
\end{equation}
where
\begin{equation} \label{Kg-oper}
\mathcal{K}_g= \mathcal{K}_g \big[\boldsymbol{\tens{f}}_g, C_g, T \big] \, , \quad
\vec{q}_g =  \vec{q}_g(T) \, ,
\end{equation}
\begin{equation} \label{f_g}
\boldsymbol{\tens{f}}_g=\boldsymbol{\tens{f}}_g\big[I_g\big] \, ,  \quad
C_g = C_g\big[I_g\big] \, ,
\end{equation}
%%%
\item the effective grey  LOQD  equations   for $E$ and $\vec{F}$ of the form
\begin{equation} \label{a-gloqd}
\frac{\partial \boldsymbol{\varphi}}{\partial t} + \mathcal{\bar K} \boldsymbol{\varphi} =  \vec{\bar q}\, , \quad
 \boldsymbol{\varphi} = \begin{pmatrix}E \\ \vec{F} \end{pmatrix} \, ,
\end{equation}
where
\begin{equation} \label{K-bar-oper}
\mathcal{\bar K}= \mathcal{\bar K} \big[\boldsymbol{\tens{\bar f}} , \bar C, \bar{\varkappa}_E, \bar{\mathbf{K}}_{R}, \bar{\bm{\eta}}, T \big] \, , \quad
  \vec{\bar q} =  \vec{\bar q}(T) \,
\end{equation}
    are defined by   Eqs. \eqref{gqd} and \eqref{grey-bc-ic},
%%%
\item the effective grey MEB equation for $T$  and $E$
 given by Eq. \eqref{grey_meb}.
\end{enumerate}

The components the group Eddington   tensor  $\boldsymbol{\mathfrak{f}}_g$
and the boundary factor $C_g$  are compressed data of the high-order solution of the BTE \eqref{a-bte}.
These data   carry all information  about the BTE solution that the hierarchy of
  the low-order equations  \eqref{a-mloqd} and \eqref{a-gloqd} needs
   to accurately describe radiative transfer physics.
In this multilevel  system of equations, the high-order BTE \eqref{a-bte}  can be interpreted
 as the  one  that generates  the shape function  for averaging $\vec{\Omega}\otimes\vec{\Omega}$ and   calculation $\boldsymbol{\mathfrak{f}}_g$ and $C_g$.
The  role of the  low-order equations is to generate the moments of the transport solution and communicate with the energy balance equation as an element of a multiphysics model.

%=================================================================================
\subsection{Discretization}

\begin{figure}[h]
\centering
	{\includegraphics[width=.4\textwidth]{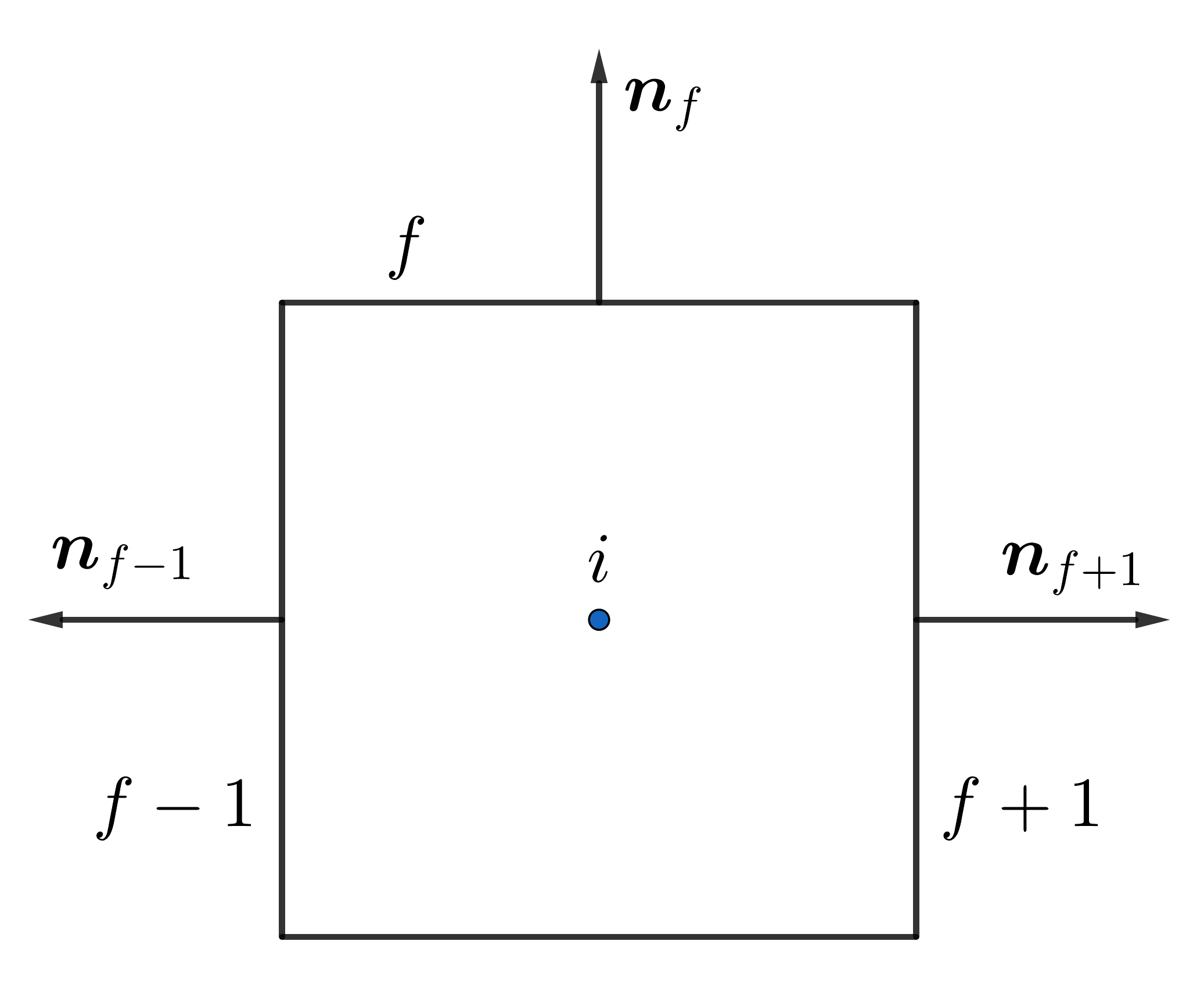}}
\caption{\label{cell-c} Notations in the cell $i$.}
\end{figure}

In this paper we consider TRT problems  in 2D Cartesian geometry.
To discretize  the multigroup LOQD equations \eqref{gp-loqd}
we apply fully implicit temporal approximation based on the Backward Euler (BE) scheme
and  a second-order  finite volume scheme  in space on orthogonal spatial grids \cite{ea-avk-1993,pg-dya-jcp-2020}.
Figure \ref{cell-c} shows  a sample  spatial cell $i$ and notations.
The multigroup radiation energy balance equation \eqref{gp_ebal} is integrated over the  cell $i$.
The multigroup  radiation  momentum balance  equations \eqref{gp_mbal}  are integrated over each half of the spatial cell.
The resultant  discretized multigroup LOQD equations
are given by
  \begin{subequations}\label{gp-loqd-d}
  \begin{equation}\label{gp-loqd0-d}
   \frac{A_i}{\Delta  t^n} \Big( E_{g,i}^n  - E_{g,i}^{n-1}\Big)
+  \sum_{f \in \omega_i} F_{g,f}^n \ell_f
+ c\varkappa_{g,i}^n E_{g,i}^n  A_i
= 4\pi \varkappa_{g,i}^n B_{g,i}^n A_i \, ,
  \end{equation}
\begin{multline}\label{gp-loqd1-d}
   \frac{A_f}{c\Delta  t^n} \Big( F_{g,f}^n  - F_{g,f}^{n-1}\Big)
+    \mathbf{e}_{\alpha}  \cdot \mathbf{n}_f  c \Big( \mathfrak{f}_{\alpha\alpha,g,f}^nE_{g,f}^n - \mathfrak{f}_{\alpha\alpha,g,i}^nE_{g,i}^n \Big) \ell_f \\
 +   \mathbf{e}_{\beta}  \cdot \mathbf{n}_{f+1}  c \Big(\mathfrak{f}_{\alpha\beta,g,f+1}^nE_{g,f+1}^n  -   \mathfrak{f}_{\alpha\beta,g,f-1}^nE_{g,f-1}^n \Big) \frac{\ell_{f+1}}{2}
+ \varkappa_{g,i}^n F_{g,f}^nA_f = 0 \, , \\
\quad
\alpha, \beta = x,y   \, , \quad
\beta  \ne \alpha \, ,
\end{multline}
\end{subequations}
where
\begin{equation}
	\boldsymbol{\mathfrak{f}}_g =
 \begin{pmatrix}
		\mathfrak{f}_{xx,g} &\mathfrak{f}_{xy,g}\\
		\mathfrak{f}_{xy,g} &\mathfrak{f}_{yy,g}
\end{pmatrix},
\end{equation}
$i$ is the cell index; $f$ is the  index of faces
 of the $i^\text{th}$ cell;
  $\ell_f$ is the length of the face $f$;
$\omega_i$ is the set  of faces of the  $i^\text{th}$ cell,
$\mathbf{n}_f$  is the unit outward normal of the cell face $f$ and
 $\mathbf{n}_f  =\mathbf{e}_{\alpha}$ for the orthogonal grids;
    $E_{g,i}$  and  $E_{g,f}$  are cell-average and face-average radiation energy densities, respectively;
  $F_{g,f} = \mathbf{n}_f  \cdot  \mathbf{F}_g $ is the normal component of the radiation flux;
    $A_i$ is the area of  the  $i^\text{th}$ cell;
  $A_f$ is the area of the half-cell associated with the edge $f$;
  $n$ is the index of the instant of time; $\Delta t^n = t^n - t^{n-1}$ is the $n^\text{th}$ time step.

The discretization of the effective grey LOQD equations \eqref{gqd} is algebraically consistent with the scheme \eqref{gp-loqd-d}.
 The discrete total radiation energy balance equation is derived by applying the projection operator $\mathcal{P}_{\nu}$ to Eq. \eqref{gp-loqd0-d} and takes the following form
 \begin{equation}\label{grey-loqd0-d}
\frac{A_i}{ \Delta t^n}\Big(E_i^n - E_i^{n-1} \Big) + \sum_{f \in \omega_i} F_f^n \ell_f +   c \bar \varkappa_{E,i}^n E_i^n A_i =
c \bar \varkappa_{B,i}^n a_R (T_i^n)^4    \, ,
 \end{equation}
 where
 \begin{equation} \label{bar-kappa-n}
\bar \varkappa_{E,i}^n  = \frac{\sum_{g=1}^{N_g} \varkappa_{g,i}^n E_{g,i}^n}{\sum_{g=1}^{N_g} E_{g,i}^n} \, , \quad
\bar \varkappa_{B,i}^n = \frac{\sum_{g=1}^{N_g} \varkappa_{g,i}^n B_{g,i}^n}{\sum_{g=1}^{N_g} B_{g,i}^n}\, .
\end{equation}
Next the discretized multigroup LOQD equation  \eqref{gp-loqd1-d} is cast in terms of $F_{g,f}^n$.
 Applying $\mathcal{P}_{\nu}$ to the resulting equation for $F_{g,f}^n$ yields an equation for the total radiation flux at each cell face
\begin{multline}\label{gloqd1-d}
 F_{f}^n =
  - c \bigg(   \pmb{e}_{\alpha}  \cdot \pmb{n}_f \Big( \mathcal{\bar D}_{\alpha\alpha,f}^n E_{f}^n  - \mathcal{\bar  D}_{\alpha\alpha,i}^n  E_{i}^n  \Big) \ell_f
 \\+ \frac{1}{2}   \pmb{e}_{\beta}  \cdot \pmb{n}_{f+1} \Big(\mathcal{\bar  D}_{\alpha\beta,f+1}^n  E_{g,f+1}^n    -   \mathcal{\bar  D}_{\alpha\beta,f-1}^n  E_{f-1}^n  \Big) \ell_{f+1}\bigg) + p_f^n \, ,
 \end{multline}
  where
  \begin{equation}
    \mathcal{\bar D}_{\alpha\beta,f'}^n  = \frac{\sum_{g=1}^{N_g}(\tilde \varkappa_{g,i}^n)^{-1} \mathfrak{f}_{\alpha\beta,i}^n E_{g,f'}^n}{\sum_{g=1}^{N_g} E_{g,f'}^n} \, \quad
\quad f'=f -1,f,f+1 \, ,
   \end{equation}
   \begin{equation}
  \tilde \varkappa_{g,i}^n  = \varkappa_{g,i}^n   + \frac{1}{c \Delta t^n} \, ,
\quad
   p_f^n = \sum_{g=1}^{N_g}\frac{F_{g,f}^{n-1}}{1+ c \Delta t^n\varkappa_{g,i}^n  } \, .
 \end{equation}

To discretize the BTE  \eqref{bte_mg}  in the angular variable we use
the method of discrete ordinates. The BE  scheme  is applied  for time integration,
 and   the simple corner-balance   method is used  for  approximation of the BTE  in space \cite{adams-1997}.
There are different methods for solving the hierarchy of equations of MLQD method
\cite{dya-aristova-vya-mm1996,aristova-vya-avk-m&c1999,dya-jcp-2019,pg-dya-jcp-2020,anistratov-2021}.

%=================================================================================
%
%=================================================================================
\section{Model Reduction with Data-Informed Closures} \label{sec:rom}

In this study, we develop  ROMs for TRT on the basis  of  the hierarchy of LOQD equations
   \eqref{a-mloqd}  and \eqref{a-gloqd}  coupled with the MEB equation \eqref{grey_meb}
   with the Eddington tensor approximated by   data-driven techniques
 using  available data.
This forms a class of ROMs
henceforth referred to as
data-driven Eddington tensor (DDET) ROMs.
The Eddington tensor data can be generated by an array of means, including (i) the full-order model (FOM) solution (i.e. Eqs. \eqref{bte_mg} \& \eqref{meb}) for some set of base-case (reference) problems the parameters of which cover a desired range, (ii) the BTE solution obtained by low-cost (coarse-mesh) calculations of  TRT problems that approximate well the radiation transport effects.

\begin{figure}[ht!]
	\centering
	\subfloat[Cell-wise grid functions]{\includegraphics[width=.4\textwidth]{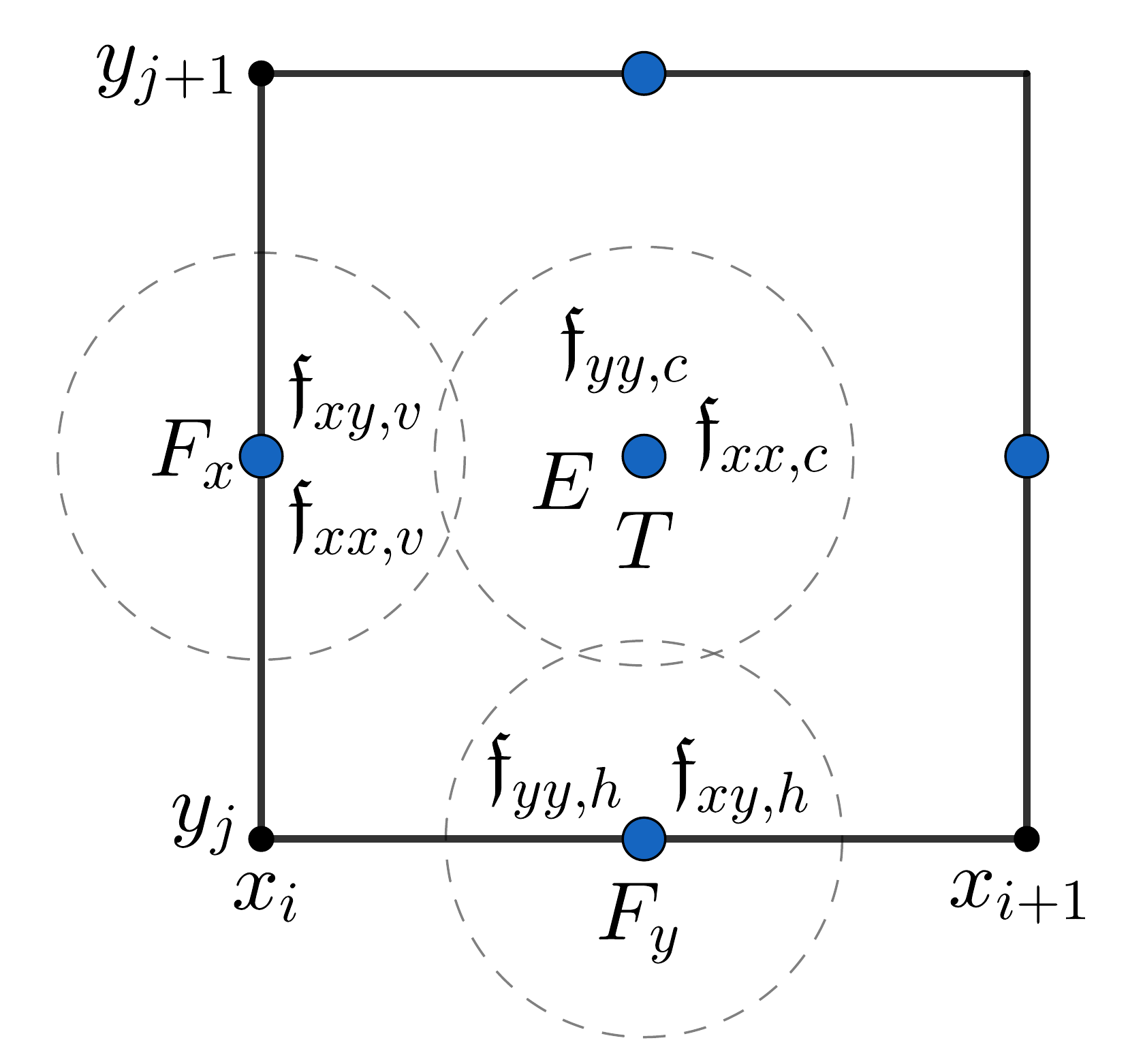} \label{subfig:qd_grid}}
	\subfloat[Boundary factors]{\includegraphics[width=.4\textwidth]{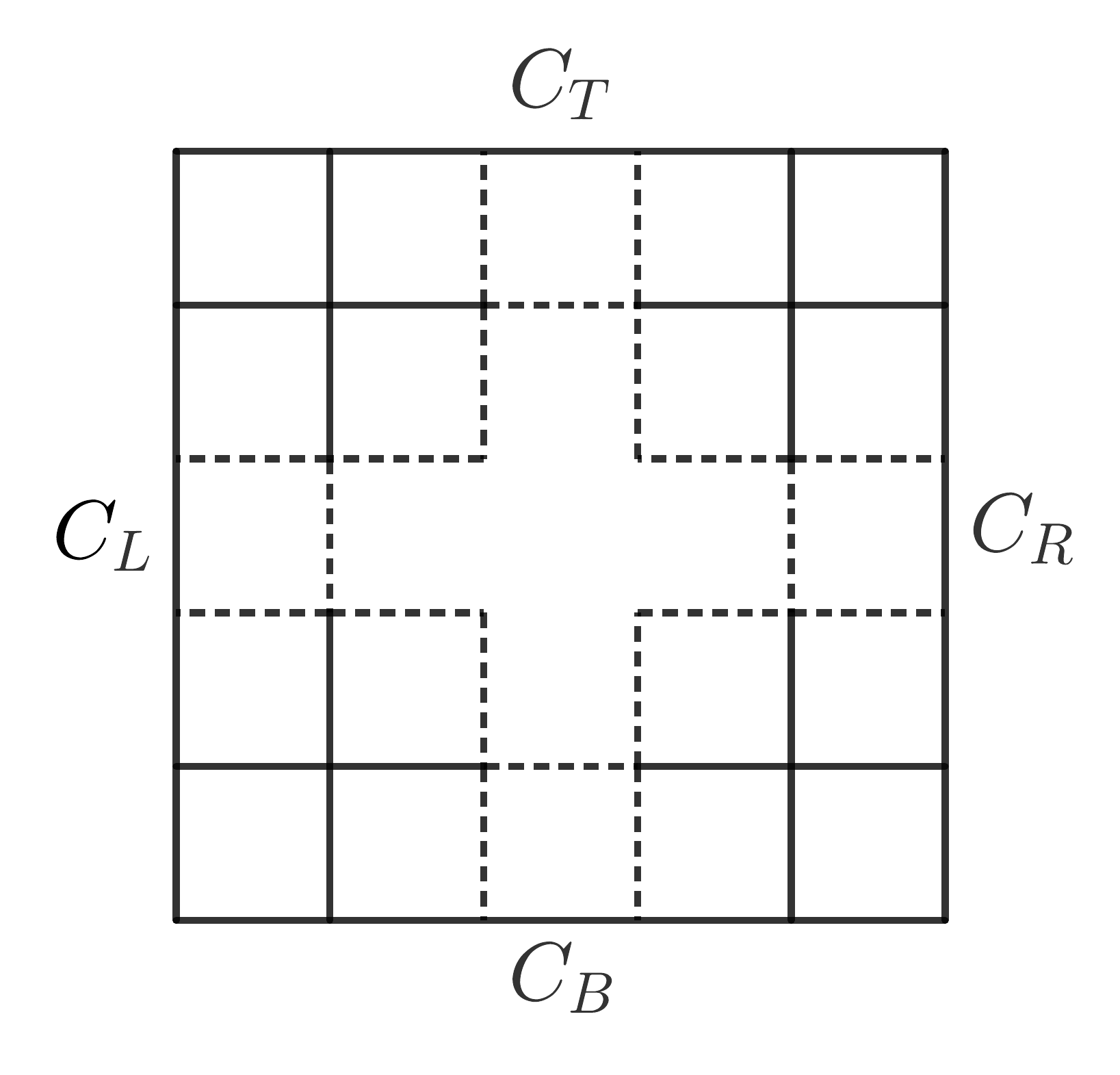} \label{subfig:bc_grid}}
	\caption{Discrete grid functions shown on a sample 2D spatial mesh}
	\label{fig:2D-cell-grid-functions}
\end{figure}

\newpage
2D TRT problems on orthogonal grids are considered, where $N_x$ and $N_y$  are the number of spatial cells in the
$\pmb{e}_{x}$- and $\pmb{e}_{y}$-directions, respectively.
Vectors of each component of the group Eddington tensor and boundary factors are constructed at every $n^\text{th}$ time step over the entire spatial grid as
\begin{subequations} \label{vec-f-group}
\begin{gather}
 	 \bm{\mathfrak{f}}_{\alpha\beta,g,\gamma}^n\in\real^{D_\gamma},\quad
\gamma=  v,h,c,
 \quad \alpha, \beta=x,y,  \\
	\vec{C}_{\theta,g}^n\in\real^{D_\theta},\quad \theta=L,B,R,T
\end{gather}
\end{subequations}
with dimensions
\begin{subequations}
\begin{gather}
	D_L=D_R=N_y,\quad D_B=D_T=N_x,\\
	D_v = (N_x+1)N_y,\quad D_h = N_x(N_y+1),\quad D_c=N_xN_y.
\end{gather}
\end{subequations}
Notations are illustrated in Figure \ref{fig:2D-cell-grid-functions}.
The vectors \eqref{vec-f-group} are subsequently `stacked' by frequency group  to construct
complete data vectors for each  quantity at the time step
\begin{subequations}
\begin{gather}
	   \mathbf{f}_{\alpha\beta,\gamma}^n = \left(\begin{array}{c}
		\bm{\mathfrak{f}}_{\alpha\beta,1,\gamma}^n\\[2.5pt]
		\bm{\mathfrak{f}}_{\alpha\beta,2,\gamma}^n\\
		\vdots\\
		\bm{\mathfrak{f}}_{\alpha\beta,N_g,\gamma}^n
	\end{array}\right) \in\real^{N_gD_\gamma},\\[5pt]
	\vec{C}^n = \left(\begin{array}{c}
		\vec{C}_1^n\\[2.5pt]
		\vec{C}_2^n\\
		\vdots\\
		\vec{C}_{N_g}^n
	\end{array}\right) \in\real^{2N_g(N_x+N_y)} ,\quad
	\vec{C}_g^n = \left(\begin{array}{c}
		\vec{C}_{L,g}^n\\[2.5pt]
		\vec{C}_{B,g}^n\\[2.5pt]
		\vec{C}_{R,g}^n\\[2.5pt]
		\vec{C}_{T,g}^n
	\end{array}\right) \in\real^{2(N_x+N_y)}.
\end{gather} \label{nvecs}
\end{subequations}

\noindent Finally the vectors \eqref{nvecs} are collected in chronological order as columns of the following snapshot matrices
\begin{equation}
	\mat{A}^{\mathfrak{f}_{\alpha\beta,\gamma}} = [ \mathbf{f}_{\alpha\beta,\gamma}^1\ \mathbf{f}_{\alpha\beta,\gamma}^2\ \dots\  \mathbf{f}_{\alpha\beta,\gamma}^{N_t}] \in\real^{N_gD_\gamma\times N_t},\quad
	%%%
	\mat{A}^c = [ \vec{C}^1\ \vec{C}^2\ \dots\ \vec{C}^{N_t} ] \in\real^{2N_g(N_x+N_y)\times N_t},
	\label{qd_snapmats}
\end{equation}

\noindent where $N_t$ is the total number of time steps. Each of these snapshot matrices is then projected onto some low-rank subspace as
\begin{equation}
	\mathcal{A}_k^\phi= \mathcal{G}^k\mat{A}^\phi,\quad k\leq \text{rank}(\mat{A}),\quad \phi=\mathfrak{f}_{\alpha\beta,\gamma},\ C
\end{equation}
where the projection operator  $\mathcal{G}^k$ is defined by the specific method of data compression.
The rank-$k$ representation of the matrix $\mat{A}^\phi$ of the form $\mathcal{A}_k^\phi\in\complx^{D_A\times (k+\delta)}$ is constructed of $k$ sets of various vectors and factors (all of which may be complex), which in total comprise $D_A$ elements. Depending on the specific method used, $\delta$ can be 0 or 1.

\newpage The method that specifies $\mathcal{G}^k$ similarly defines a map $\mathcal{M}^n:\complx^{D_A\times (k+\delta)}\rightarrow \real^{D_\phi}$, where
\begin{equation}
	D_\phi=\left\{\begin{array}{ll}
		N_gD_\gamma, & \phi=\mathfrak{f}_{\alpha\beta,\gamma}   \\
		2N_g(N_x+N_y), & \phi=C
\end{array}\right. ,
\end{equation}
such that
\begin{equation}
 \vec{\Phi}^n = \mathcal{M}^n\mathcal{A}_k^\phi ,\quad \vec{\Phi}^n\in\real^{D_\phi},\quad \vec{\Phi}^n\approx\left\{\begin{array}{ll}
 	\mathbf{f}^n_{\alpha\beta,\gamma}, & \phi=\mathfrak{f}_{\alpha\beta,\gamma} \\
 	\vec{C}^n, &\phi=C
 \end{array}\right.
\end{equation}
 The data used to form each $\mathcal{A}_k^\phi$ is generated from the FOM solution obtained by means of the MLQD method on a given grid in the phase space and time. This yields an  Eddington tensor and boundary factors over the discrete spatial domain and frequency groups at time $t^n$.

The method for solving TRT problems with the DDET ROMs is summarized in Algorithm \ref{alg:rom_alg}.
At each $n^\text{th}$ time step, the corresponding grid functions of the Eddington tensor and boundary factors are generated by applying the operator $\mathcal{M}^n$ to the input compressed data $\mathcal{A}^\phi_k$.
Then the material temperature and total radiation energy density are calculated iteratively.
 At each $\ell$ iteration the  compressed representations of the Eddington tensor and boundary factor data are used
 to define the multigroup LOQD equations. The solution to the multigroup LOQD equations is used to generate all effective grey opacities and factors. The effective grey problem formed by the coupled effective grey LOQD and MEB equations is subsequently solved via Newton's method to obtain the material temperature, total radiation energy density and total radiation flux.

\begin{algorithm}[ht!]
	\SetAlgoLined
	Input: $\mathcal{A}^{\mathfrak{f}_{\alpha\beta,\gamma}}_k, \mathcal{A}^c_k$\\
	$n=0$\\
	\While{$t_n \leq t^{\text{end}}$}{
		$n = n+1$  \\
		$\mathcal{M}^n \mathcal{A}^{\mathfrak{f}_{\alpha\beta,\gamma}}_k   \ \Rightarrow \  \bm{\tens{f}}_g^n$, \
		$\mathcal{M}^n \mathcal{A}^c_k  \ \Rightarrow \ \vec{C}_g^n $   \\\vspace*{.05cm}
		$\ell=-1$ \\
		$T^{n \,(0)} = T^{n-1}$ \\
		\While{ $\|T^{n \, (\ell+1)} - T^{n \,(\ell)}\| > \epsilon_1\|T^{n \,(\ell+1)}\| + \epsilon_2 \, \lor \, \|E^{n \,(\ell+1)} - E^{n \,(\ell)}\| > \epsilon_1\|E^{n \,(\ell+1)}\| + \epsilon_2$ }{
			$\ell=\ell+1$\\
			Update $B_g^{n \, (\ell)}, \varkappa_g^{n \, (\ell)}, \bar{\varkappa}_B^{n \, (\ell)}$ with $T^{n  \, (\ell)}$\\
			Solve multigroup LOQD equations  \eqref{mgqd} for $E_g^{n \, (\ell)}, \ \vec{F}_g^{n \, (\ell)}$\\
			Compute spectrum-averaged coefficients $ \bar{\varkappa}_E^{n \, (\ell)},\ \bar{\varkappa}_B^{n \, (\ell)}, \
			\bar{\vec{C}}^{n \, (\ell)},\ \bar{\bm{\tens{f}}}^{n \, (\ell)},\ \bar{\mathcal{D}}^{n \, (\ell)}$ \\
			Solve effective grey problem \eqref{gqd} and \eqref{grey_meb} for $T^{n \, (\ell+1)}, \ E^{n \, (\ell+1)}, \ \vec{F}^{n \, (\ell+1)}$\\
		}
		$T^n \leftarrow T^{n \, (\ell+1)}$
		
	}
	\caption{Algorithm for solving equations of the DDET ROMs for TRT problems
		\label{alg:rom_alg}}
\end{algorithm}

%=================================================================================
%
%=================================================================================
\section{Projection of Data onto Low-Rank Subspaces} \label{sec:data_methods}

%=================================================================================
%
%=================================================================================
\subsection{Proper Orthogonal Decomposition}
Let $\{\vec{a}^n\}_{n=0}^{m}$ be a set of data-vectors $\vec{a} \in \real^d$ such that $\vec{a}^n=\vec{a}(t^n)$ for $n=0,\dots,m$
at some set of instances $\{t^n\}_{n=0}^m$.
The POD seeks an orthonormal basis $\{\vec{u}_1, \vec{u}_2, \dots, \vec{u}_k\}$ onto which the set of zero-mean data $\{\hat{\vec{a}}^n\}_{n=0}^m$, defined as
\begin{equation}
	\hat{\vec{a}}^n=\vec{a}^n-\bar{\vec{a}},\quad \bar{\vec{a}} = \frac{1}{m+1}\sum_{n=0}^{m}\vec{a}^n
\end{equation}

\noindent can be projected in a way that POD modes optimally capture the energy associated with the given data set \cite{aubry-1991}.
 This optimality condition amounts to finding the projection of $\{\vec{a}^n\}_{n=0}^m$ onto $\text{span}\{\bar{\vec{a}}, \vec{u}_1, \vec{u}_2, \dots, \vec{u}_k\}$ with minimal error in the Frobenius norm.
The rank-$k$ POD of the data $\{\vec{a}^n\}_{n=0}^m$ (i.e. $\mathcal{A}_k$) is thus constructed from the data mean $\bar{\vec{a}}$, POD basis $\{\vec{u}_i\}_{i=1}^k$ and coefficients
\begin{equation}
	\alpha_{j}^n = \ip{\hat{\vec{a}}^n}{\vec{u}_j},\quad  n=0,\dots,m,\quad j=1,\dots,k,
\end{equation}

\noindent and approximates each data point $\vec{a}^n$ through the operator $\mathcal{M}^n$ as follows:
\begin{equation}
	\vec{a}^n\approx \mathcal{M}^n\mathcal{A}_k = \bar{\vec{a}} + \sum_{j=1}^{k}\alpha^n_{j}\vec{u}_j.
\end{equation}

Let us define the matrices $\mat{A}\in\real^{d\times (m+1)}$ whose columns are given by the vectors $\{\hat{\vec{a}}^n\}_{n=0}^m$.
The thin (reduced) singular value decomposition (SVD) of $\mat{A}$ is
\begin{equation}
	\mat{A} = \mat{U}\mat{S}\mat{V}^\top,
\end{equation}
\noindent where $\mat{U}\in\real^{d\times r}$ and $\mat{V}\in\real^{(m+1)\times r}$ hold the first $r$ left and right singular vectors of $\mat{A}$ in their columns, respectively, and $\mat{S}\in\real^{r\times r}$ holds the $r$ nonzero singular values of $\mat{A}$ along its diagonal in descending order, where $r=\text{rank}(\mat{A})$. The POD basis vectors are found as the first $k$ columns of $\mat{U}$ and the projection coefficients are $\alpha^n_{j} = v_{n,j}\sigma_j$, where $\sigma_j$ is the $j^\text{th}$ singular value of $\mat{A}$ and $v_{n,j}$ is the $(n,j)^\text{th}$ element of $\mat{V}$. An efficient compression of $\mat{A}$ is then constructed with the rank-$k$ truncated SVD (TSVD) of $\mat{A}$
\begin{equation}
	\mat{A} \approx \mat{A}_k = \mat{U}_k\mat{S}_k\mat{V}_k^\top, \quad k\leq r
\end{equation}

\noindent where $\mat{U}_k\in\real^{d\times k}$ and $\mat{V}_k\in\real^{(m+1)\times k}$ hold the first $k$ left and right singular vectors of $\mat{A}$ in their columns, respectively, and $\mat{S}_k\in\real^{k\times k}$ holds the first $k$ singular values of $\mat{A}$ along its diagonal in descending order.
$\mat{A}_k$ is actually the orthogonal projection of $\mat{A}$ onto $\{\vec{u}_i\}_{i=1}^k$, written as $\mat{A}_k= \mat{U}_k\mat{U}_k^\top\mat{A}$.
The error introduced by this orthogonal projection is
given by \cite{ipsen-2009}
\begin{equation}
	\xi^2 = \|\mat{A} - \mat{U}_k\mat{U}_k^\top\mat{A}\|_F^2 = \sum_{i=k+1}^r\sigma_i^2.
\end{equation}
The relative error of the POD approximation in the Frobenius norm is therefore
\begin{equation}
	\xi^2_{\text{rel}} = \frac{\|\mat{A}-\mat{A}_k\|_F^2}{\|\mat{A}\|_F^2} = \frac{\sum_{i=k+1}^r\sigma_i^2}{\sum_{i=1}^r\sigma_i^2}. \label{PODerr_rel}
\end{equation}

\noindent $\xi^2_{\text{rel}}$ can be interpreted as the ratio of energy encompassed by the first $k$ POD modes to the total energy comprised by all POD modes of the given data \cite{smith-moehlis-holmes-2005}. When the POD is performed, we choose some desired value for $\xi_{\text{rel}}$ and find the rank $k$ that satisfies the expression \eqref{PODerr_rel}.

%=================================================================================
%
%=================================================================================
\subsection{Dynamic Mode Decomposition}

Let us consider the case when the time instances $\{t^n\}_{n=0}^m$  are uniformly spaced such that
\begin{equation}
  	t^{n+1}=t^n+\Delta t,\quad n=0,\dots,m-1 \, .
\end{equation}

\noindent The DMD constructs the best-fit linear operator $\mat{B}$ to the data $\{\vec{a}^n\}_{n=0}^m$, generating the following dynamic system:
\begin{equation}
	\frac{d \tilde{\vec{a}}(t)}{dt} = \mat{B}\tilde{\vec{a}}(t),
\end{equation}

\noindent whose solution $\tilde{\vec{a}}(t)$  approximates $\vec{a}(t)$ and  is given by
\begin{equation}
	\tilde{\vec{a}}(t) = \sum_{j=1}^k\beta_j\vec{\varphi}_je^{\omega_j t},
	\label{a_tilde_exp}
\end{equation}
where $(\vec{\varphi}_i,\omega_i)$ are the eigenpairs of $\mat{B}$ and $\{\beta_i\}_{i=1}^k$ is some set of coefficients. In this case $\mathcal{A}_k$ is characterized by the set of eigenpairs and coefficients $\{(\vec{\varphi}_i,\omega_i,\beta_i)\}_{i=1}^k$.
The original function $\vec{a}(t)$ is then reconstructed through the map $\mathcal{M}^t$ such that
\begin{equation}
	\vec{a}(t)\approx \mathcal{M}^t\mathcal{A}_k = \sum_{j=1}^k\beta_j\vec{\varphi}_je^{\omega_j t}.
	\label{dmd_exp}
\end{equation}

\noindent To find the eigenpairs $(\vec{\varphi}_i,\omega_i)$, let us define the orbital data matrices
\begin{equation}
	\mat{X} = [ \vec{a}_0 \ \vec{a}_1 \, \dots \, \vec{a}_{m-1} ] \in \real^{n\times m},\quad \hat{\mat{X}} = [ \vec{a}_1 \ \vec{a}_2 \, \dots \, \vec{a}_{m} ] \in \real^{n\times m},
\end{equation}

\noindent then $\tilde{\mat{B}}=\hat{\mat{X}}\mat{X}^+$ is the closest approximation to $\mat{B}$ in the Frobenius norm where $+$ signifies the Moore-Penrose pseudo inverse \cite{ipsen-2009}. The eigenpairs of $\tilde{\mat{B}}$, written as $(\tilde{\vec{\varphi}}_i,\lambda_i)$, are closely related to the eigenpairs of $\mat{B}$ \cite{tu-rowley-2014}
and each eigenvector $\vec{\varphi}_{i}$ can be calculated from the corresponding reduced eigenvector $\tilde{\vec{\varphi}}_{i}$.
The pairs $(\vec{\varphi}_i,\lambda_i)$
can construct $\tilde{\vec{a}}(t)$ at the specific points $\{t^n\}_{n=0}^m$ as
\begin{equation}
	\tilde{\vec{a}}(t^n) = \sum_{j=1}^k\beta_j\vec{\varphi}_j\lambda_j^n,\quad  n=0,\dots,m.
\end{equation}

\noindent This expression yields the DMD expansion \eqref{dmd_exp} with the transformation $\omega_j = \frac{\ln(\lambda_{j})}{\Delta t}$. The pairs $(\vec{\varphi}_i,\lambda_i)$ are called DMD modes and eigenvalues and are in practice calculated via the projected linear operator $\tilde{\mat{B}}_k=\mat{U}_k^\top\tilde{\mat{B}}\mat{U}_k$, whose eigenpairs
are written as $(\tilde{\vec{\varphi}}_i^{(k)},\lambda_i)$.
Here $\mat{U}_k$ holds the left singular vectors of $\mat{X}$ in its columns.
Note that the eigenvalues of $\tilde{\mat{B}}_k$ are the DMD eigenvalues.
The process of calculating the eigenpairs $(\vec{\varphi}_i,\lambda_i)$ is outlined in \linebreak Algorithm \ref{alg:dmd_tsvd} \cite{tu-rowley-2014}.

\begin{algorithm}[ht!]
	\SetAlgoLined
	\vspace*{.1cm}
	Input: solution data $\{\vec{a}^n\}_{n=0}^{m}$ and $\xi_{\text{rel}}$\\
	\begin{enumerate}
		\item Construct data matrices $\mat{X}, \hat{\mat{X}} \leftarrow \{{\vec{a}}^n\}_{n=0}^{m}$
		\item Compute truncated SVD $\mat{X} \approx \mat{U}_k\mat{S}_k\mat{V}_k^\top$ with $k\leq \text{rank}(\mat{X})$ satisfying equation \eqref{PODerr_rel} given $\xi_{\text{rel}}$
		\item Compute reduced DMD matrix $\tilde{\mat{B}}_k = \mat{U}_k^\top \hat{\mat{X}}\mat{V}_k\mat{S}_k^{-1}$
		\item Find eigenpairs $\{(\tilde{\vec{\varphi}}_{i}^{(k)},\lambda_{i})\}_{i=1}^{k}$ of $\tilde{\mat{B}}_k$
		\item Compute DMD modes:
		\begin{itemize}
			\item (Exact DMD) $\vec{\varphi}_{i} \leftarrow \frac{1}{\lambda_{i}}\hat{\mat{X}}\mat{V}_k\mat{S}_k^{-1}\tilde{\vec{\varphi}}_{i}^{(k)}, \ \ \lambda_{i}\ne 0, \ \ i=1,\dots,k$
			\item (Projected DMD) $\hat{\vec{\varphi}}_{i} \leftarrow \mat{U}_k\tilde{\vec{\varphi}}_{i}^{(k)}, \ \ i=1,\dots,k$
		\end{itemize}
		
	\end{enumerate}
	\vspace*{-.2cm}
	Output: DMD modes $\{\vec{\varphi}_{i}\}_{i=1}^{k}$ or $\{\hat{\vec{\varphi}}_{i}\}_{i=1}^{k}$ and DMD eigenvalues $\{\lambda_{i}\}_{i=1}^{k}$
	\vspace*{.2cm}
	\caption{Algorithm for computing DMD modes and eigenvalues  \cite{tu-rowley-2014}\label{alg:dmd_tsvd}}
\end{algorithm}

In Algorithm \ref{alg:dmd_tsvd} there are two types of DMD modes that can be calculated: (i) \textit{exact} DMD modes and (ii) \textit{projected} DMD modes. In practice the \textit{exact} DMD modes are preferred, as they can be shown to be the eigenvalues of the linear operator $\mat{B}$ that lie in the image of $\hat{\mat{X}}$. The \textit{projected} DMD modes have been shown to be simply the projection of the exact modes onto the image of $\mat{X}$ \cite{tu-rowley-2014}. Because the \textit{exact} DMD modes are generally regarded as the default in literature we find it important to note that for the ROMs developed in this paper, when the DMD is invoked we actually use the \textit{projected} modes.

The projected DMD modes were used in the original formulation of the DMD, which can be interpreted as a method that approximates the last data-vector as a linear combination of all former vectors, i.e.
\begin{equation}
	\vec{a}_m = \sum_{i=0}^{m-1}c_i\vec{a}_i + r,
\end{equation}

\noindent where $c_i$ are coefficients and $r$ is the residual incurred by the DMD approximation \cite{chen-tu-2012}. It comes naturally then, that when the exact DMD modes are used instead of the projected DMD modes, the DMD can be interpreted as approximating the \textit{first} data-vector as a linear combination of all latter vectors,
\begin{equation*}
	\vec{a}_0 = \sum_{i=1}^{m}c_i\vec{a}_i + r,
\end{equation*}
since the exact DMD modes lie in the image of $\hat{\mat{X}}$. This effective `shift' of the DMD residual to the first data-vector can come at a large cost to the time-dependent problems we consider, where the initial transients tend to be more difficult to capture compared to later times. It is with this in mind that we choose to utilize the projected DMD modes in this paper when applying the expansion \eqref{dmd_exp}.

%=================================================================================
\subsection{Equilibrium-Subtracted DMD}

In this paper a variant of the DMD is also considered which we will refer to as the equilibrium-subtracted DMD, or DMD-E. The DMD-E differs from the DMD by constructing the  linear operator $\mat{B}$ to fit the equilibrium-subtracted data $\{\check{\vec{a}}^n\}_{n=0}^{m'}$,
where $\check{\vec{a}}^n=\vec{a}^n-\vec{a}_b$ and $\vec{a}_b$ is
the
equilibrium solution vector \cite{chen-tu-2012,alla-kutz-2017}. The same Algorithm \ref{alg:dmd_tsvd} is used to calculate the DMD-E eigenvectors and modes, only replacing $\{{\vec{a}}^n\}_{n=0}^m$ with $\{\check{\vec{a}}^n\}_{n=0}^{m'}$. Thus for the DMD-E $\mathcal{A}_k$ is characterized by the set of eigenpairs and coefficients $\{(\vec{\varphi}_i,\omega_i,\beta_i)\}_{i=1}^k$, along with the vector $\vec{a}_b$.
The original function $\vec{a}(t)$ is reconstructed through the map $\mathcal{M}^t$, similarly to Eq. \eqref{dmd_exp}, as
\vspace{-.4cm}
\begin{equation}
	\vec{a}(t)\approx \mathcal{M}^t\mathcal{A}_k = \vec{a}_b + \sum_{j=1}^k\beta_j\vec{\varphi}_je^{\omega_j t}.
	\label{dmdb_exp}
\end{equation}

The vector $\vec{a}_b$ is chosen from any equilibrium solution of the underlying system that determines $\vec{a}(t)$ \cite{alla-kutz-2017}. The time-dependent TRT problems under consideration here possess a steady-state solution that is approached as $t\rightarrow\infty$. The most natural choice for this application is then to let $\vec{a}_b=\lim\limits_{t\rightarrow\infty} \vec{a}(t)$. In
 this study  we use $\vec{a}_b=\vec{a}(t^m)$ to approximate the steady-state solution, so that $m'=m-1$ and the equilibrium subtracted data is
\begin{equation}
	\check{\vec{a}}^n = \vec{a}^n - \vec{a}^m,\quad  n=0,\dots,m-1.
\end{equation}

%=================================================================================
%
%=================================================================================
\section{Numerical Results} \label{sec:results}

\subsection{Test Problem}

To analyze the accuracy of the DDET ROMs, we use a 2-dimensional extension of the well-known Fleck-Cummings (F-C) test problem \cite{fleck-1971}. This F-C test takes the form of a square homogeneous domain in the $x-y$ plane, 6 cm in length on both sides. The domain is initially at a temperature of $T^0 = 1\ \text{eV}$, the left boundary of the domain is subject to incoming radiation with blackbody spectrum at a temperature of $T^\text{in}=1\ \text{KeV}$, and there is no  incoming radiation at  other boundaries. The material is characterized by an opacity of
\begin{equation}
	\varkappa_\nu = \frac{27}{\nu^3}\bigg(1-e^{-\frac{\nu}{T}}\bigg),
\end{equation}

\noindent and a material energy density that is linear in temperature
\vspace{-.2cm}
\begin{equation}
	\varepsilon(T) = c_vT,
	\vspace*{-.2cm}
\end{equation}

\noindent
with material heat capacity $c_v = 0.5917 a_R (T^\text{in})^3$.

A uniform grid of  $20\times 20$ cells (i.e. $N_x=N_y=20$) with side lengths of $\Delta x = \Delta y = 0.3\ \text{cm}$  is used to discretize the slab. $N_g=17$ frequency groups are defined as shown in Table \ref{tab:freq_grps}.
The Abu-Shumays angular quadrature set q461214 with
36 discrete directions per quadrant  is used  \cite{abu-shumays-2004}. The total number of angular directions is $N_\Omega=144$.
The F-C problem is solved for the time interval
$0\leq t\leq 6 \ \text{ns}$ with $N_t=300$ uniform time steps $\Delta t = 2\times 10^{-2} \ \text{ns}$. When generating ROM solutions to the F-C problem, the following convergence criteria are used
(ref. Algorithm \ref{alg:rom_alg}):
$\epsilon_1=10^{-14}$ and  $\epsilon_2=10^{-15}$.

\begin{table}[ht!]
	\centering
	\caption{Upper boundaries for each frequency group}
	\begin{tabular}{|l|l|l|l|l|l|l|l|l|l|}
		\hline
		$g$ &1&2&3&4&5&6&7&8&9 \\ \hline
		$\nu_{g}$ [KeV]
		& 0.7075
		& 1.415
		& 2.123
		& 2.830
		& 3.538
		& 4.245
		& 5.129
		& 6.014
		& 6.898 \\ \hline\hline
		$g$&10&11&12&13&14&15&16&17& \\ \hline
		$\nu_{g}$ [KeV]
		& 7.783
		& 8.667
		& 9.551
		& 10.44
		& 11.32
		& 12.20
		& 13.09
		& 1$\times 10^{7}$ & \\ \hline
	\end{tabular}
	\label{tab:freq_grps}
\end{table}

\begin{figure}[ht!]
	\centering
	\begin{tabular}{c|c|c|c}
		& t=1ns & t=2ns & t=3ns \\ \hline &&& \\[-.3cm]
		%%%
		$T$ &
		\raisebox{-.5\height}{\includegraphics[trim=1cm 0cm 4cm .5cm,clip,height=3.5cm]{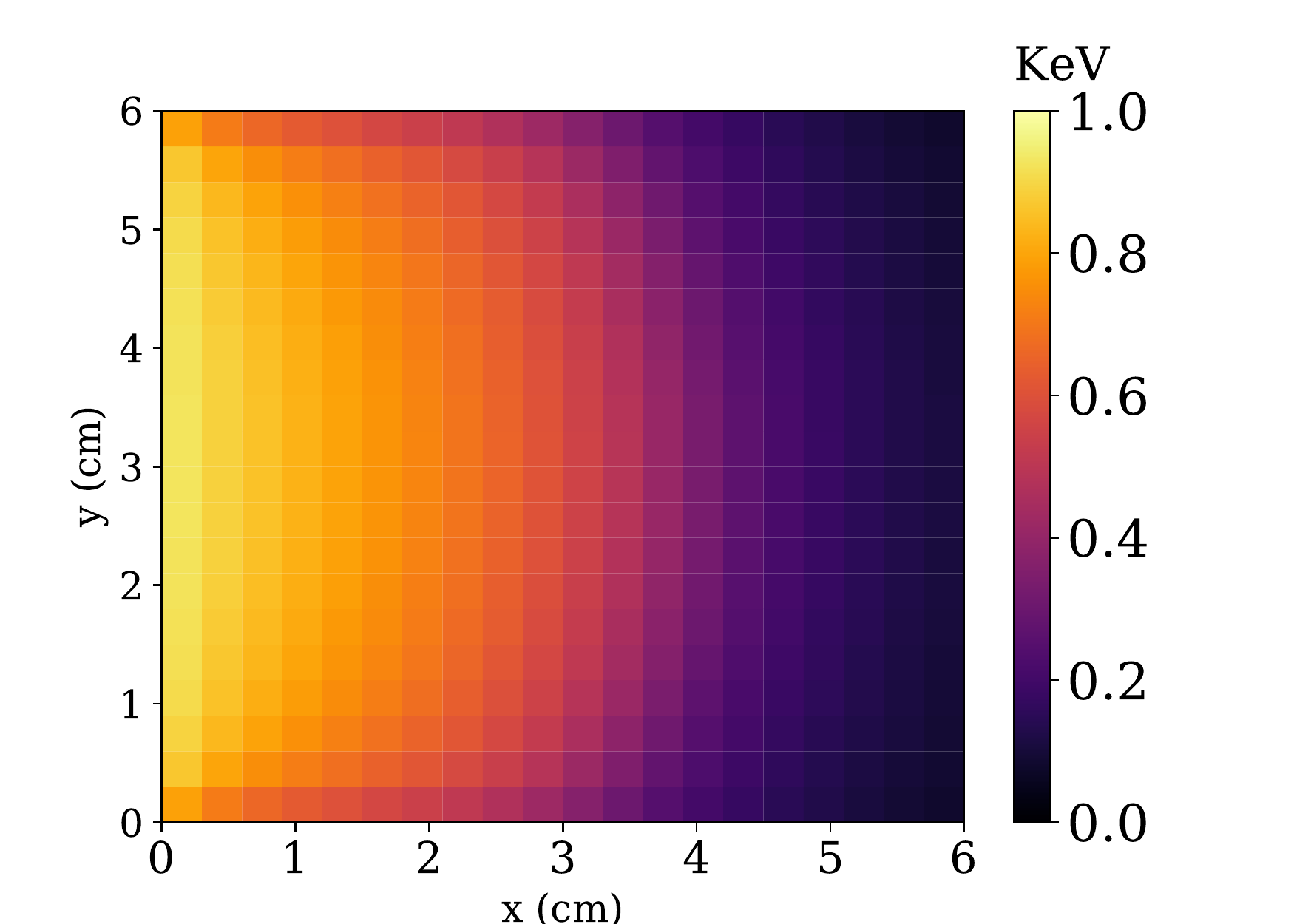}}&
		\raisebox{-.5\height}{\includegraphics[trim=1cm 0cm 4cm .5cm,clip,height=3.5cm]{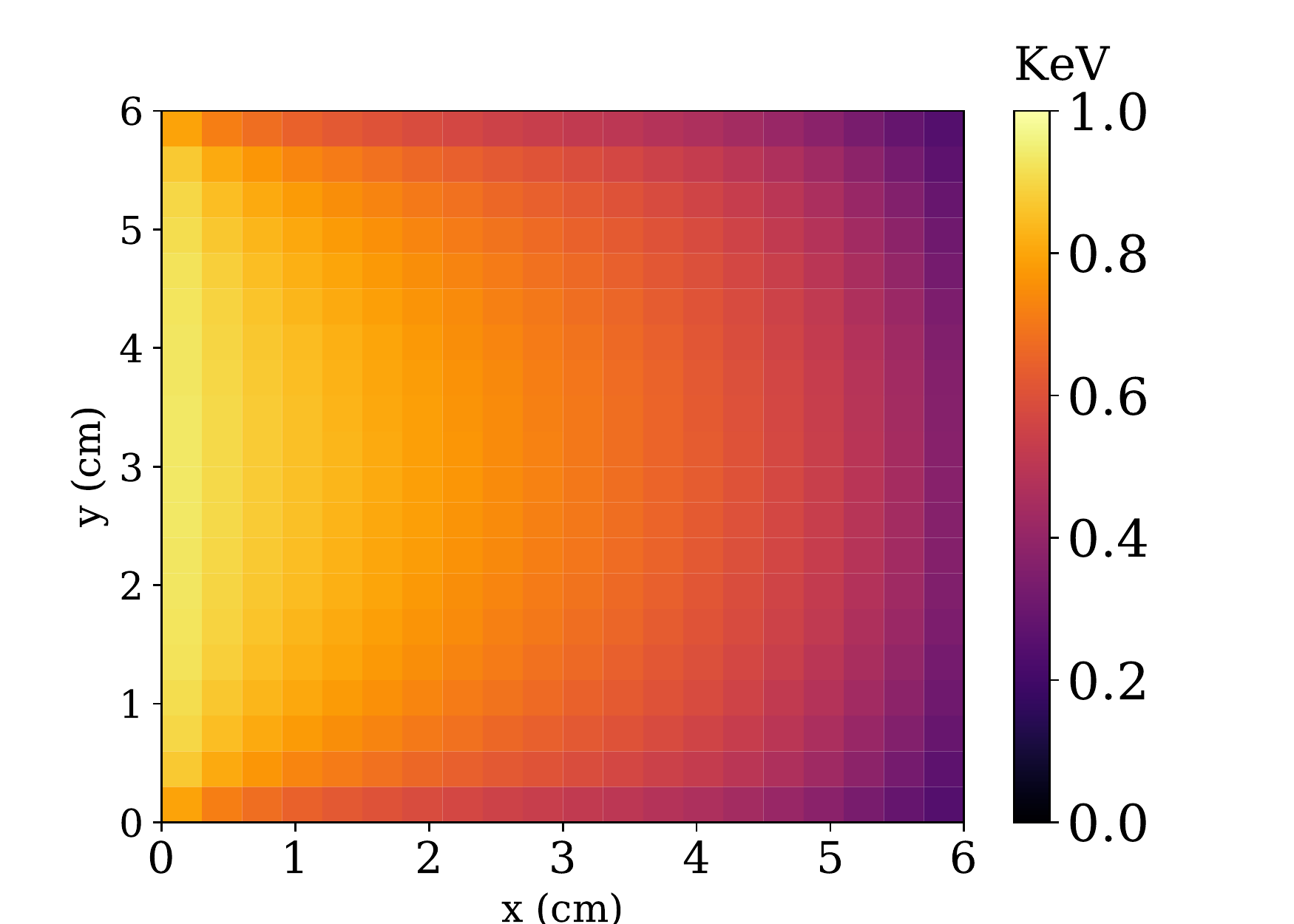}}&
		\raisebox{-.5\height}{\includegraphics[trim=1cm 0cm 0cm .5cm,clip,height=3.5cm]{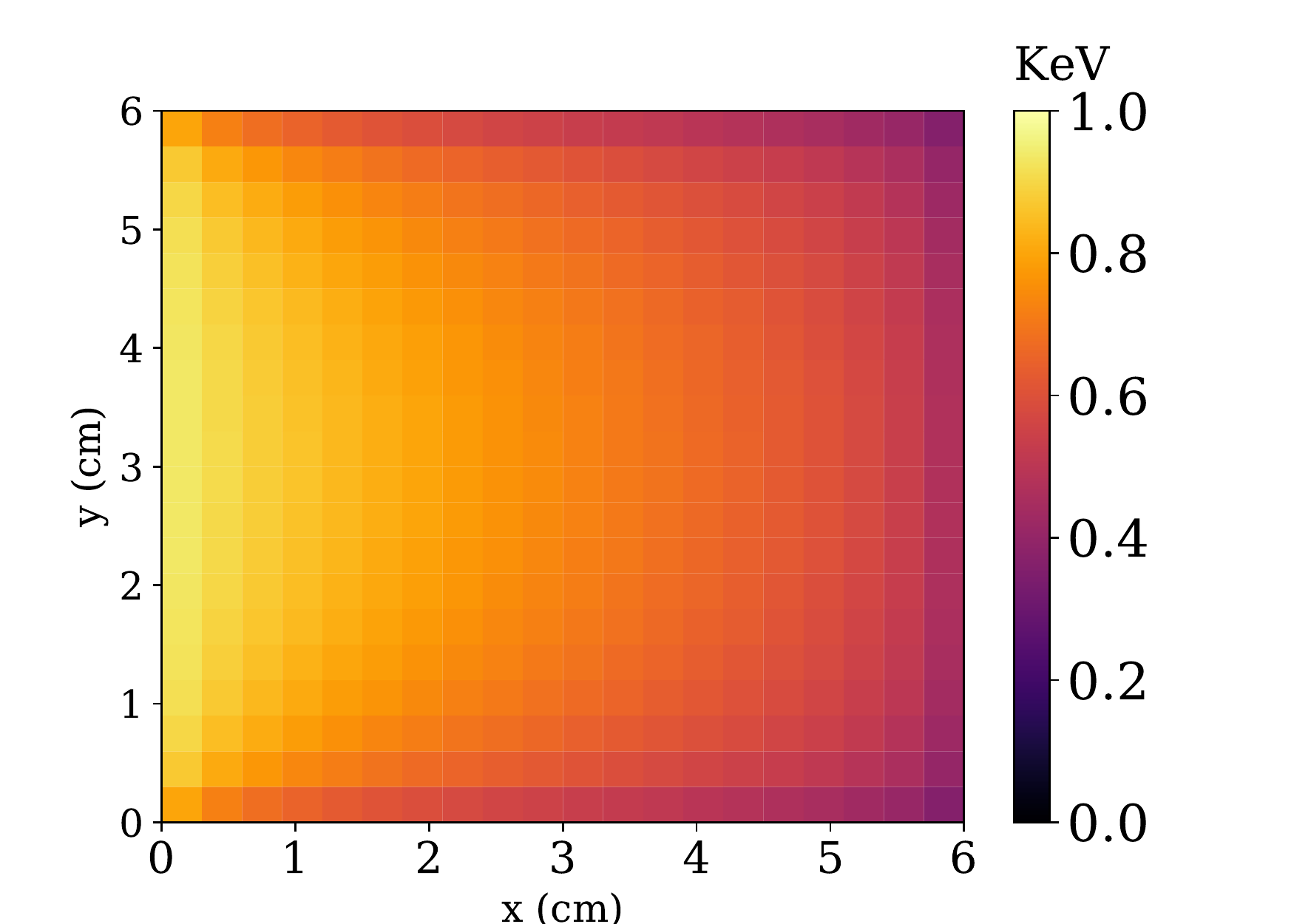}}\\[.1cm] \hline &&& \\[-.3cm]
		%%%
		$E$ &
		\raisebox{-.5\height}{\includegraphics[trim=1cm 0cm 4cm 0cm,clip,height=3.65cm]{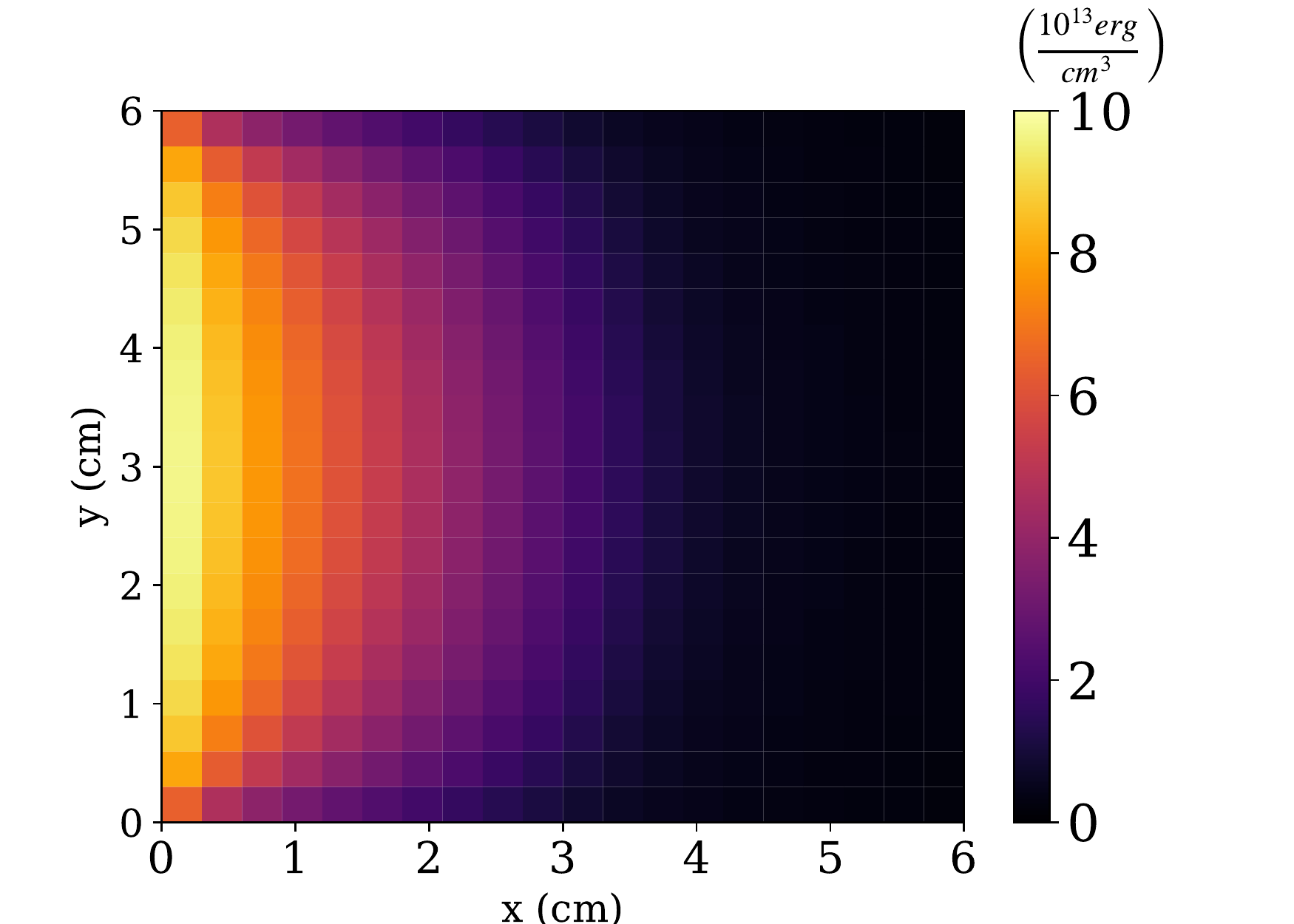}}&
		\raisebox{-.5\height}{\includegraphics[trim=1cm 0cm 4cm 0cm,clip,height=3.65cm]{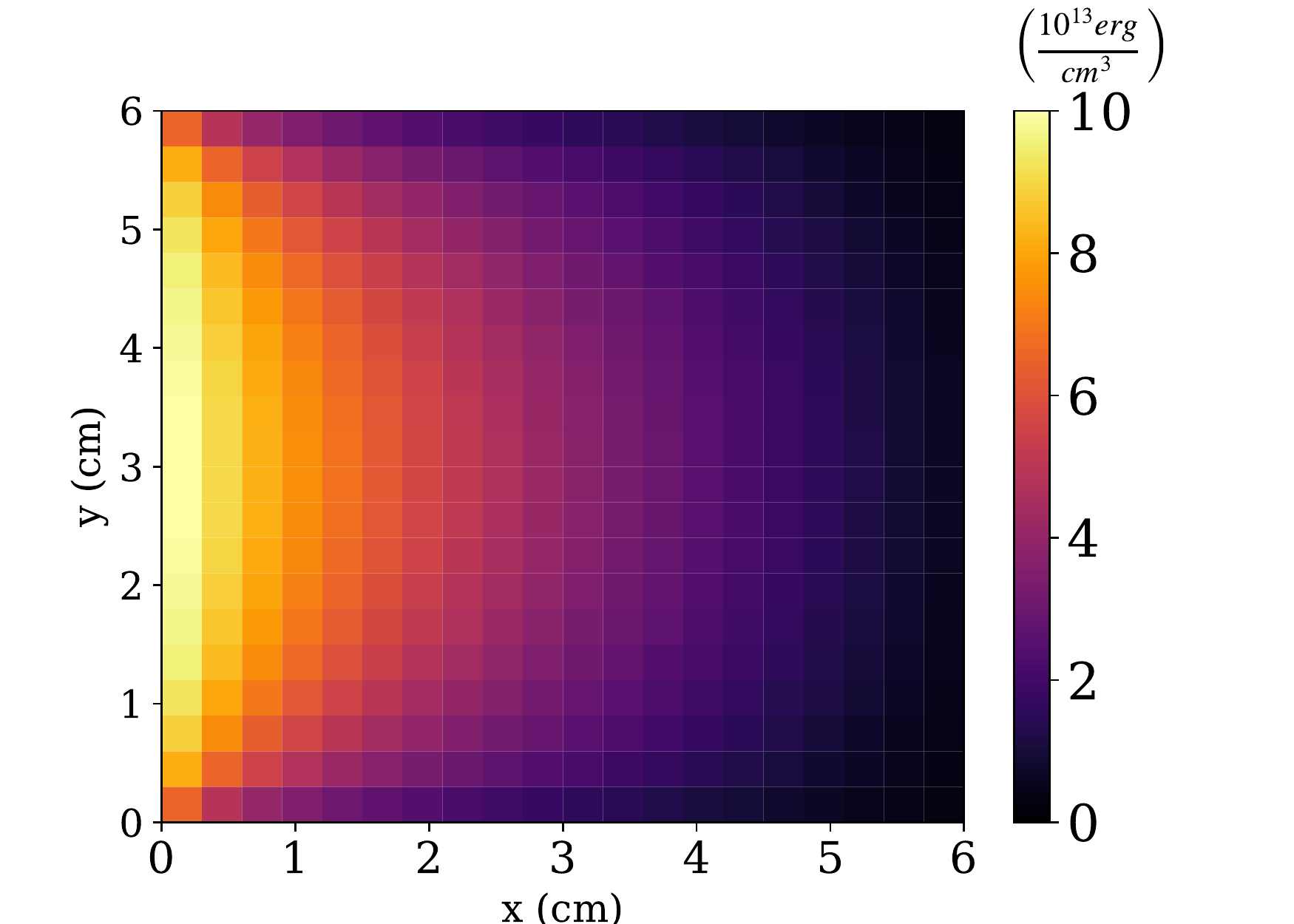}}&
		\raisebox{-.5\height}{\includegraphics[trim=1cm 0cm 0cm 0cm,clip,height=3.65cm]{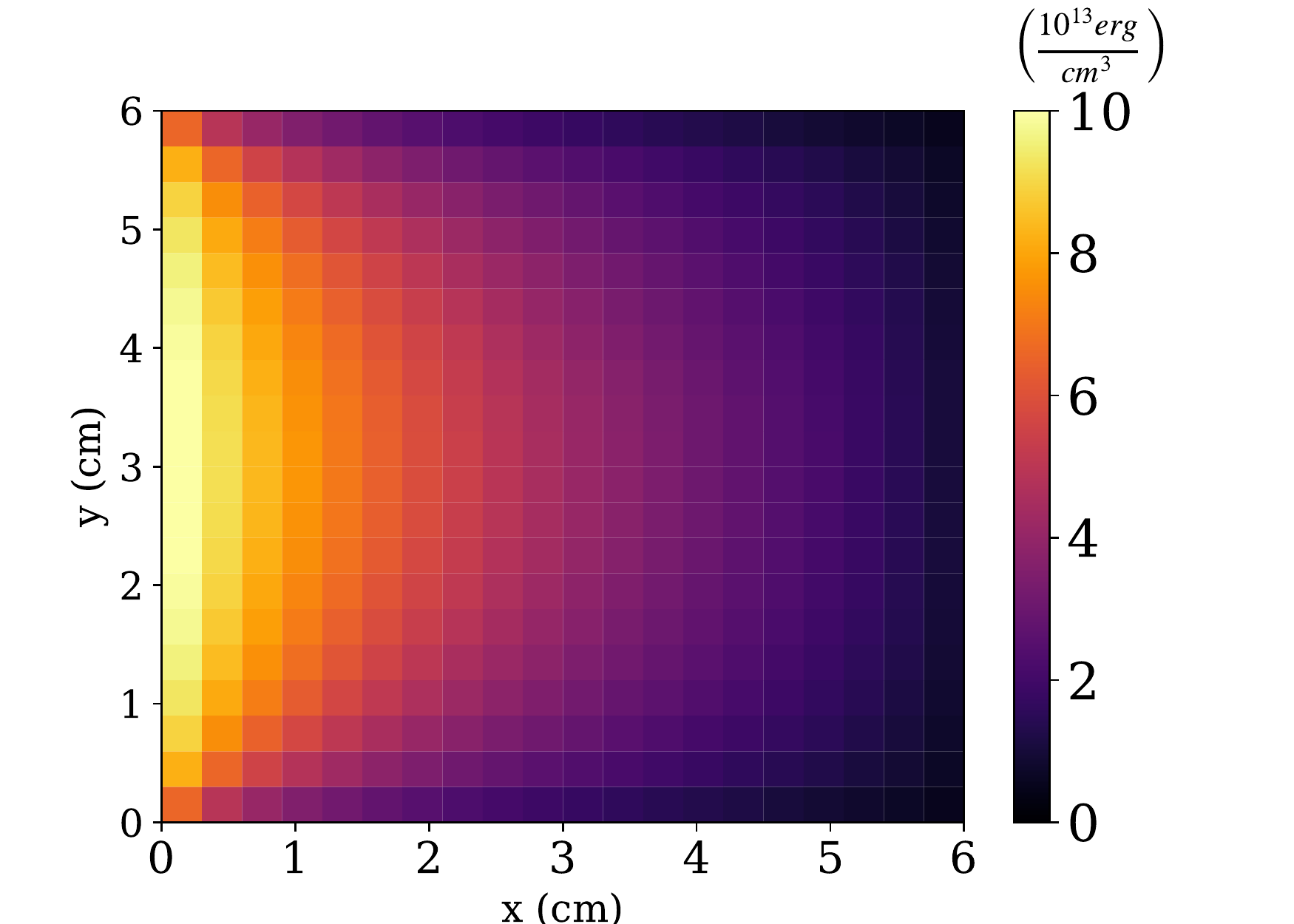}}
	\end{tabular}
	\vspace*{.25cm}
	\caption{F-C test solution for the material temperature $(T)$ and total radiation energy density $(E)$ over the spatial domain at times t=1,\,2,\,3\,ns.}
	
	\label{fig:fom_solution}
\end{figure}

The solution to this F-C test for the material temperature and total radiation energy density at times
$t=1,2,3$ ns is depicted in Figure \ref{fig:fom_solution}. The solutions of both $T$ and $E$ take the form of a wave that first rapidly forms on the left boundary before propagating to the right. After this the domain is continuously heated.
Eventually the solution reaches a regime close to steady state.

The total number of degrees of freedom occupied by the Eddington tensor $\boldsymbol{\mathfrak{f}}_g$ at a single instant of time is
$D_\tens{f} = 2(D_v+D_h+D_c)N_g = 4.216  \times 10^4$. In comparison, the degrees of freedom occupied by the radiation intensities from the simple corner balance scheme equals $D_I=4N_xN_yN_gN_\Omega=3.9168 \times 10^6$. This means that even before compressing the Eddington tensor with the POD or DMD, the required memory occupation is $\frac{D_I}{D_\tens{f}}=93$ times smaller than for the radiation intensities.

%=================================================================================
%
%=================================================================================
\subsection{Data Analysis}

\begin{figure}[ht!]
	\centering
	\vspace{-.5cm}
	\subfloat[$\mat{A}^{\tens{f}_{xx,c}}$]{\includegraphics[trim=2cm 1.5cm 2cm 2cm,clip,width=.4\textwidth]{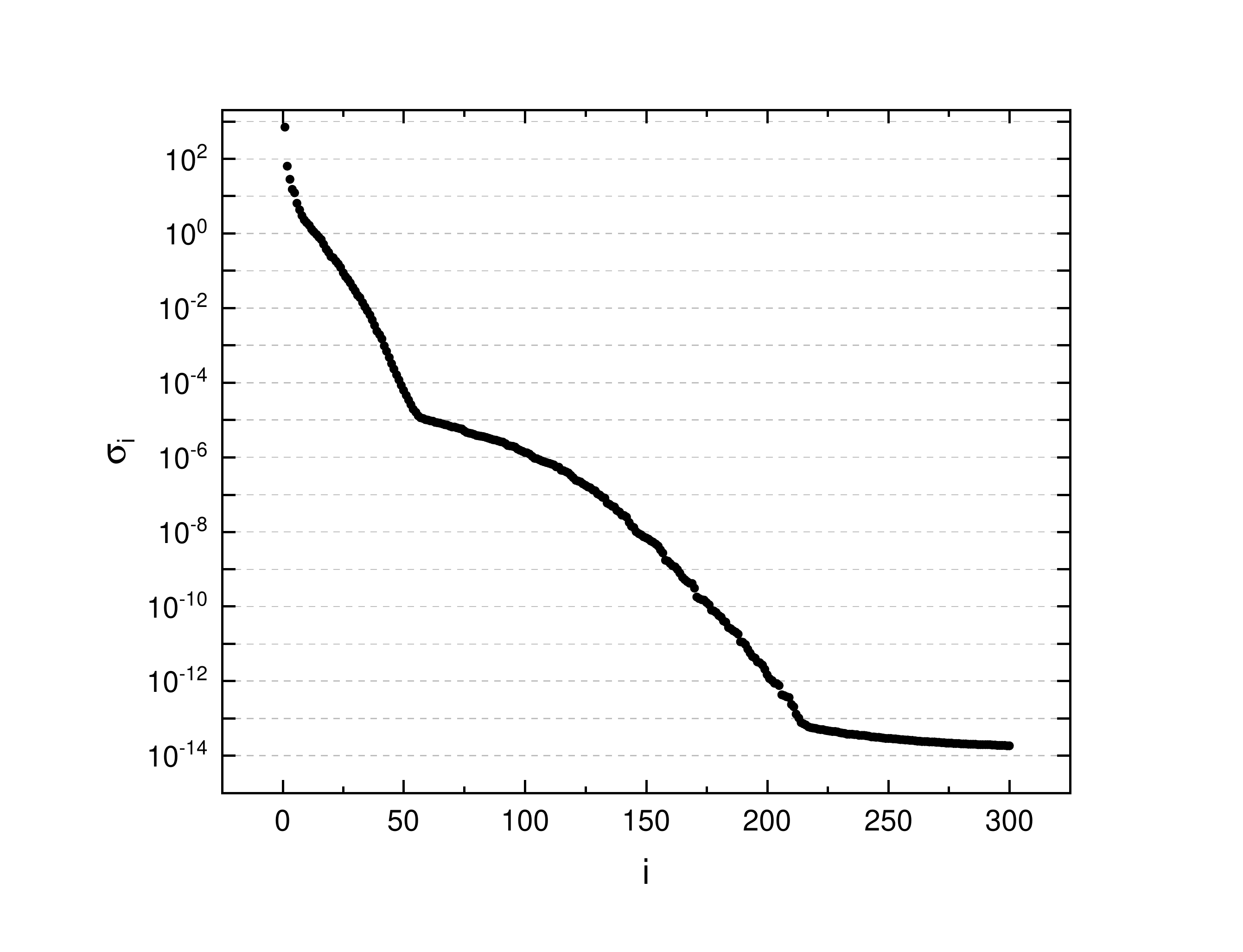}}
	\subfloat[$\mat{A}^{\tens{f}_{yy,h}}$]{\includegraphics[trim=2cm 1.5cm 2cm 2cm,clip,width=.4\textwidth]{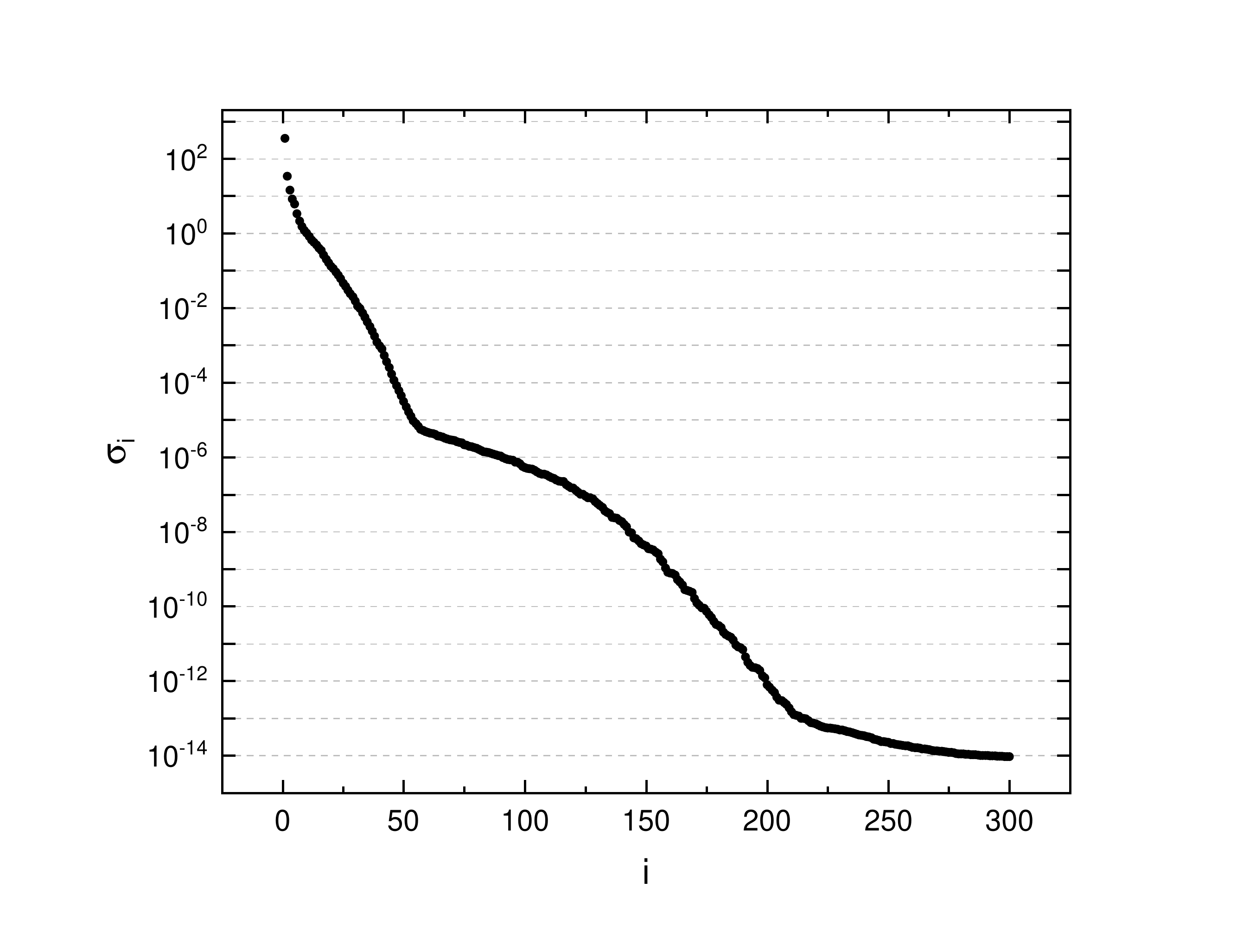}}\\
	\subfloat[$\mat{A}^{\tens{f}_{xy,v}}$]{\includegraphics[trim=2cm 1.5cm 2cm 2cm,clip,width=.4\textwidth]{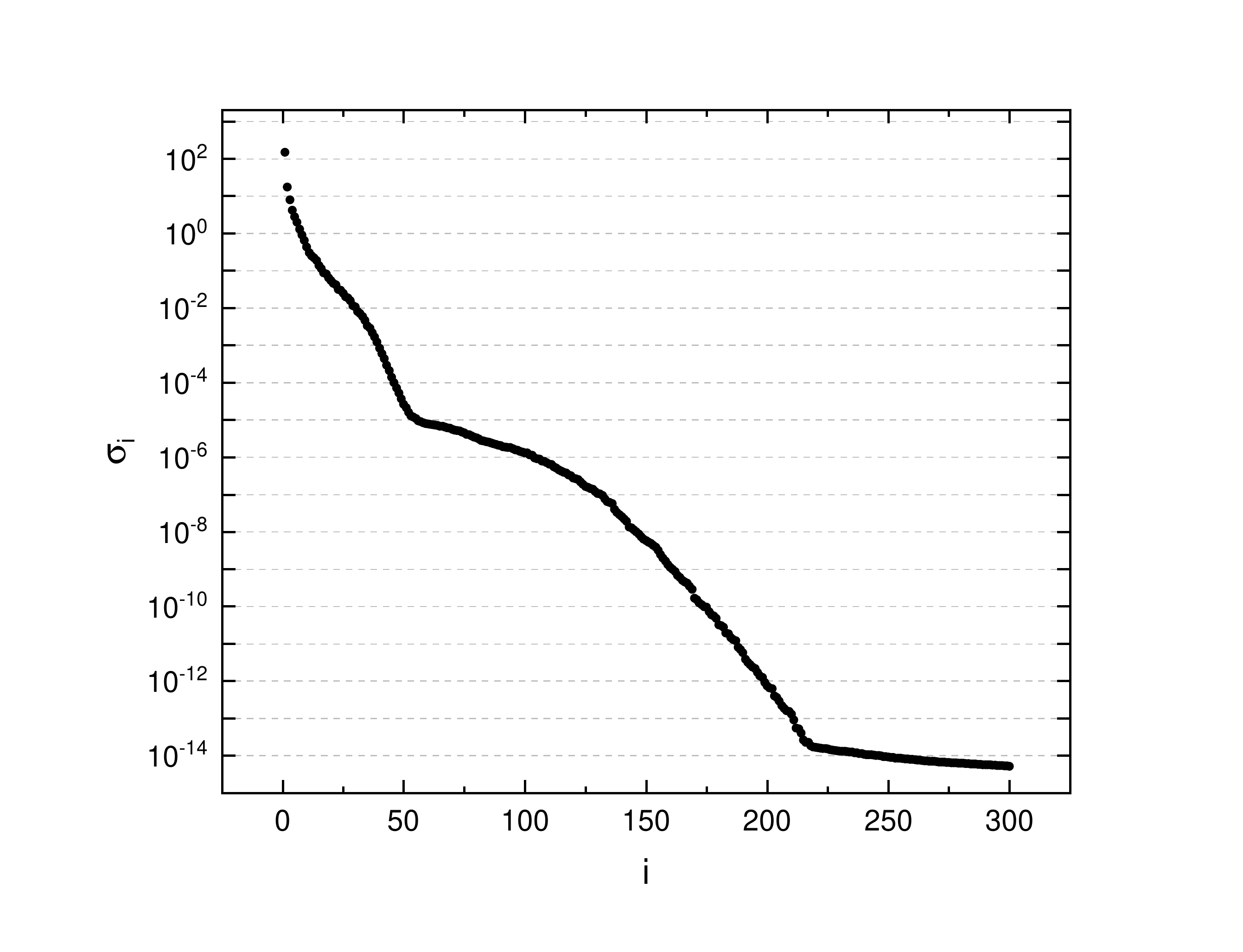}}
	\subfloat[$\mat{A}^c$]{\includegraphics[trim=2cm 1.5cm 2cm 2cm,clip,width=.4\textwidth]{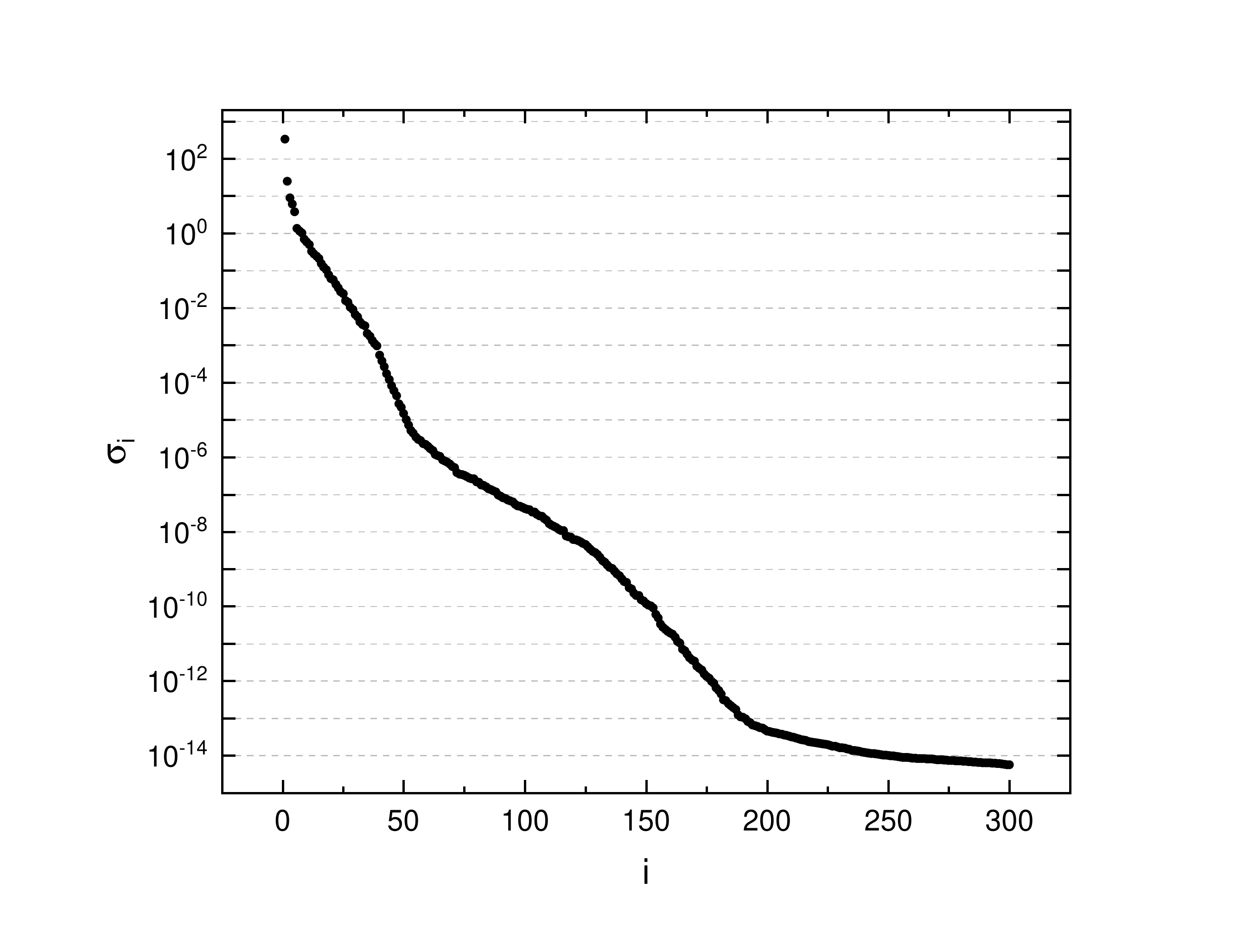}}
	
	\caption{Singular value distributions of select snapshot matrices
		of   grid functions of the Eddington tensor and boundary factor for
		the F-C test (included matrices are: (a) $\mat{A}^{\tens{f}_{xx,c}}$, (b) $\mat{A}^{\tens{f}_{yy,h}}$, (c) $\mat{A}^{\tens{f}_{xy,v}}$, (d) $\mat{A}^c$)}
	
	\label{fig:FC_svals}
\end{figure}

The snapshot data used to construct the matrices  $\mat{A}^{\mathfrak{f}_{\alpha\beta,\gamma}}$  and $\mat{A}^c$ (see Eq. \eqref{qd_snapmats})
is obtained by solving the TRT problem (Eqs. \eqref{bte_mg} \& \eqref{meb})  on the given grid in phase-space and time
by means of the MLQD method.
Convergence criteria for this simulation was set to $\epsilon=10^{-14}$.
This is the FOM solution of the test problem that is used as the reference solution.
The singular values of a select few of these snapshot matrices are depicted in Figure \ref{fig:FC_svals}. The singular values of those databases not shown here do not deviate significantly from the chosen plots.  The singular values for each of the
databases decay in a similar manner with 3 distinct sharp drops in magnitudes before reaching a value of approximately $10^{-14}$ where decay halts. The singular values that have a value at or below $10^{-14}$ have reached the limit of machine precision and can be considered numerically zero.

\newpage Although the POD, DMD and DMD-E make use of slight variations on these snapshot matrices, the singular value distributions of these variant matrices are very similar to those pictured. For the databases without their final column, used for the DMD, their SVD is almost exactly the same as for the full matrices since the final column holds near steady-state data and does not add much new information to the span of the columns. For the POD when the databases are centered about their column-mean, the only significant difference from the shown plots is in the first singular value which decreases by roughly an order of magnitude.  The second singular value is also seen to decrease by roughly half. The equilibrium-subtracted databases used for the DMD-E acquire singular value spectra very similar to those obtained through the POD.

\begin{figure}[ht!]
	%\vspace*{-.5cm}
	\centering
	\begin{minipage}{.6\textwidth}
		\captionof{table}{Ranks $k$ for each approximate database \\corresponding to different values of $\xi_{\text{rel}}$ for the POD \label{tab:pod_ranks}}
		\hspace{-.35cm}
		\resizebox{.95\columnwidth}{!}{
			\def\arraystretch{1.2}
		\begin{tabular}{|l|l|l|l|l|l|l|l|}
			\hline
			$\xi_{\text{rel}}$ & $\mat{A}^{\tens{f}_{xx,c}}$ & $\mat{A}^{\tens{f}_{xx,v}}$ & $\mat{A}^{\tens{f}_{yy,c}}$ & $\mat{A}^{\tens{f}_{yy,h}}$ & $\mat{A}^{\tens{f}_{xy,v}}$ & $\mat{A}^{\tens{f}_{xy,h}}$ & $\mat{A}^c$ \\ \hline
			$10^{-2}$ & 15 & 17 & 15 & 15 & 14 & 17 & 14 \\ \hline
			$10^{-4}$ & 34 & 36 & 34 & 34 & 37 & 36 & 35 \\ \hline
			$10^{-6}$ & 49 & 49 & 52 & 49 & 65 & 68 & 48 \\ \hline
			$10^{-8}$ & 115 & 110 & 120 & 115 & 129 & 127 & 87 \\ \hline
			$10^{-10}$ & 152 & 148 & 154 & 153 & 159 & 158 & 132 \\ \hline
			$10^{-12}$ & 179 & 178 & 180 & 180 & 185 & 184 & 160 \\ \hline
			$10^{-14}$ & 203 & 203 & 205 & 205 & 207 & 207 & 188 \\ \hline
			$10^{-16}$ & 300 & 300 & 300 & 300 & 300 & 300 & 300 \\ \hline
		\end{tabular}
	    }
	\end{minipage}%
	\begin{minipage}{.4\textwidth}
		\hspace*{-.5cm}
		\includegraphics[trim=0cm .5cm 0cm -1.5cm,clip,width=1.15\textwidth]{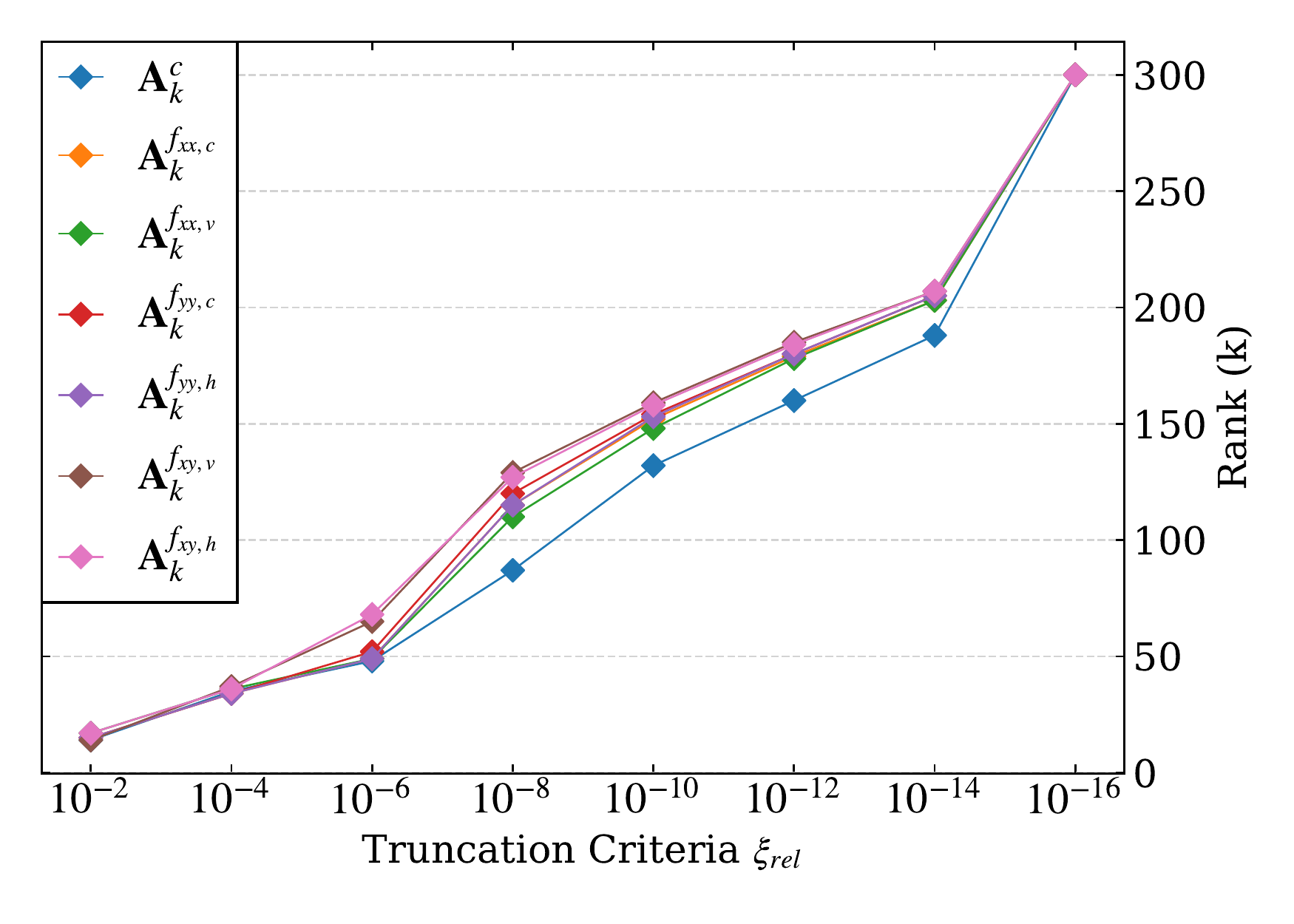}
		\captionof{figure}{Plotted ranks $k$ for the POD (see Table \ref{tab:pod_ranks}) \label{fig:pod_ranks}
		}
	\end{minipage}
\end{figure}

\begin{figure}[ht!]
	\centering
	\begin{minipage}{.6\textwidth}
		\captionof{table}{Ranks $k$ for each approximate database \\corresponding to different values of $\xi_{\text{rel}}$ for the DMD
			\label{tab:dmd_ranks}}
		\hspace{-.35cm}
		\resizebox{.95\columnwidth}{!}{
			\def\arraystretch{1.2}
		\begin{tabular}{|l|l|l|l|l|l|l|l|}
			\hline
			$\xi_{\text{rel}}$ & $\mat{A}^{\tens{f}_{xx,c}}$ & $\mat{A}^{\tens{f}_{xx,v}}$ & $\mat{A}^{\tens{f}_{yy,c}}$ & $\mat{A}^{\tens{f}_{yy,h}}$ & $\mat{A}^{\tens{f}_{xy,v}}$ & $\mat{A}^{\tens{f}_{xy,h}}$ & $\mat{A}^c$ \\ \hline
			$10^{-2}$ & 6 & 7 & 6 & 6 & 7 & 9 & 5 \\ \hline
			$10^{-4}$ & 28 & 30 & 28 & 28 & 30 & 30 & 25 \\ \hline
			$10^{-6}$ & 43 & 44 & 44 & 43 & 46 & 46 & 42 \\ \hline
			$10^{-8}$ & 90 & 79 & 100 & 87 & 111 & 111 & 61 \\ \hline
			$10^{-10}$ & 138 & 136 & 142 & 139 & 148 & 147 & 112 \\ \hline
			$10^{-12}$ & 168 & 165 & 170 & 169 & 175 & 175 & 147 \\ \hline
			$10^{-14}$ & 195 & 194 & 196 & 196 & 199 & 199 & 173 \\ \hline
			$10^{-16}$ & 286 & 286 & 287 & 287 & 292 & 291 & 274 \\ \hline
			$10^{-18}$ & 299 & 299 & 299 & 299 & 299 & 299 & 299 \\ \hline
		\end{tabular}
	    }
	\end{minipage}%
	\begin{minipage}{.4\textwidth}
		\hspace*{-.5cm}
		\includegraphics[trim=0cm .5cm 0cm -1.5cm,clip,width=1.15\textwidth]{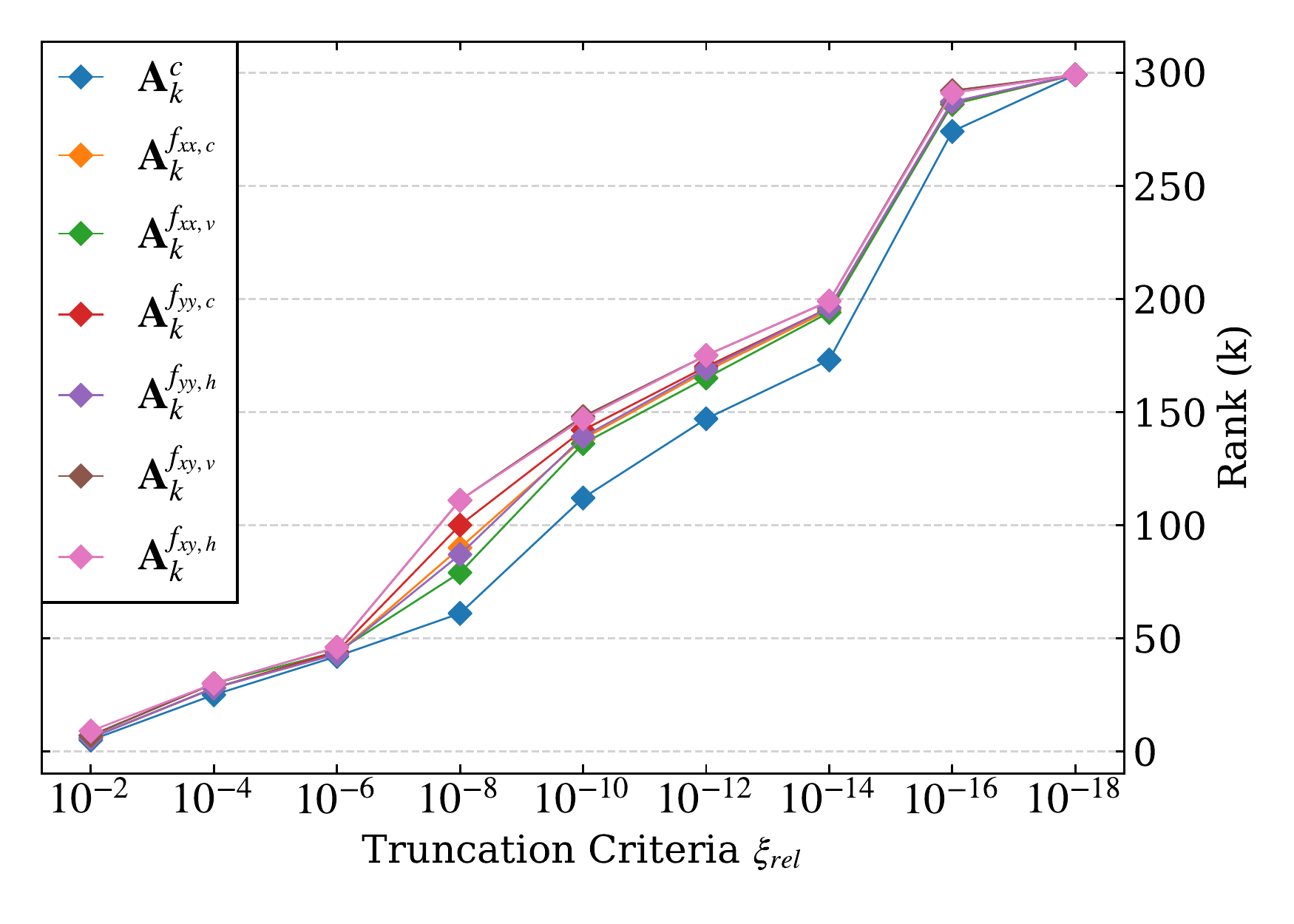}
		\captionof{figure}{Plotted ranks $k$ for the DMD (see Table \ref{tab:dmd_ranks}) \label{fig:dmd_ranks}
		}
	\end{minipage}
	\centering
	\begin{minipage}{.6\textwidth}
		\captionof{table}{Ranks $k$ for each approximate database \\corresponding to different values of $\xi_{\text{rel}}$ for the DMD-E
			\label{tab:dmdb_ranks}}
		\hspace{-.35cm}
		\resizebox{.95\columnwidth}{!}{
			\def\arraystretch{1.2}
		\begin{tabular}{|l|l|l|l|l|l|l|l|}
			\hline
			$\xi_{\text{rel}}$ & $\mat{A}^{\tens{f}_{xx,c}}$ & $\mat{A}^{\tens{f}_{xx,v}}$ & $\mat{A}^{\tens{f}_{yy,c}}$ & $\mat{A}^{\tens{f}_{yy,h}}$ & $\mat{A}^{\tens{f}_{xy,v}}$ & $\mat{A}^{\tens{f}_{xy,h}}$ & $\mat{A}^c$ \\ \hline
			$10^{-2}$ & 14 & 16 & 15 & 15 & 14 & 16 & 14 \\ \hline
			$10^{-4}$ & 34 & 36 & 34 & 34 & 36 & 35 & 34 \\ \hline
			$10^{-6}$ & 48 & 49 & 51 & 48 & 62 & 64 & 48 \\ \hline
			$10^{-8}$ & 114 & 109 & 119 & 114 & 125 & 127 & 85 \\ \hline
			$10^{-10}$ & 151 & 148 & 154 & 152 & 158 & 157 & 131 \\ \hline
			$10^{-12}$ & 179 & 177 & 180 & 179 & 184 & 183 & 160 \\ \hline
			$10^{-14}$ & 203 & 202 & 204 & 204 & 207 & 206 & 186 \\ \hline
			$10^{-16}$ & 298 & 298 & 298 & 298 & 298 & 298 & 298 \\ \hline
		\end{tabular}
		}
	\end{minipage}%
	\begin{minipage}{.4\textwidth}
		\hspace*{-.5cm}
		\includegraphics[trim=0cm .5cm 0cm -1.5cm,clip,width=1.15\textwidth]{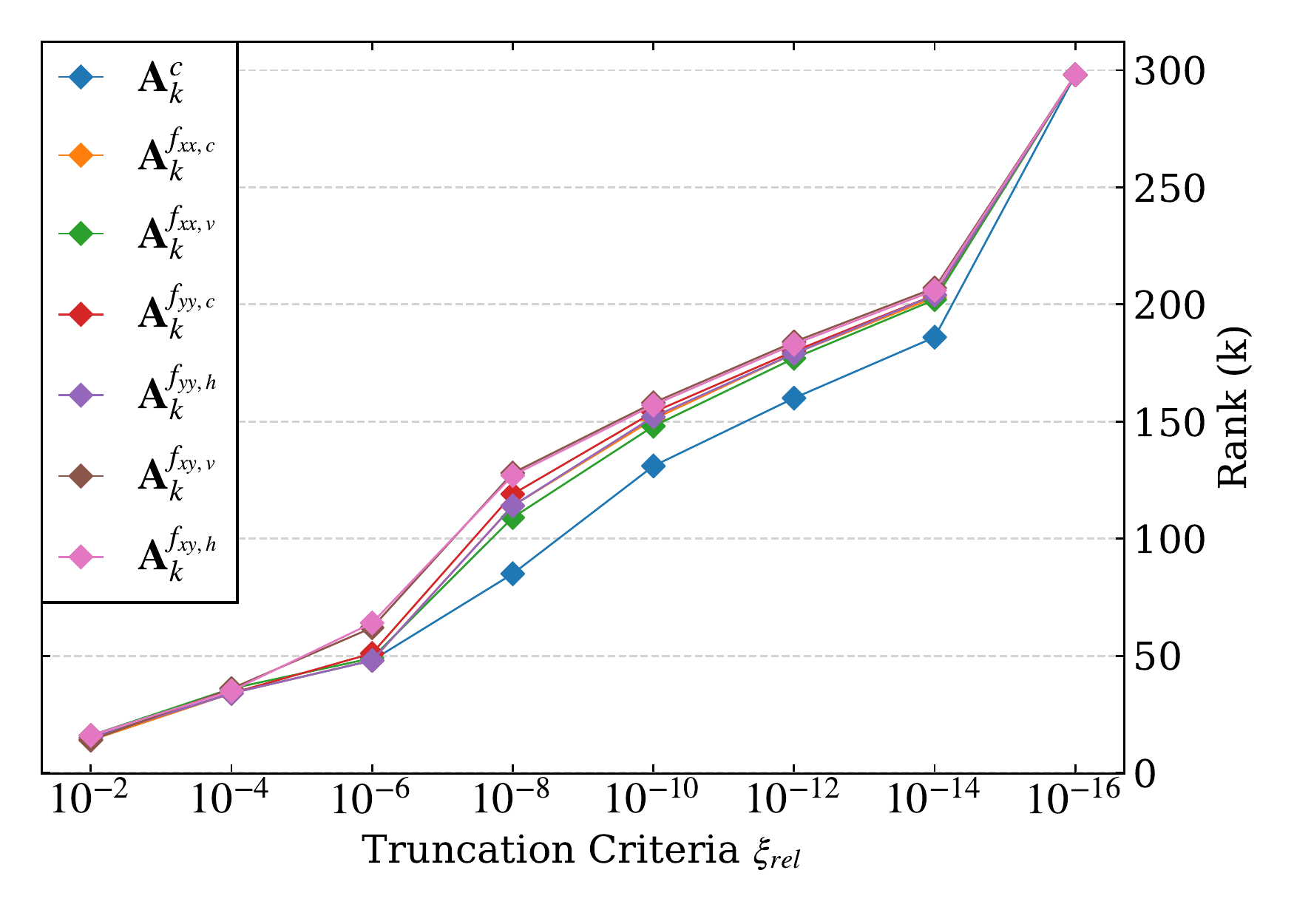}
		\captionof{figure}{Plotted ranks $k$ for the DMD-E (see Table \ref{tab:dmdb_ranks}) \label{fig:dmdb_ranks}
		}
	\end{minipage}
\vspace{-.5cm}
\end{figure}

The POD, DMD and DMD-E are applied to the databases $\mat{A}^{\tens{f}_{\alpha\beta,\gamma}}$ and $\mat{A}^c$ to generate several different rank-$k$ approximations of
the FOM Eddington tensor data, henceforth denoted as $\mat{A}^{\tens{f}_{\alpha\beta,\gamma}}_k$ and $\mat{A}^c_k$. For each method of approximation, several ranks $k$ were determined to satisfy a spectrum of chosen values for $\xi_{\text{rel}}$ while calculating the TSVD
(Eq. \eqref{PODerr_rel}). Tables \ref{tab:pod_ranks}, \ref{tab:dmd_ranks} and \ref{tab:dmdb_ranks} display the ranks used to approximate each individual database for every $\xi_{\text{rel}}$. Figures \ref{fig:pod_ranks}, \ref{fig:dmd_ranks} and \ref{fig:dmdb_ranks} plot these ranks against $\xi_{\text{rel}}$. The ranks used for the POD, DMD and DMD-E behave similarly with changes in $\xi_{\text{rel}}$ for each snapshot matrix, gradually increasing with decreases in $\xi_{\text{rel}}$ until $\xi_{\text{rel}}=10^{-16}$ where each database's rank
 increases by roughly 100.
This is due to the singular value decay structures shown in Figure \ref{fig:FC_svals} where decay stops after about
 200  singular values. The only significant difference in the used ranks between each of these methods given the same $\xi_{\text{rel}}$ is that the DMD always uses a lower rank than the POD and DMD-E. This is
an artifact of the centering and equilibrium-subtraction operations done on the databases prior to the calculation of each TSVD for the POD and DMD-E. Here these operations only significantly decreased the first  and second singular values  of each matrix. This has the effect of reducing only the denominator of equation \eqref{PODerr_rel} for all  $k>1$  and therefore inflating the rank required to satisfy a given $\xi_{\text{rel}}$.


\subsection{Performance of Low-Rank ROMs}

\begin{figure}[ht!]
	\vspace{.25cm}
	\centering
	\subfloat[Material Temperature]{\includegraphics[width=.5\textwidth]{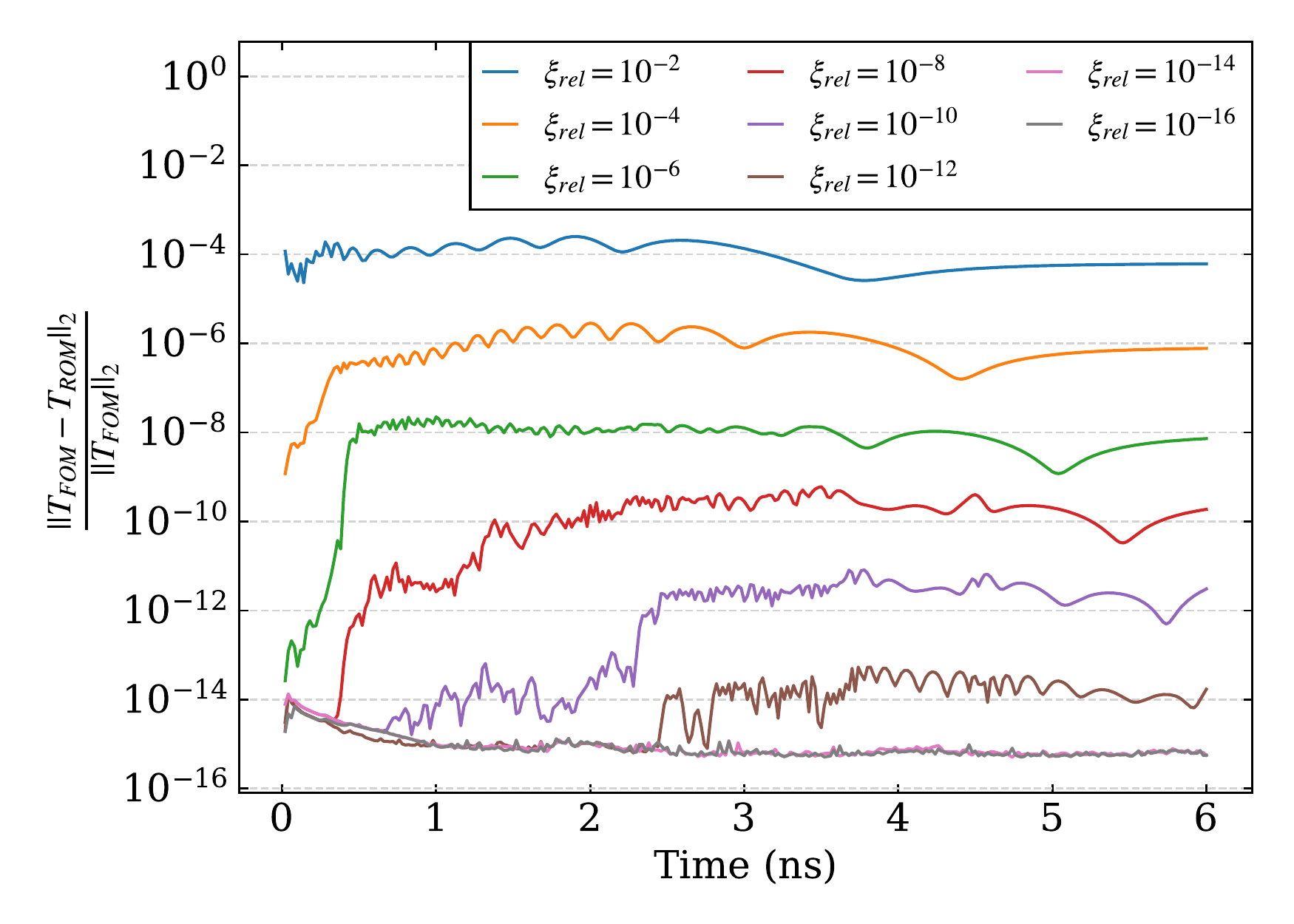}}
	\subfloat[Radiation Energy Density]{\includegraphics[width=.5\textwidth]{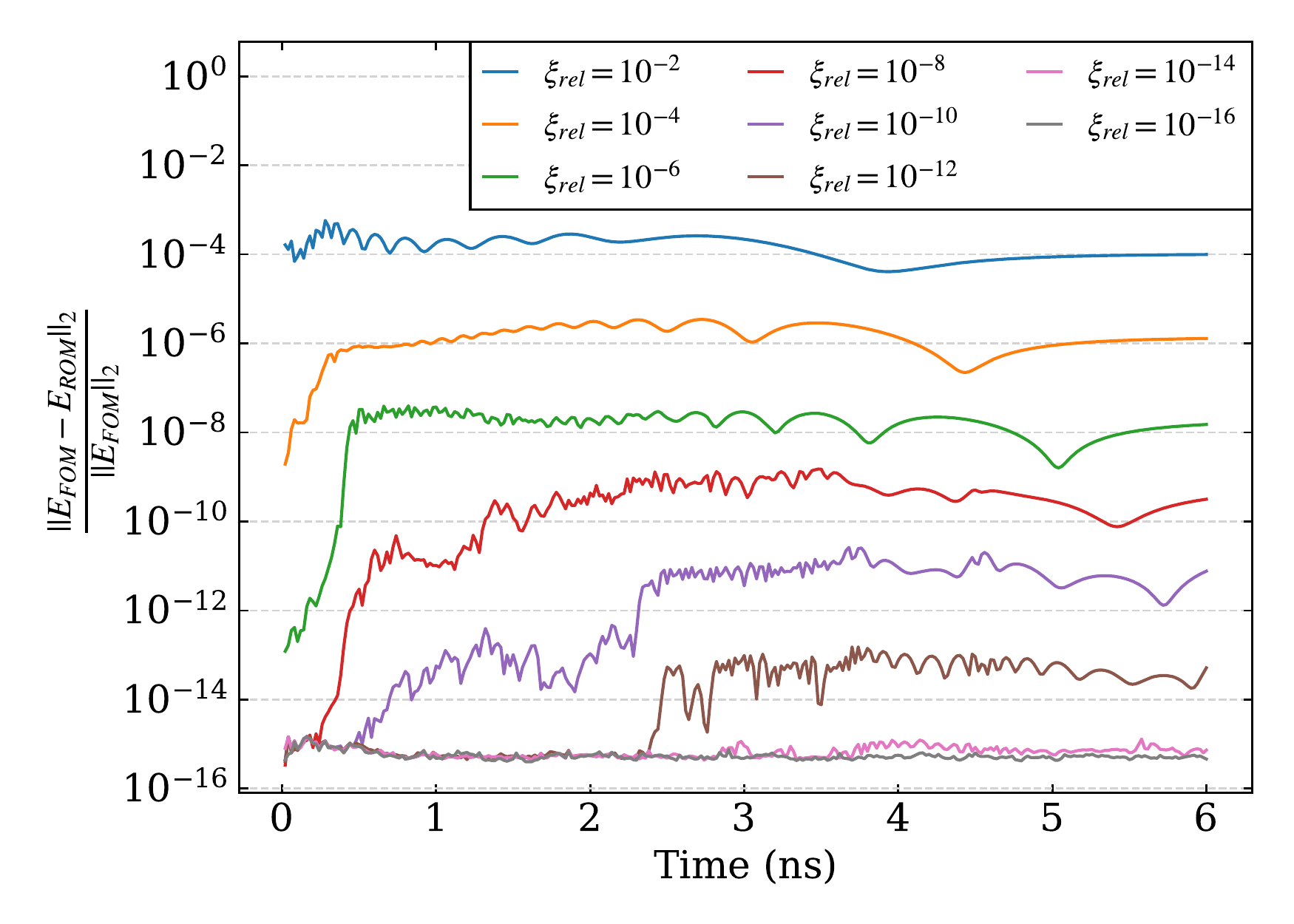}}
	\caption{Relative errors in the 2-norm of the DDET ROM using the POD using several $\xi_{\text{rel}}$, plotted vs time}
	\label{fig:2nrm-errs_POD}
	\vspace{.5cm}
	\centering
	\subfloat[Material Temperature]{\includegraphics[width=.5\textwidth]{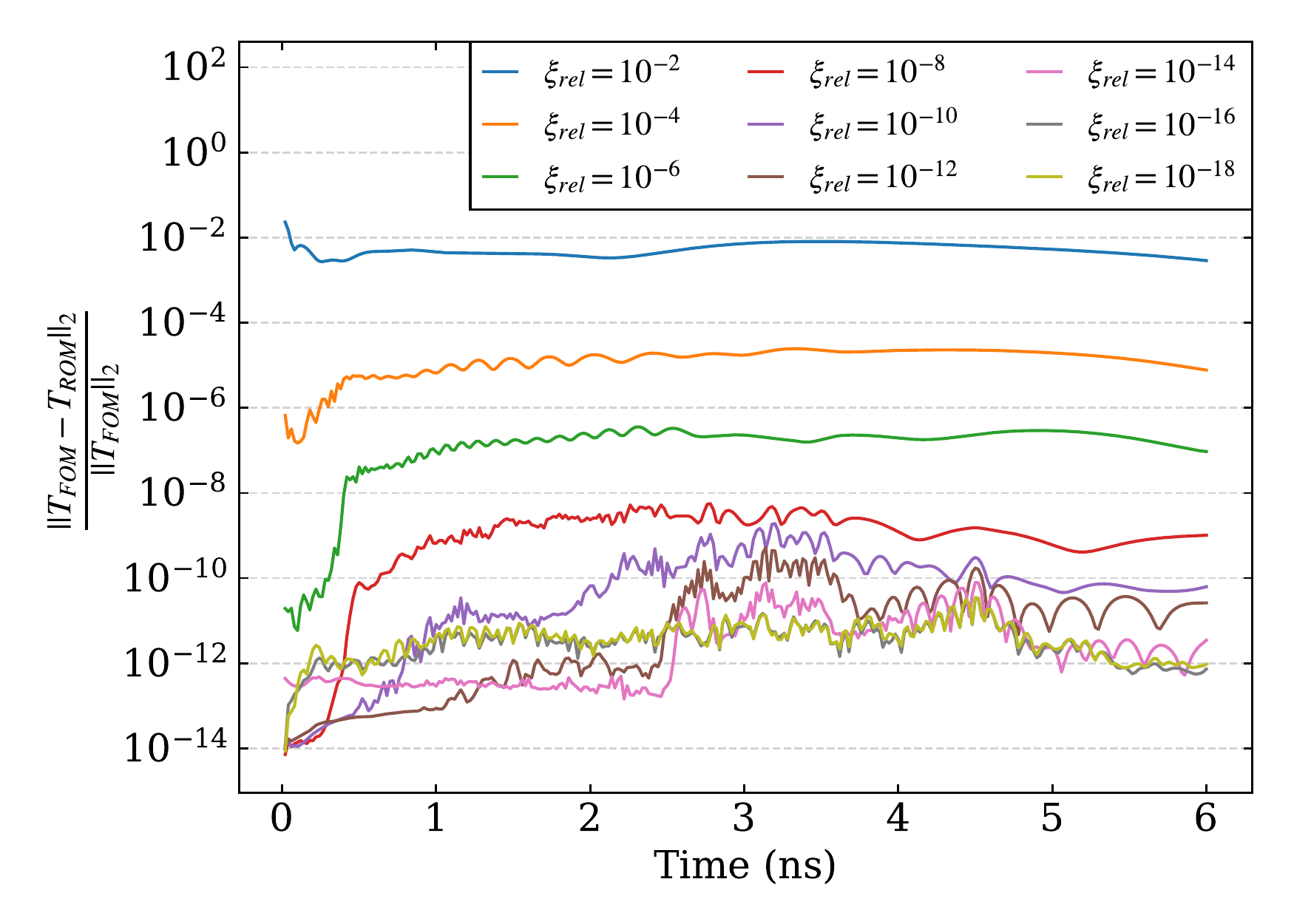}}
	\subfloat[Radiation Energy Density]{\includegraphics[width=.5\textwidth]{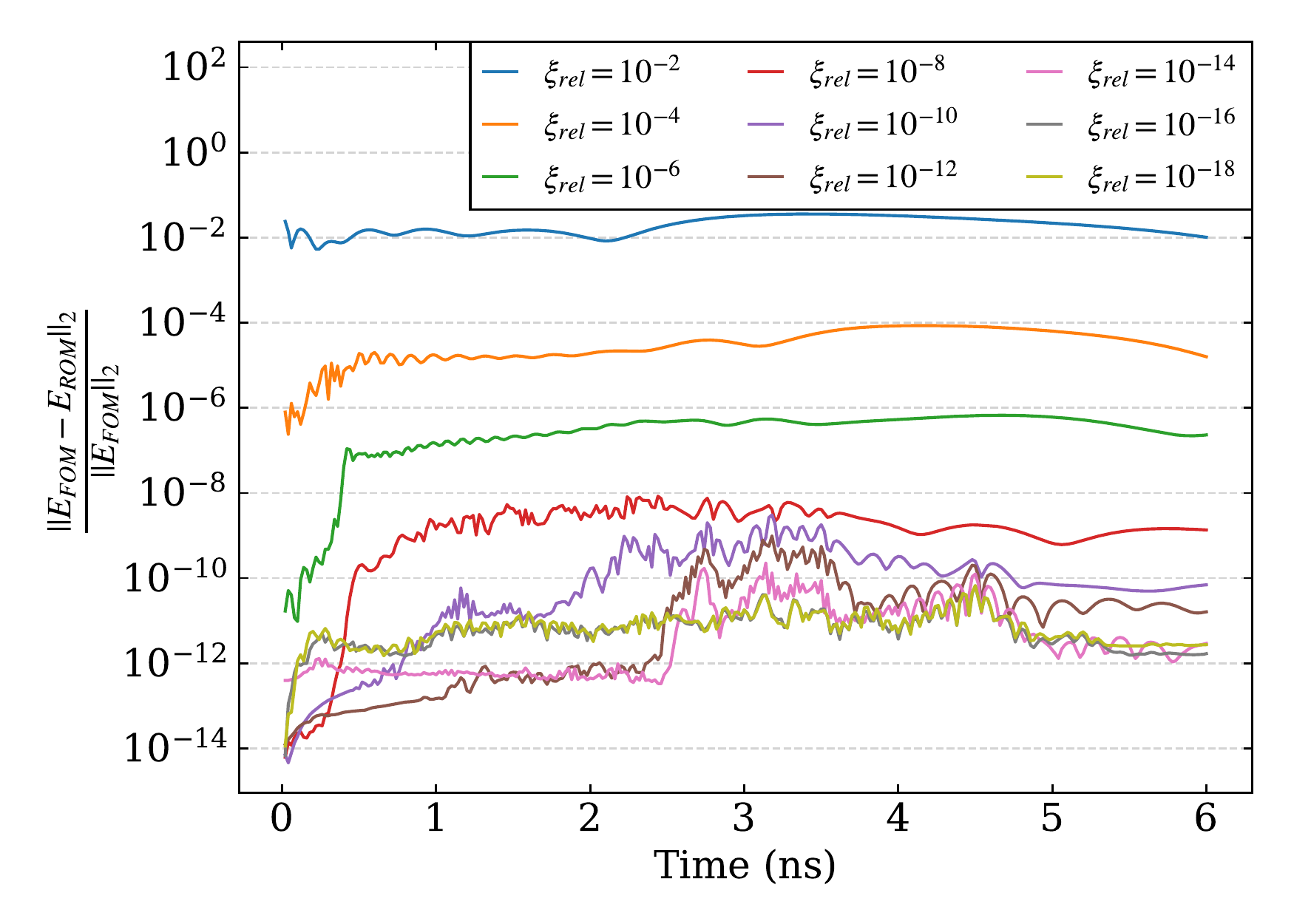}}
	\caption{Relative errors in the 2-norm of the DDET ROM using the DMD using several $\xi_{\text{rel}}$, plotted vs time}
	\label{fig:2nrm-errs_DMD}
	\vspace{.25cm}
\end{figure}

\begin{figure}[ht!]
	\vspace{.25cm}
	\centering
	\subfloat[Material Temperature]{\includegraphics[width=.5\textwidth]{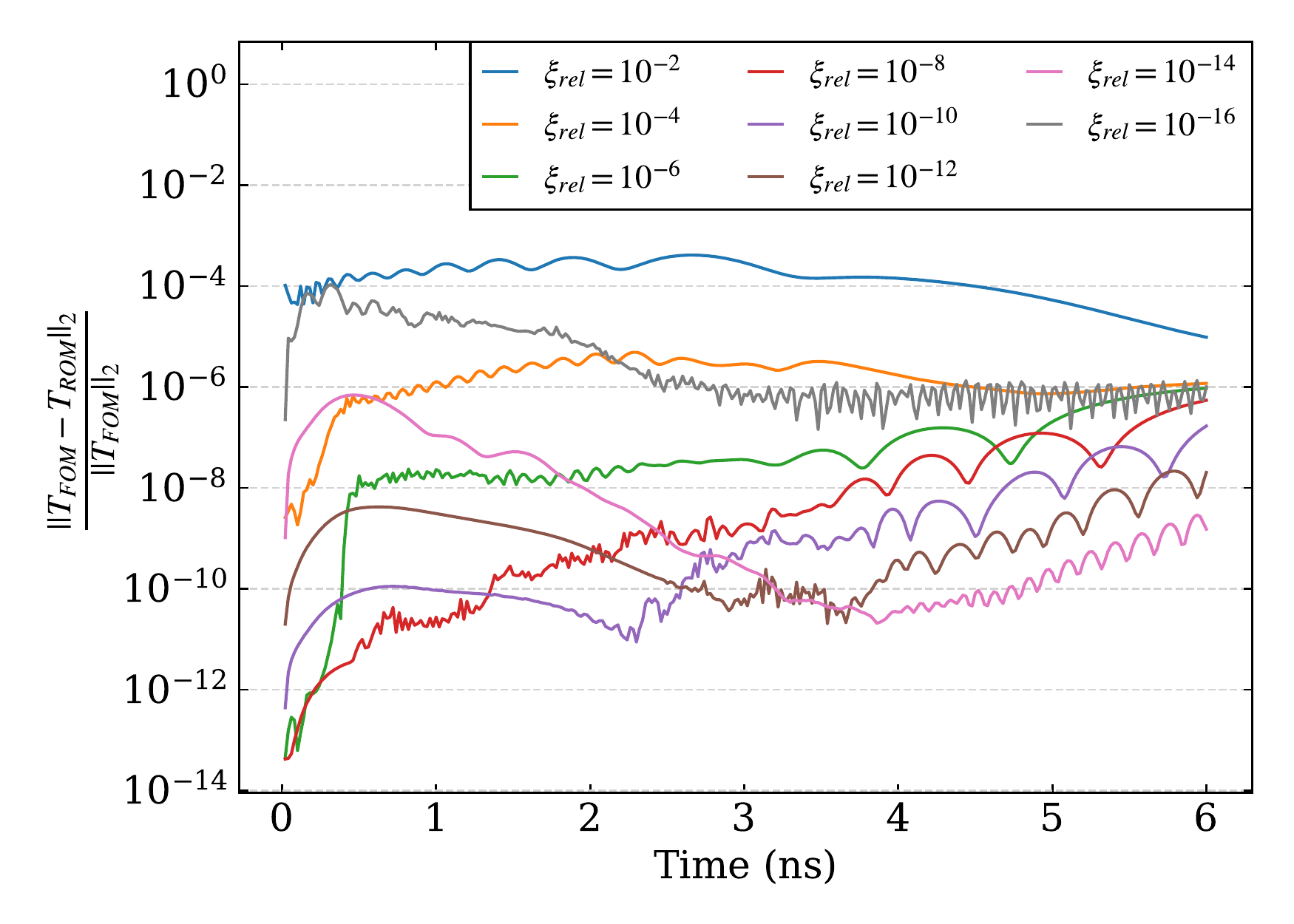}}
	\subfloat[Radiation Energy Density]{\includegraphics[width=.5\textwidth]{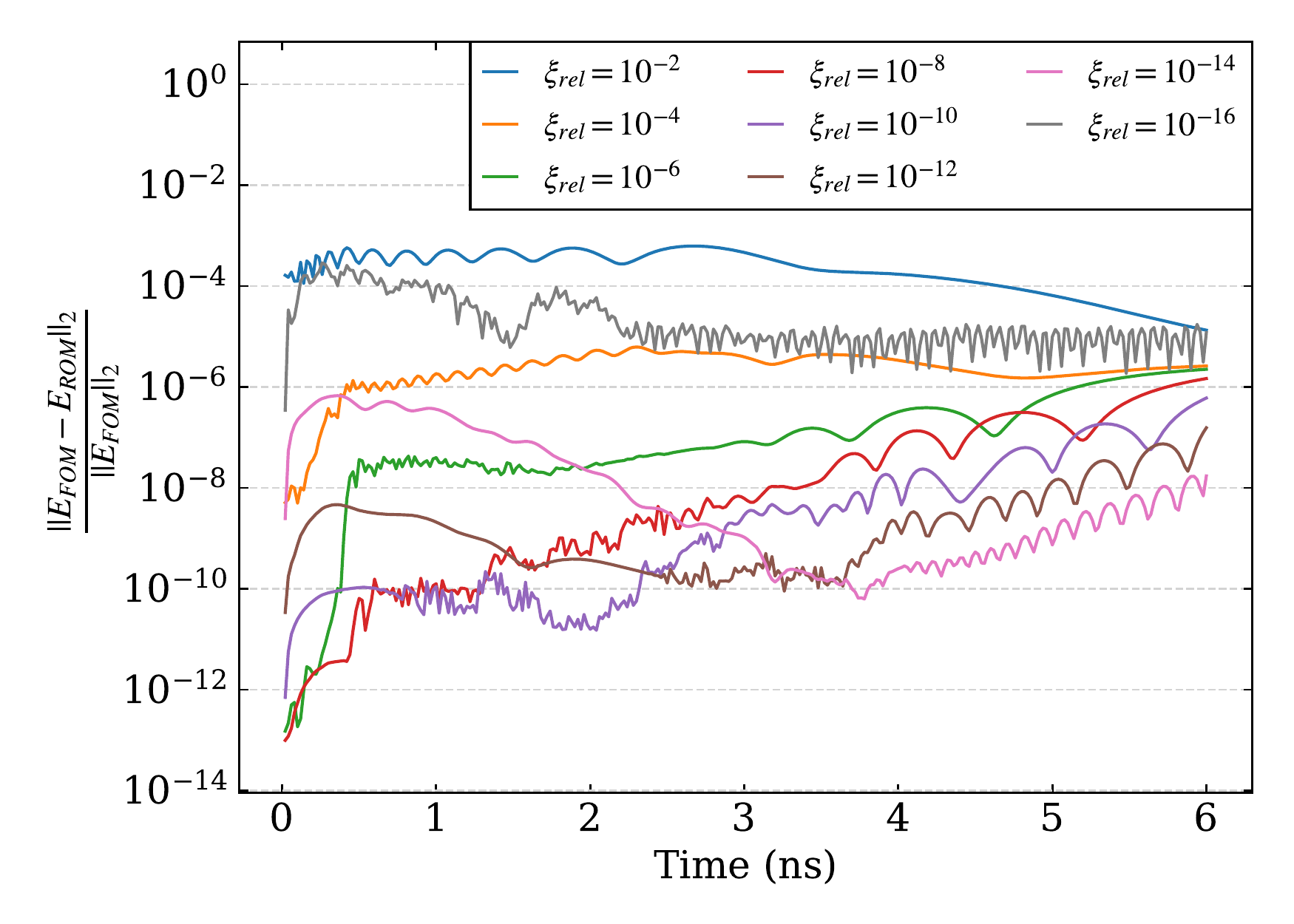}}
	\caption{Relative errors in the 2-norm of the DDET ROM using the DMD-E using several $\xi_{\text{rel}}$, plotted vs time}
	\label{fig:2nrm-errs_DMDB}
	\vspace{.25cm}
\end{figure}

We now analyze the solutions of the F-C test computed by DDET ROMs
with the reduced-rank databases $\mat{A}^{\tens{f}_{\alpha\beta,\gamma}}_k$ and $\mat{A}^c_k$.
Figures \ref{fig:2nrm-errs_POD}, \ref{fig:2nrm-errs_DMD} and \ref{fig:2nrm-errs_DMDB} show the relative error for the material temperature $(T)$ and total radiation energy density $(E)$ calculated in the 2-norm over space at each instant of time in $t\in[0,6\text{ns}]$ where each unique curve corresponds to the ROM solution generated for a given value of $\xi_{\text{rel}}$. The discrete FOM solution is chosen as the reference to compute errors against to determine how the ROM solution converges to its training data.
The only errors incurred by our FOM are due to discretization and as such our FOM  will converge to the multigroup TRT solution in the limit $N_x,N_y,N_\Omega,N_t\rightarrow\infty$.
Therefore we
postulate that if the solution of the DDET model converges to the discrete FOM solution then it will too converge to the continuous solution given a database generated on a fine-enough grid.

The relative error
of the ROM using the POD for a given $\xi_{\text{rel}}$ (Figure \ref{fig:2nrm-errs_POD})
first increases during the initial stage of wave evolution and then stabilizes in time slightly varying around some value.
The errors  exhibit uniform convergence with decreasing $\xi_{\text{rel}}$.
The lowest-rank POD approximation ($\xi_{\text{rel}}=10^{-2}$) yields errors on the order of $10^{-4}$ . Given the full-rank POD (i.e. $\xi_{\text{rel}}=10^{-16}$), this DDET ROM successfully reproduces the reference FOM solution within the numerical convergence bounds.

The ROM using the DMD (Figure \ref{fig:2nrm-errs_DMD}) shows similar performance to the ROM with POD, although with lower accuracy for each $\xi_{\text{rel}}$. For instance, the lowest-rank DMD approximation yields errors on the order of $10^{-2}$. Despite this increase in error, the ROM with DMD still approaches the FOM solution when using full-rank representations of each database. Here the FOM solution is recreated within
an error of $10^{-12}$ only with the exception of times roughly within the range $t\in[2.5,5\text{ns}]$ where relative error still remains very low, on the order of $10^{-10}$. We attribute this effect to numerical noise of which the DMD is especially susceptible to as rank is increased \cite{tu-rowley-2014,chen-tu-2012}. This effect at very small $\xi_{\text{rel}}$ is also present for the ROM equipped with the DMD-E, but amplified by a considerable level. The ROM with low-rank DMD-E approximations of each snapshot matrix (Figure \ref{fig:2nrm-errs_DMDB}) actually yields
very similar levels of accuracy to the ROM with POD, achieving errors on the order of $10^{-4}$ with $\xi_{\text{rel}}=10^{-2}$. These errors are only observed to decrease with $\xi_{\text{rel}}$ up until  $\xi_{\text{rel}}=10^{-10}$  however, and as $\xi_{\text{rel}}$ is decreased further the errors are observed to increase instead of stagnate as seen with the DMD. In fact, the relative error is comparable for the ROM with the DMD-E for both $\xi_{\text{rel}}=10^{-2}$ and $\xi_{\text{rel}}=10^{-16}$.

\begin{figure}[ht!]
	\vspace{-.25cm}
	\centering
	\subfloat[Material Temperature]{\includegraphics[width=.5\textwidth]{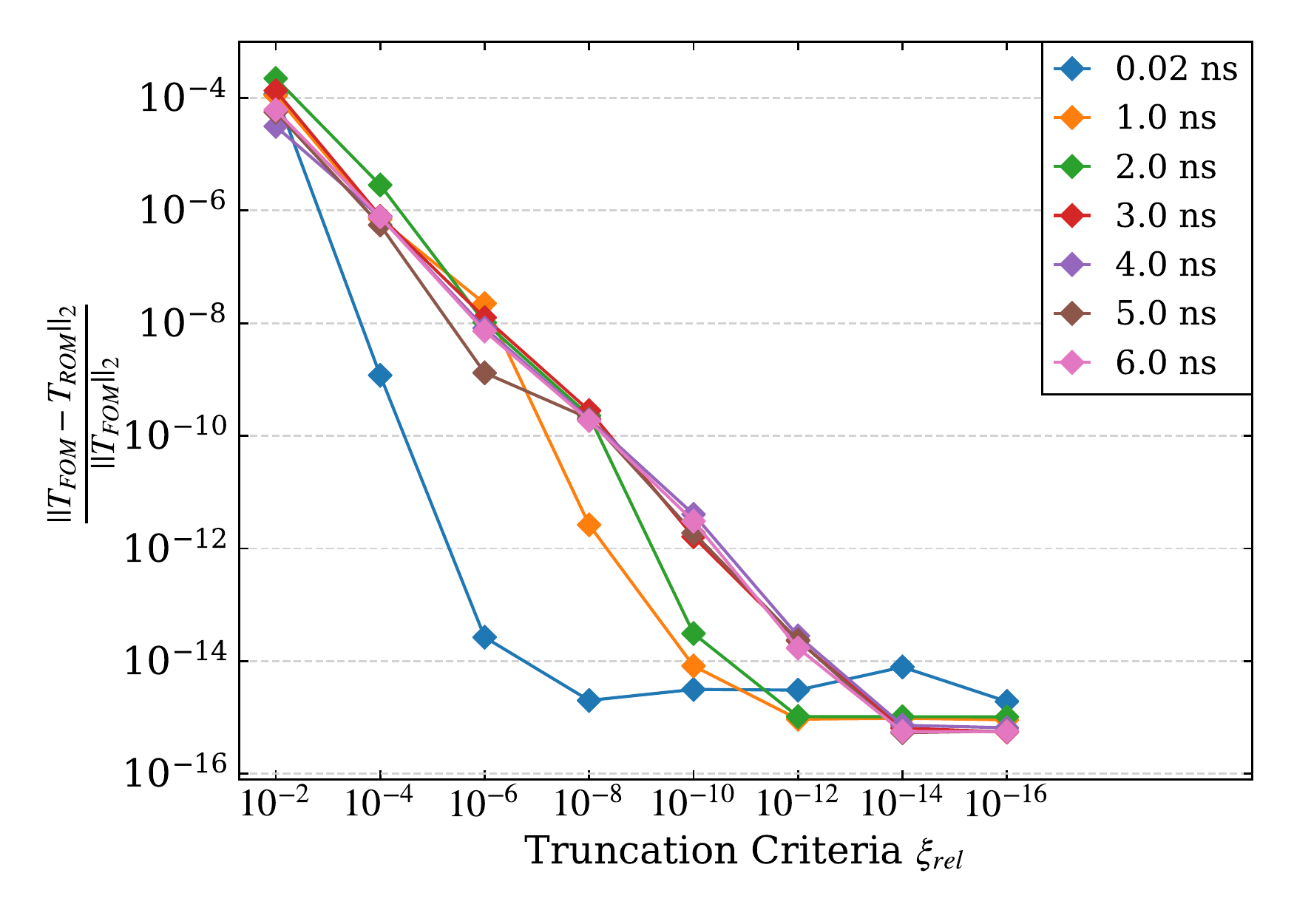}}
	\subfloat[Radiation Energy Density]{\includegraphics[width=.5\textwidth]{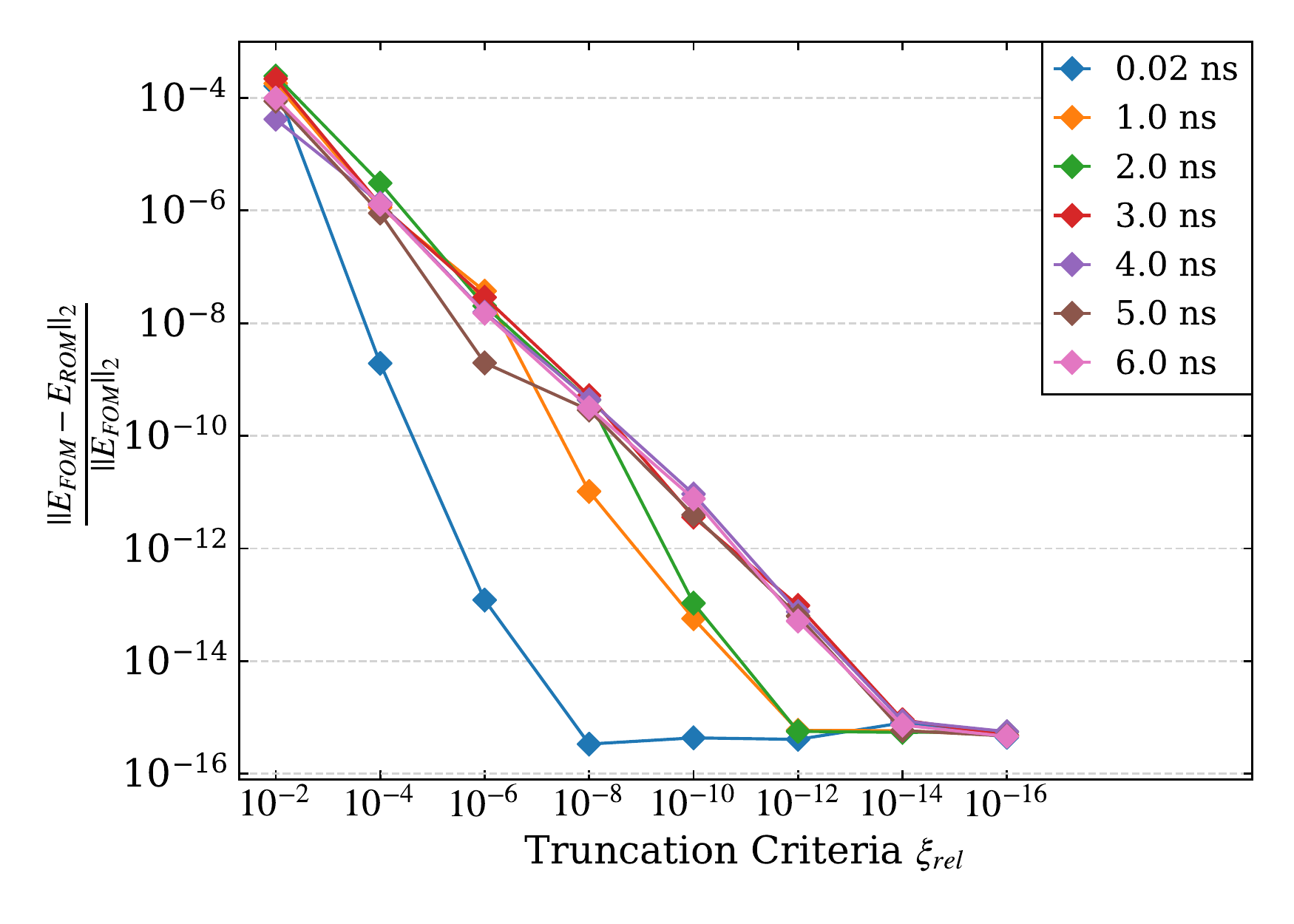}}
	\caption{Relative errors in the 2-norm of the DDET ROM using the POD at several times, plotted vs $\xi_{\text{rel}}$}
	\label{fig:2nrm-errs_trend_POD}
	\centering
	\subfloat[Material Temperature]{\includegraphics[width=.5\textwidth]{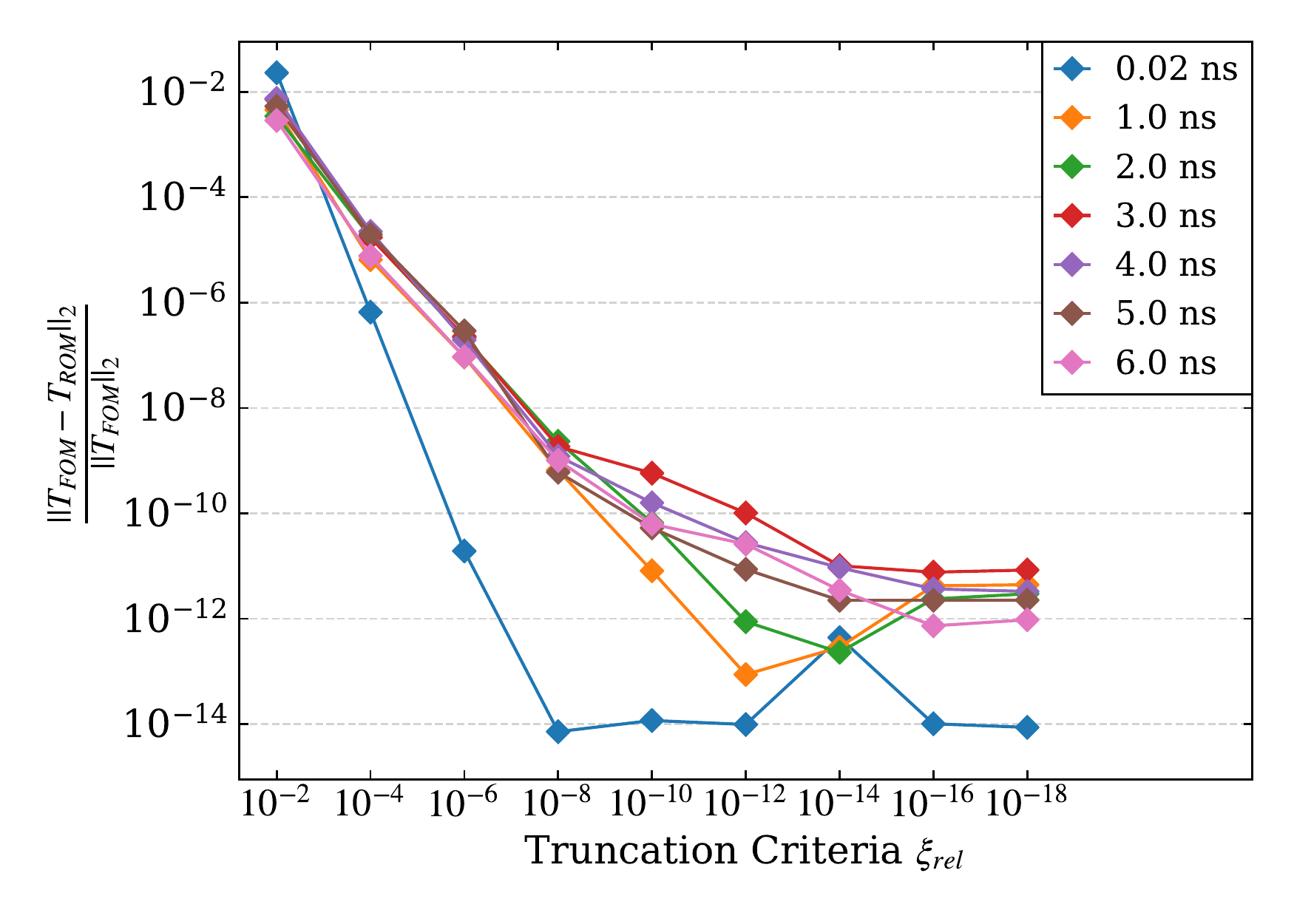}}
	\subfloat[Radiation Energy Density]{\includegraphics[width=.5\textwidth]{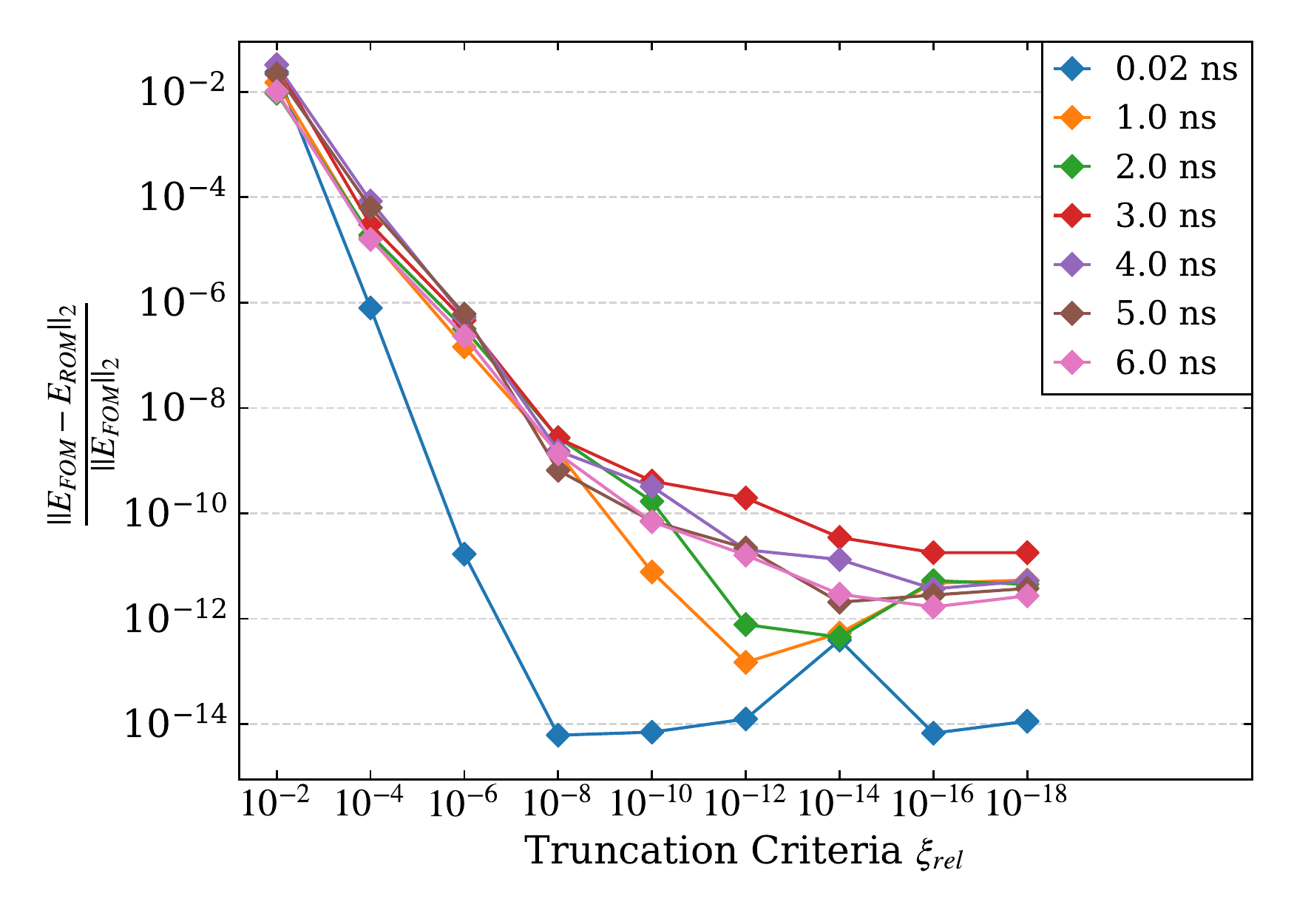}}
	\caption{Relative errors in the 2-norm of the DDET ROM using the DMD at several times, plotted vs $\xi_{\text{rel}}$}
	\vspace{.5cm}
	\label{fig:2nrm-errs_trend_DMD}
	\vspace{-.25cm}
\end{figure}

\begin{figure}[ht!]
	\vspace*{-.5cm}
	\centering
	\subfloat[Material Temperature]{\includegraphics[width=.5\textwidth]{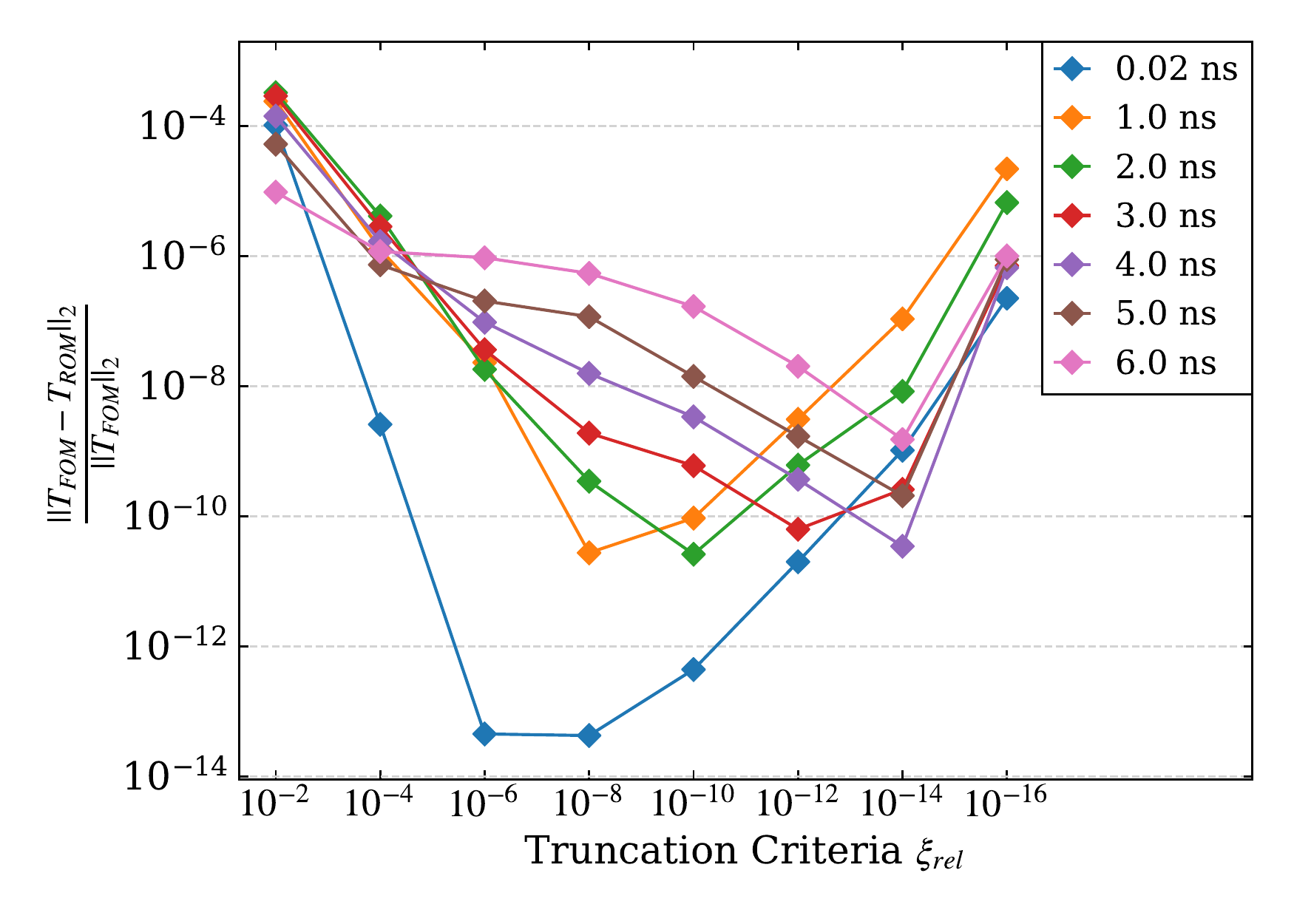}}
	\subfloat[Radiation Energy Density]{\includegraphics[width=.5\textwidth]{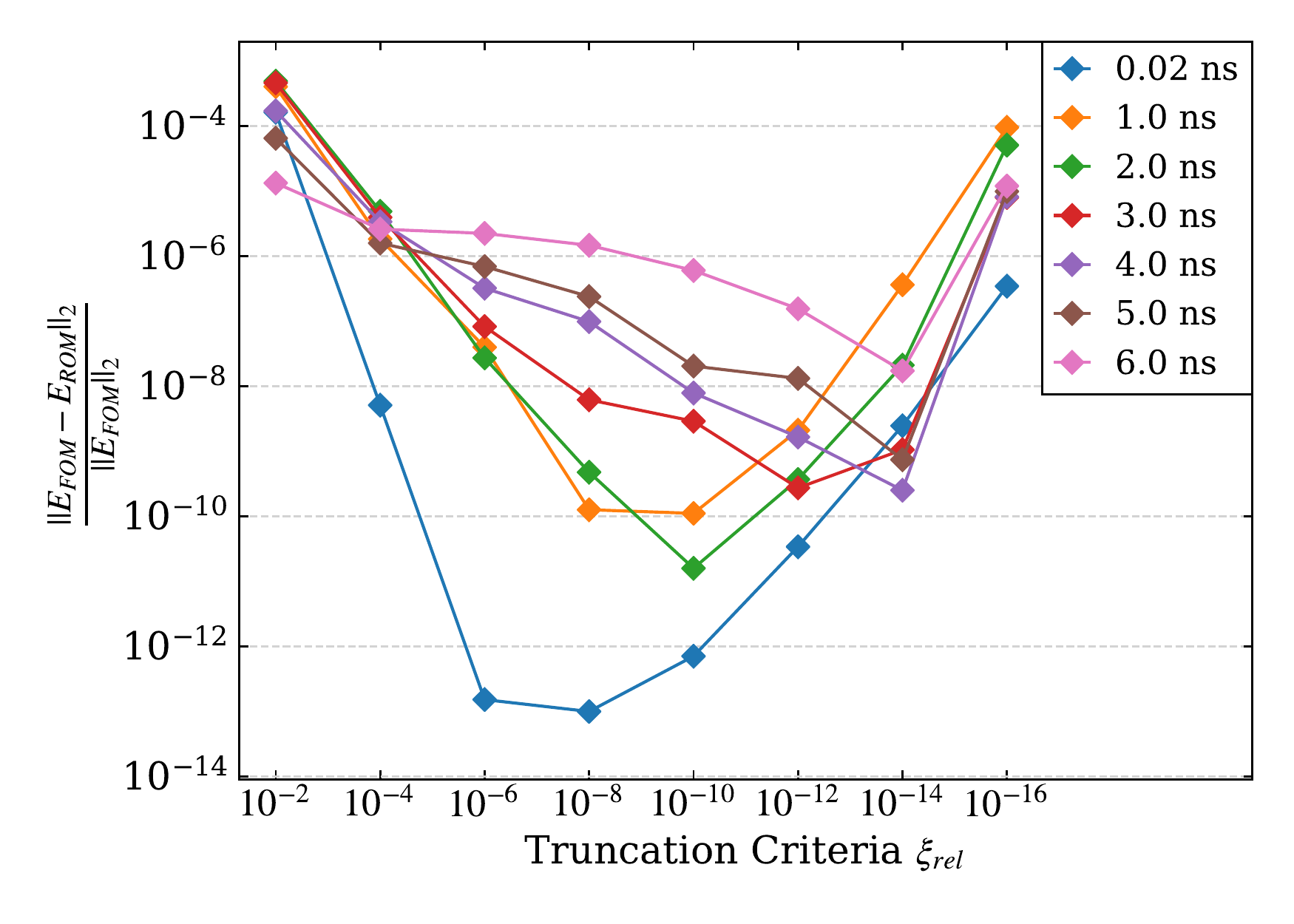}}
	\caption{Relative errors in the 2-norm of the DDET ROM using the DMD-E at several times, plotted \linebreak vs $\xi_{\text{rel}}$}
	\label{fig:2nrm-errs_trend_DMDB}
\end{figure}

Figures \ref{fig:2nrm-errs_trend_POD}, \ref{fig:2nrm-errs_trend_DMD} and \ref{fig:2nrm-errs_trend_DMDB}
demonstrate  the ROMs  convergence with increases in the rank of approximation of data.
Here each curve corresponds to a specific instant of time, showing how the error of the ROMs  changes  with respect to $\xi_{\text{rel}}$.
The results show that the solution of the POD-based ROM converges to the reference FOM solution
 linearly with respect to $\xi_{\text{rel}}$.
There is stagnation of the relative error around $10^{-14}$ due to limitation of the finite precision of calculations.
 Furthermore
there is an apparent correlation that the relative error of the ROM in the 2-norm is proportional to $10^{-2}\xi_{\text{rel}}$
for both $T$ and $E$ until the limits of finite precision become overwhelming.
The same convergence behavior with $\xi_{\text{rel}}$ is also demonstrated for
the solution of the ROM with the DMD.
With the DMD, the relative errors in the 2-norm of the ROM are proportional  to $\xi_{\text{rel}}$ for both $T$ and $E$ while $\xi_{\text{rel}}\geq 10^{-10}$.
The numerical noise affects convergence  for $\xi_{\text{rel}} <10^{-10}$ and limit the relative error.
The errors largely stagnate
in the case of $\xi_{\text{rel}}\in [10^{-10},10^{-18}]$.
Increases in the errors are seen at $\xi_{\text{rel}} = 10^{-14},10^{-16}$ for times $t=1,2$ns, and a temporary increase in error for time $t=0.02$ns is observed at $\xi_{\text{rel}} = 10^{-14}$.

 Finally considering the ROM with the DMD-E, Figure \ref{fig:2nrm-errs_trend_DMDB} clearly demonstrates numerical instability for small $\xi_{\text{rel}}$.
An initial increase in error level is seen at  times $t=.02,1$ns for
$\xi_{\text{rel}}=10^{-10}$, and the errors at
 $t=2$ns increase  at $\xi_{\text{rel}}=10^{-12}$. When $\xi_{\text{rel}}=10^{-16}$, the errors in the ROM solution have increased back to the level observed for $\xi_{\text{rel}}=10^{-2}$.
 This behavior is accredited  to a large magnification of numerical errors as seen with the DMD.
The DMD-E can be interpreted as the DMD on a set of residual vectors representing the distance of the decomposed data to the steady-state solution. The residual vectors for near steady-state data are then expected to have elements of very small magnitude which can contribute to numerical issues. This combined with the inherent sensitivity of the DMD to numerics can lead to large amplifications of error. Such numerical problems do not necessitate an abandonment of the DMD-E however, as only the DDET ROMs with especially high-rank data approximations are impaired.

\begin{figure}[ht!]
	\centering
	\subfloat[Material Temperature]{
	\begin{tabular}{c|c|c|c}
		$\xi_\text{rel} $& t=1ns & t=2ns & t=3ns \\ \hline &&& \\[-.3cm]
		%%%
		$10^{-2}$ &
		\raisebox{-.5\height}{\includegraphics[trim=1cm 0cm 4cm .5cm,clip,height=3.5cm]{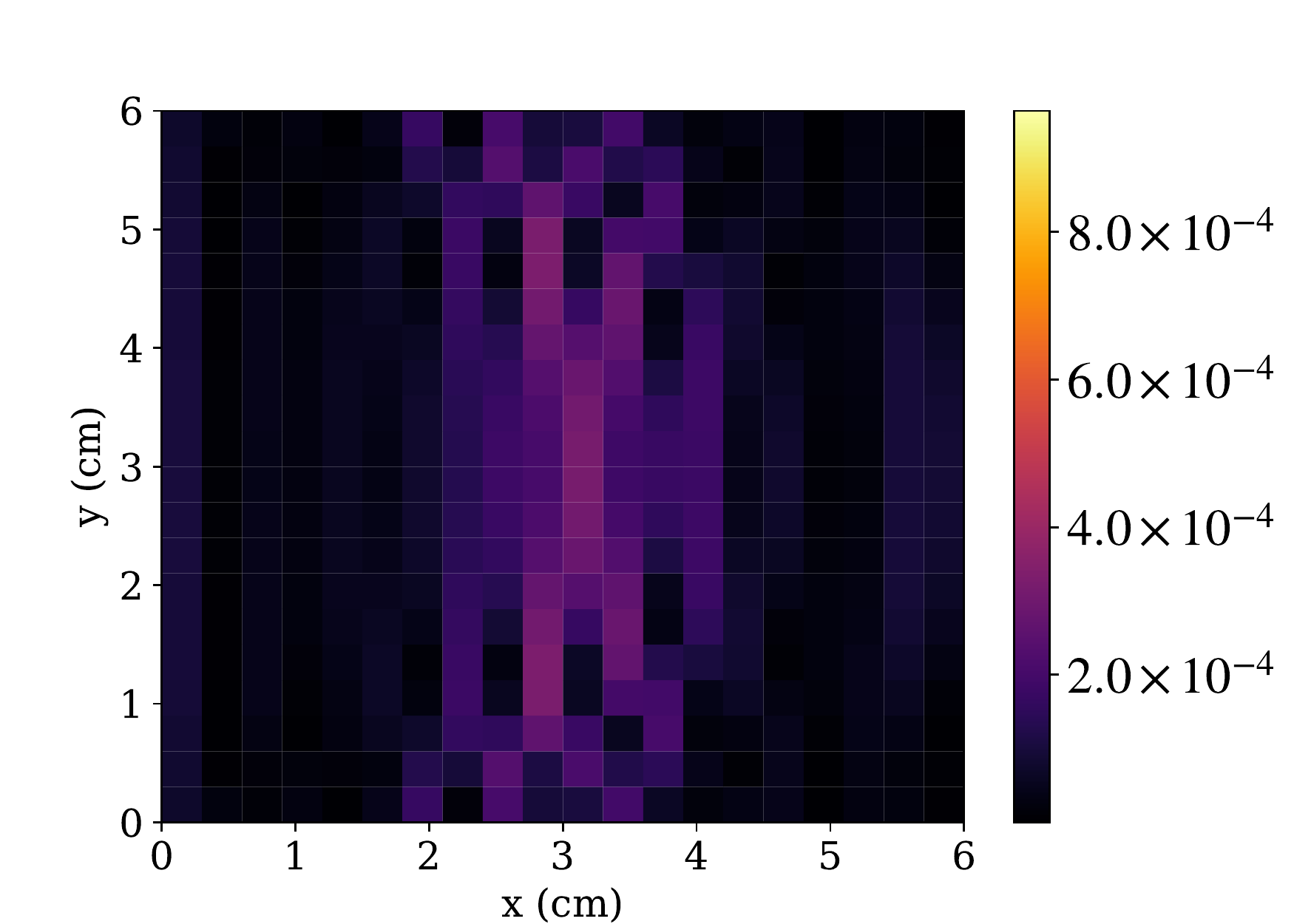}}&
		\raisebox{-.5\height}{\includegraphics[trim=1cm 0cm 4cm .5cm,clip,height=3.5cm]{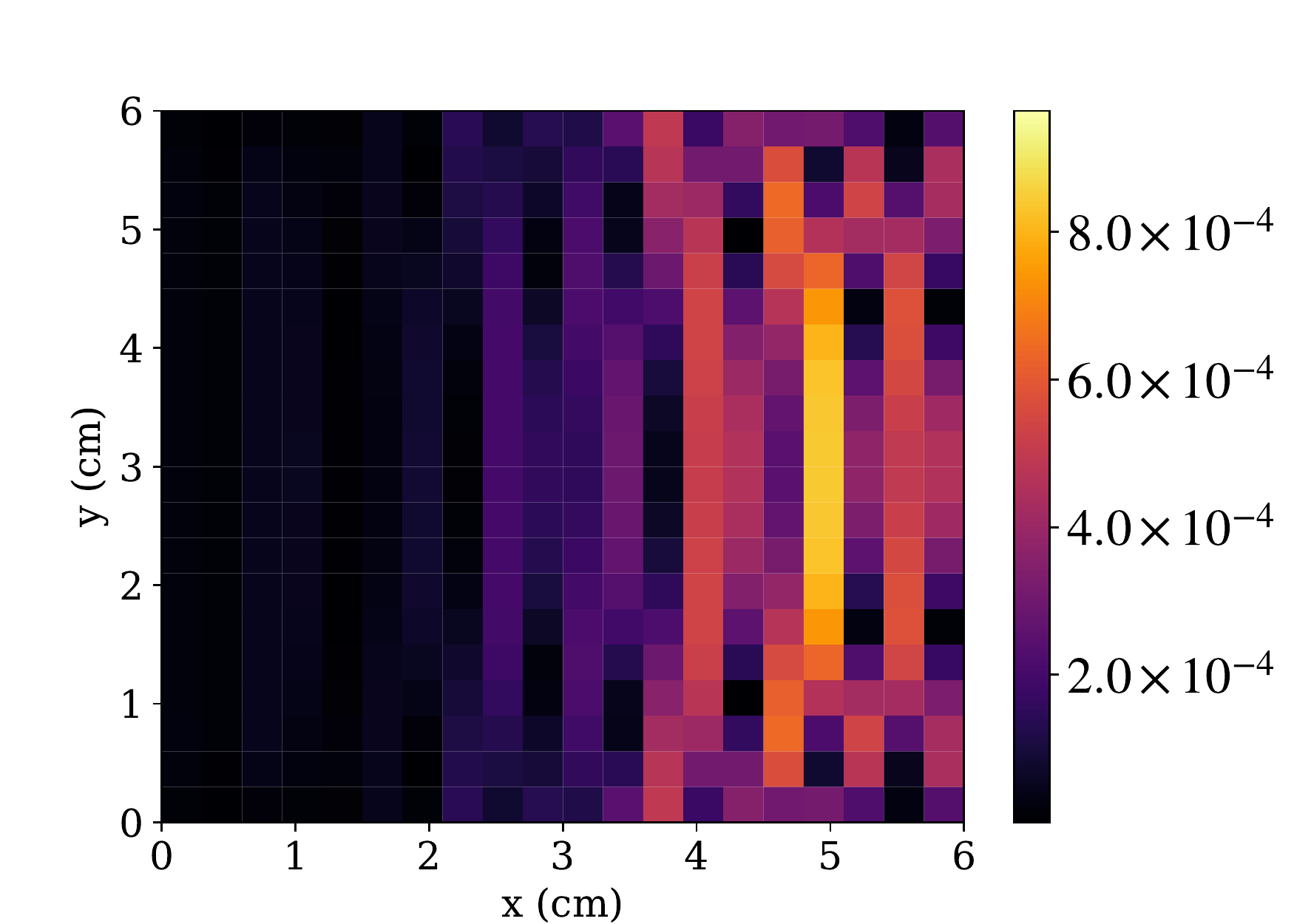}}&
		\raisebox{-.5\height}{\includegraphics[trim=1cm 0cm 0cm .5cm,clip,height=3.5cm]{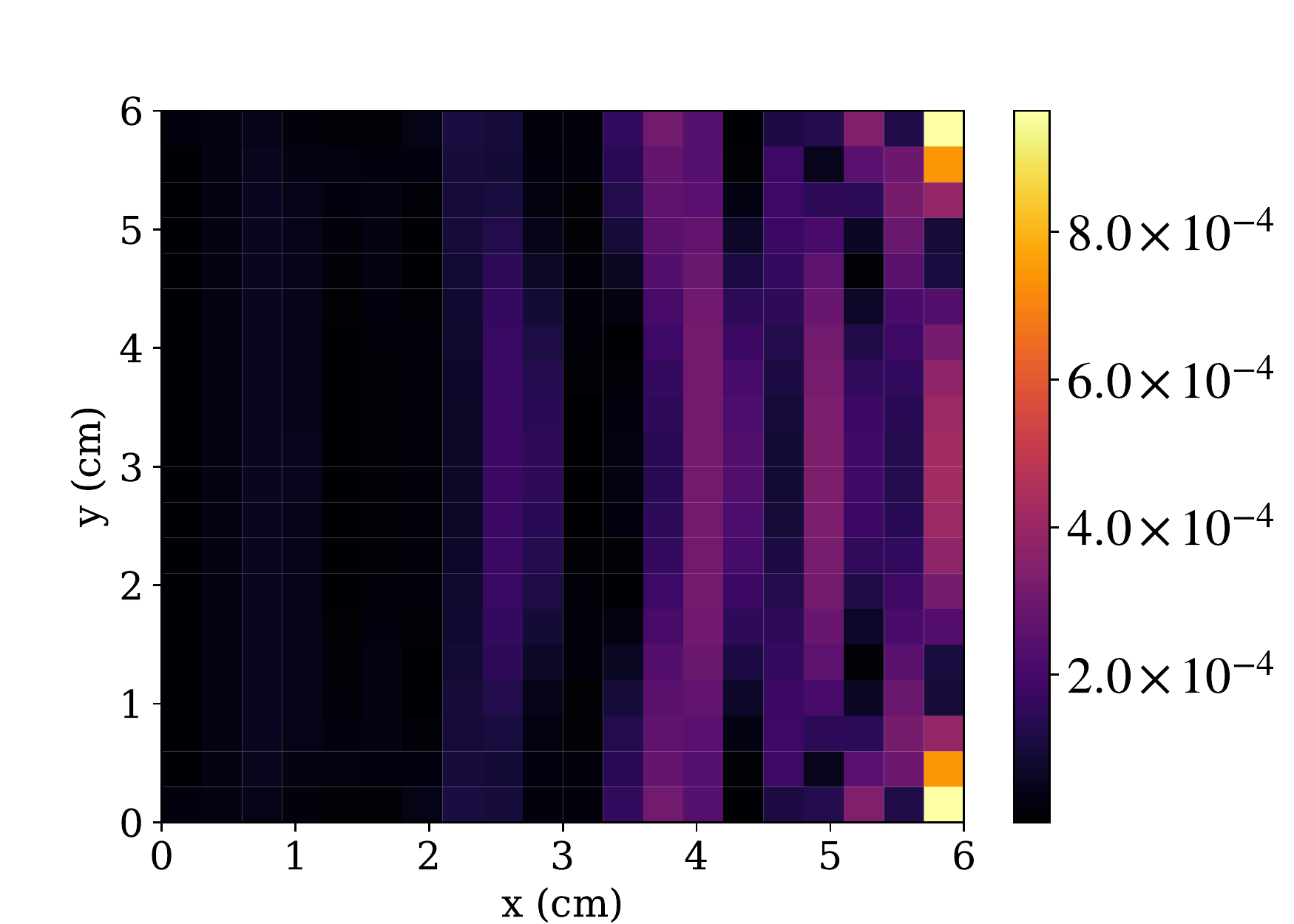}}\\[.1cm] \hline &&& \\[-.3cm]
		%%%
		$10^{-4}$ &
		\raisebox{-.5\height}{\includegraphics[trim=1cm 0cm 4cm .5cm,clip,height=3.5cm]{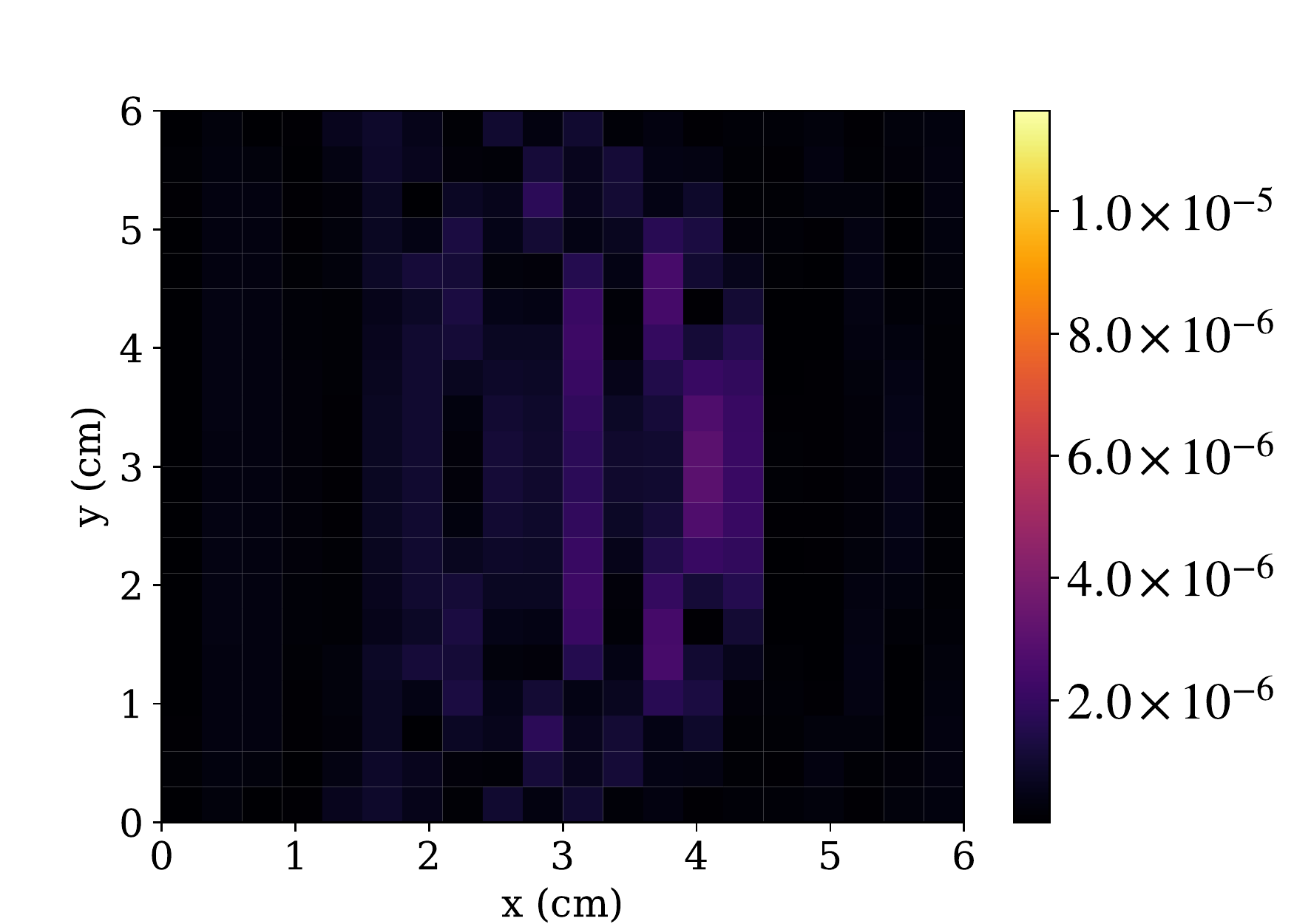}}&
		\raisebox{-.5\height}{\includegraphics[trim=1cm 0cm 4cm .5cm,clip,height=3.5cm]{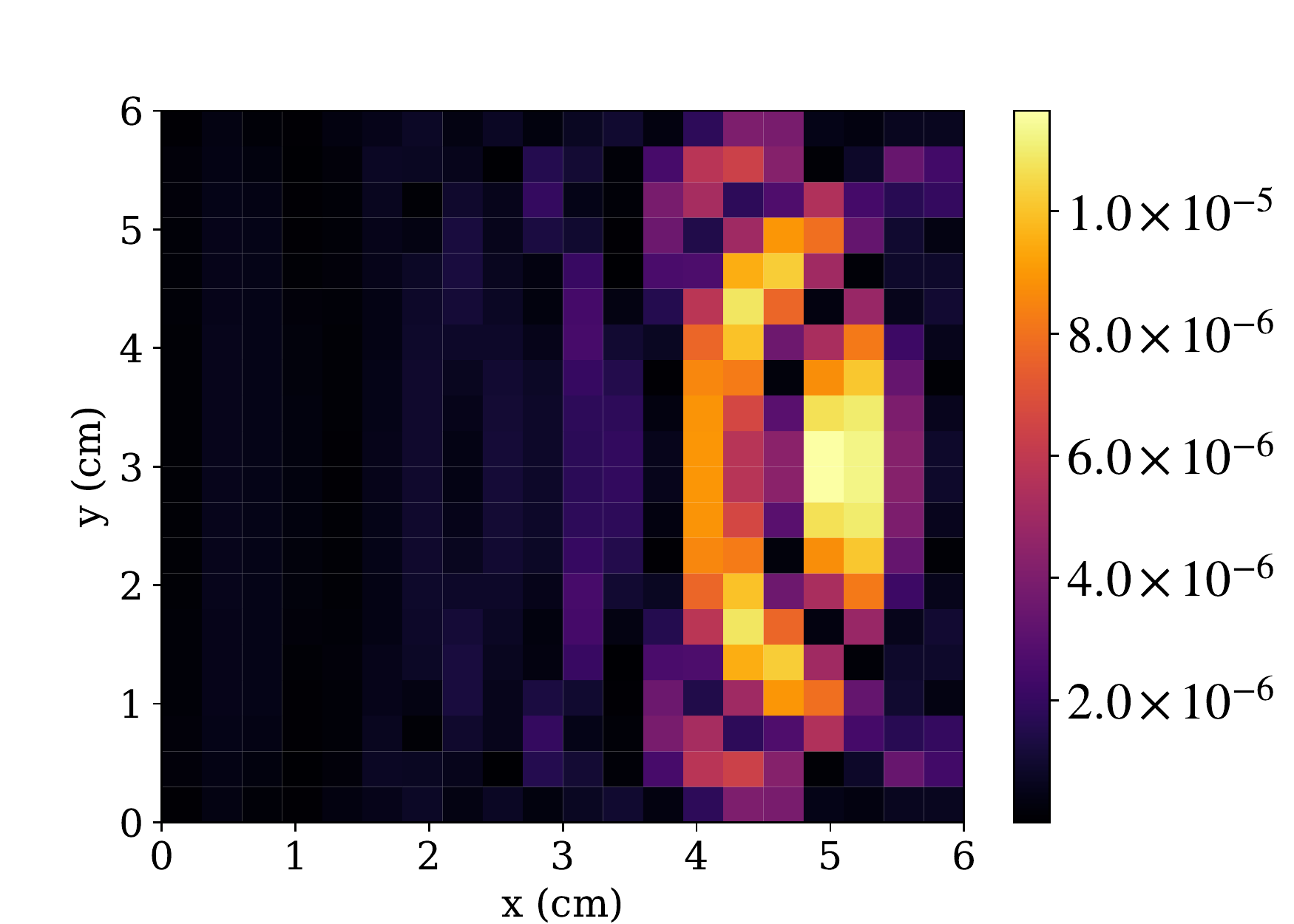}}&
		\raisebox{-.5\height}{\includegraphics[trim=1cm 0cm 0cm .5cm,clip,height=3.5cm]{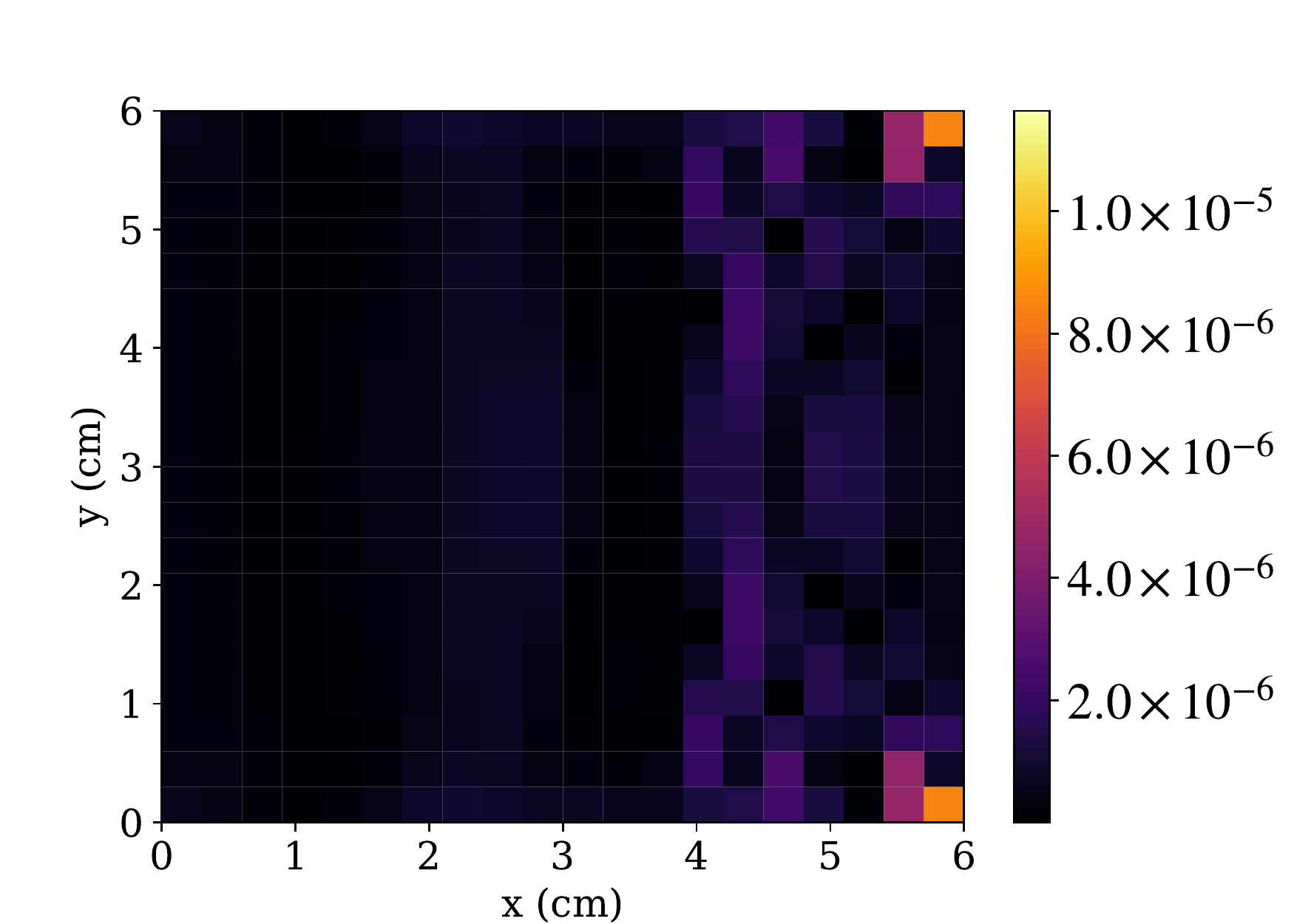}}
	\end{tabular} }
	\vspace*{.25cm}
	\subfloat[Radiation Energy Density]{
	\begin{tabular}{c|c|c|c}
		$\xi_\text{rel} $& t=1ns & t=2ns & t=3ns \\ \hline &&& \\[-.3cm]
		%%%
		%%%
		$10^{-2}$ &
		\raisebox{-.5\height}{\includegraphics[trim=1cm 0cm 4cm .5cm,clip,height=3.5cm]{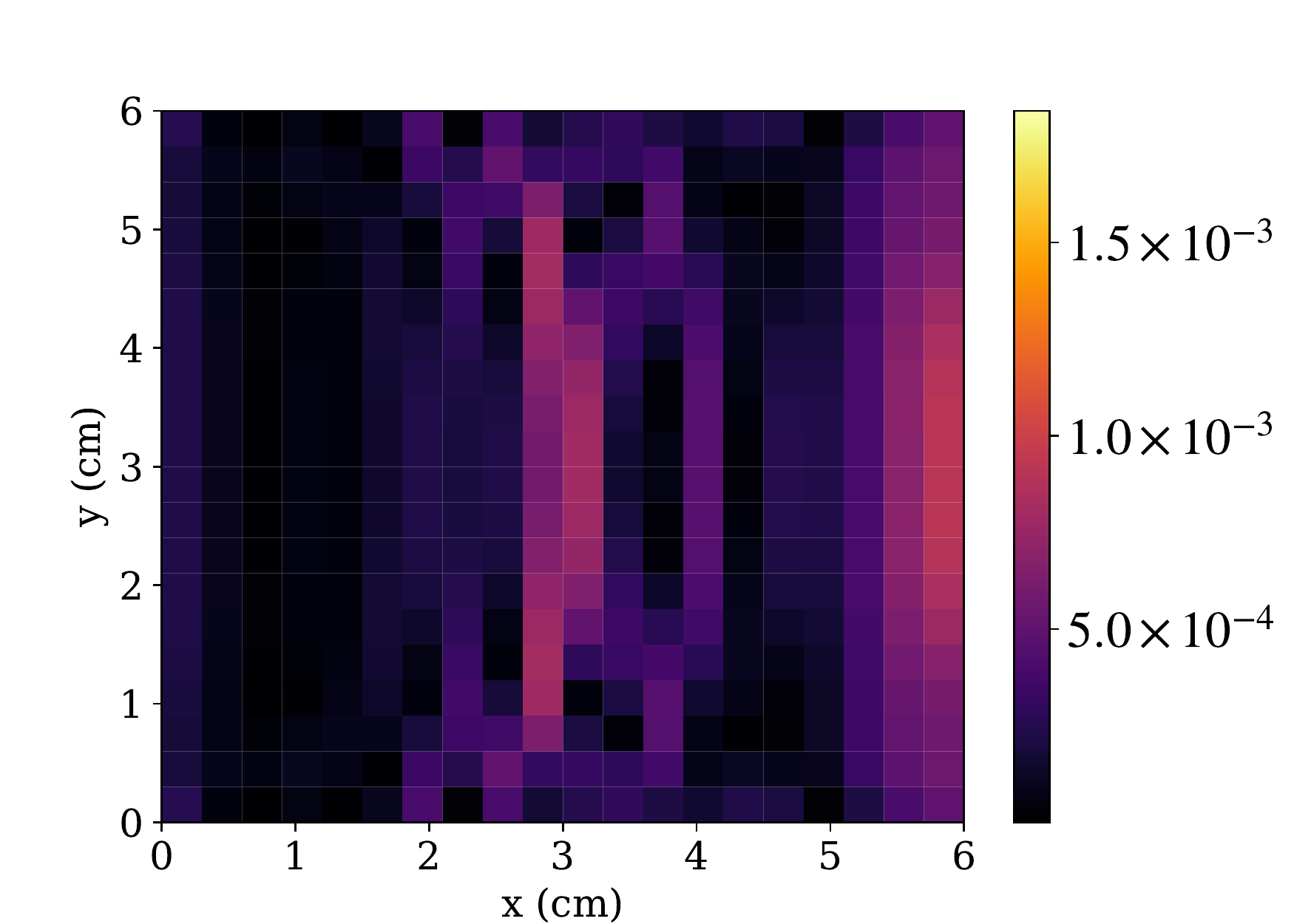}}&
		\raisebox{-.5\height}{\includegraphics[trim=1cm 0cm 4cm .5cm,clip,height=3.5cm]{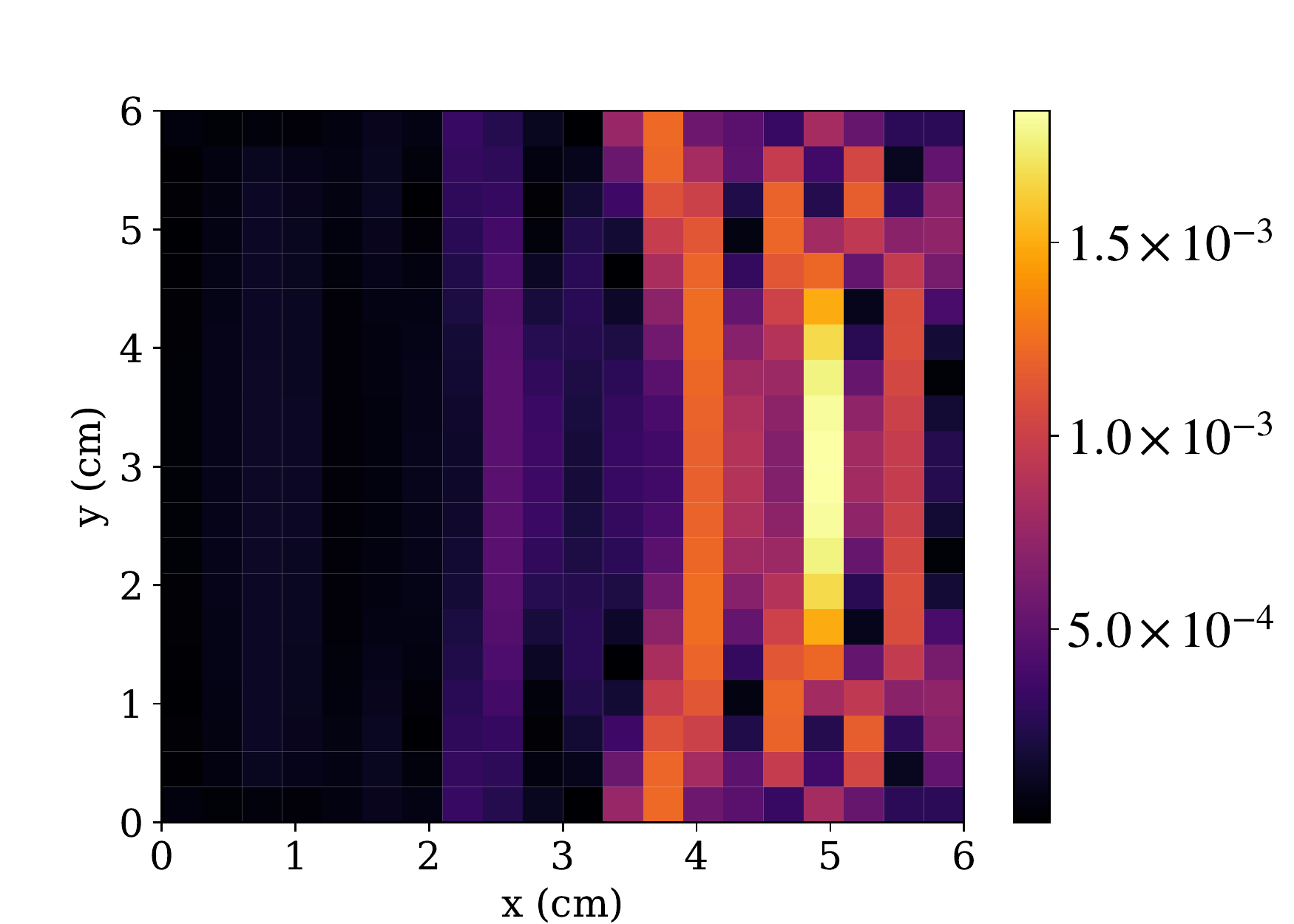}}&
		\raisebox{-.5\height}{\includegraphics[trim=1cm 0cm 0cm .5cm,clip,height=3.5cm]{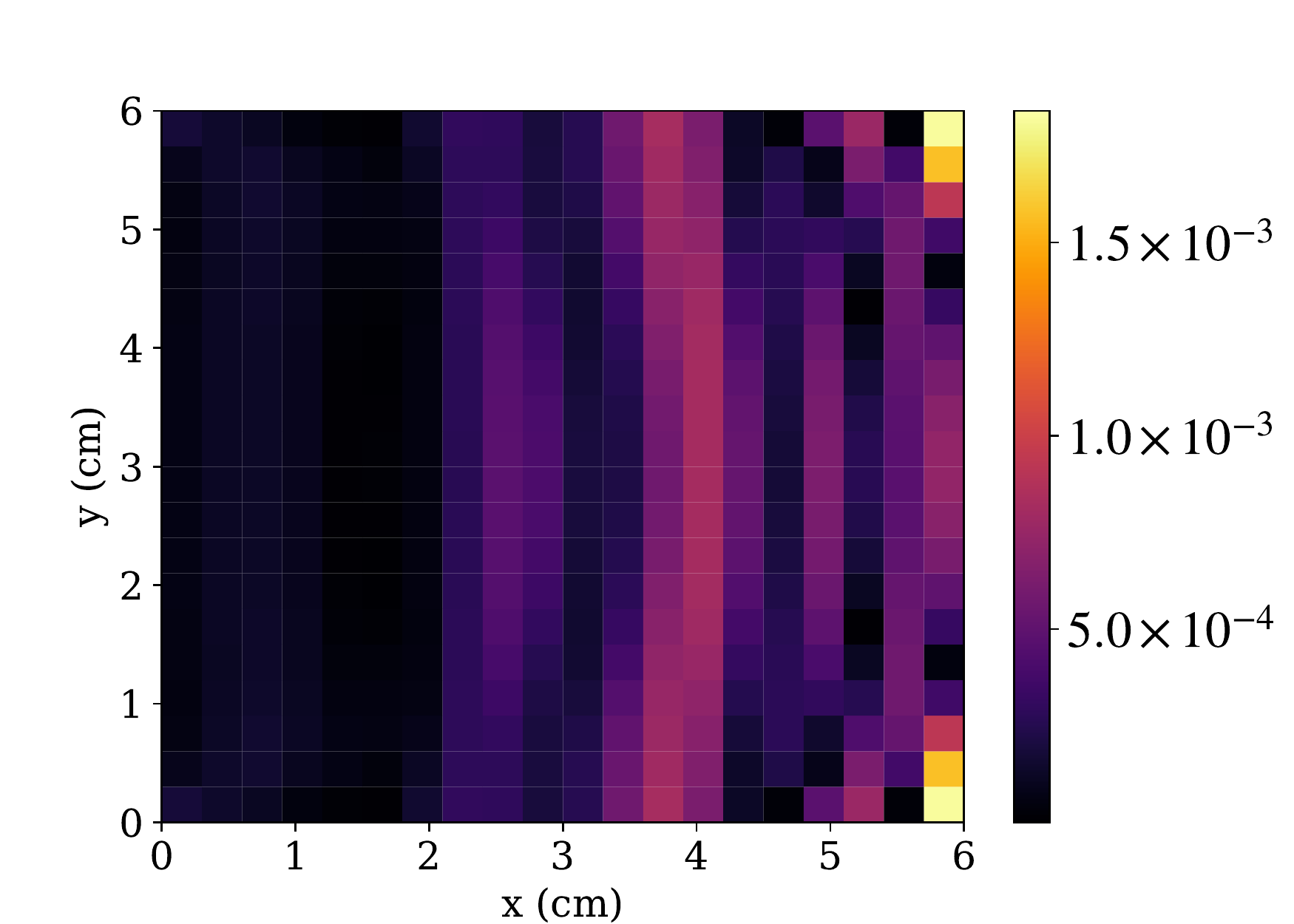}}\\[.1cm] \hline &&& \\[-.3cm]
		%%%
		$10^{-4}$ &
		\raisebox{-.5\height}{\includegraphics[trim=1cm 0cm 4cm .5cm,clip,height=3.5cm]{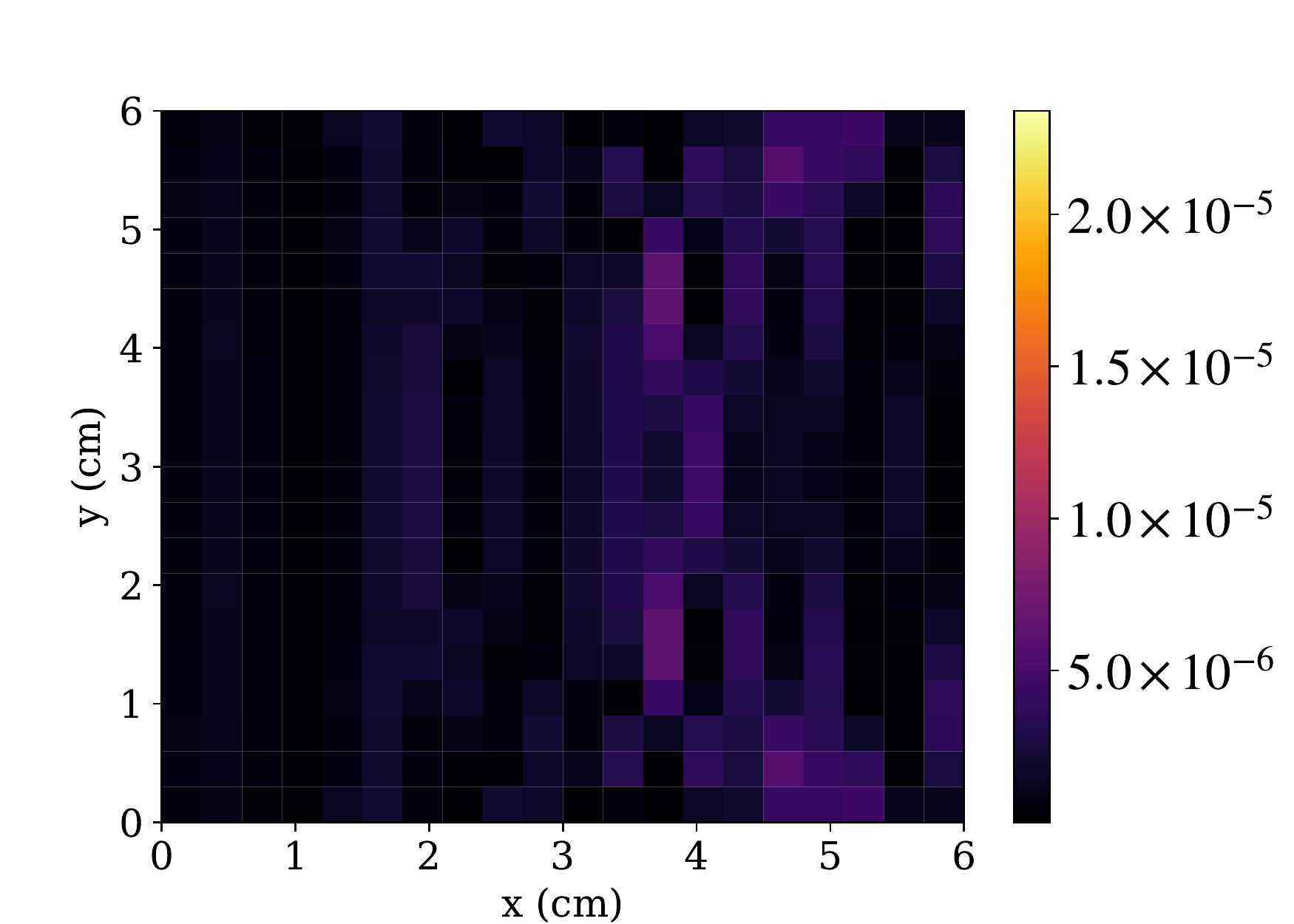}}&
		\raisebox{-.5\height}{\includegraphics[trim=1cm 0cm 4cm .5cm,clip,height=3.5cm]{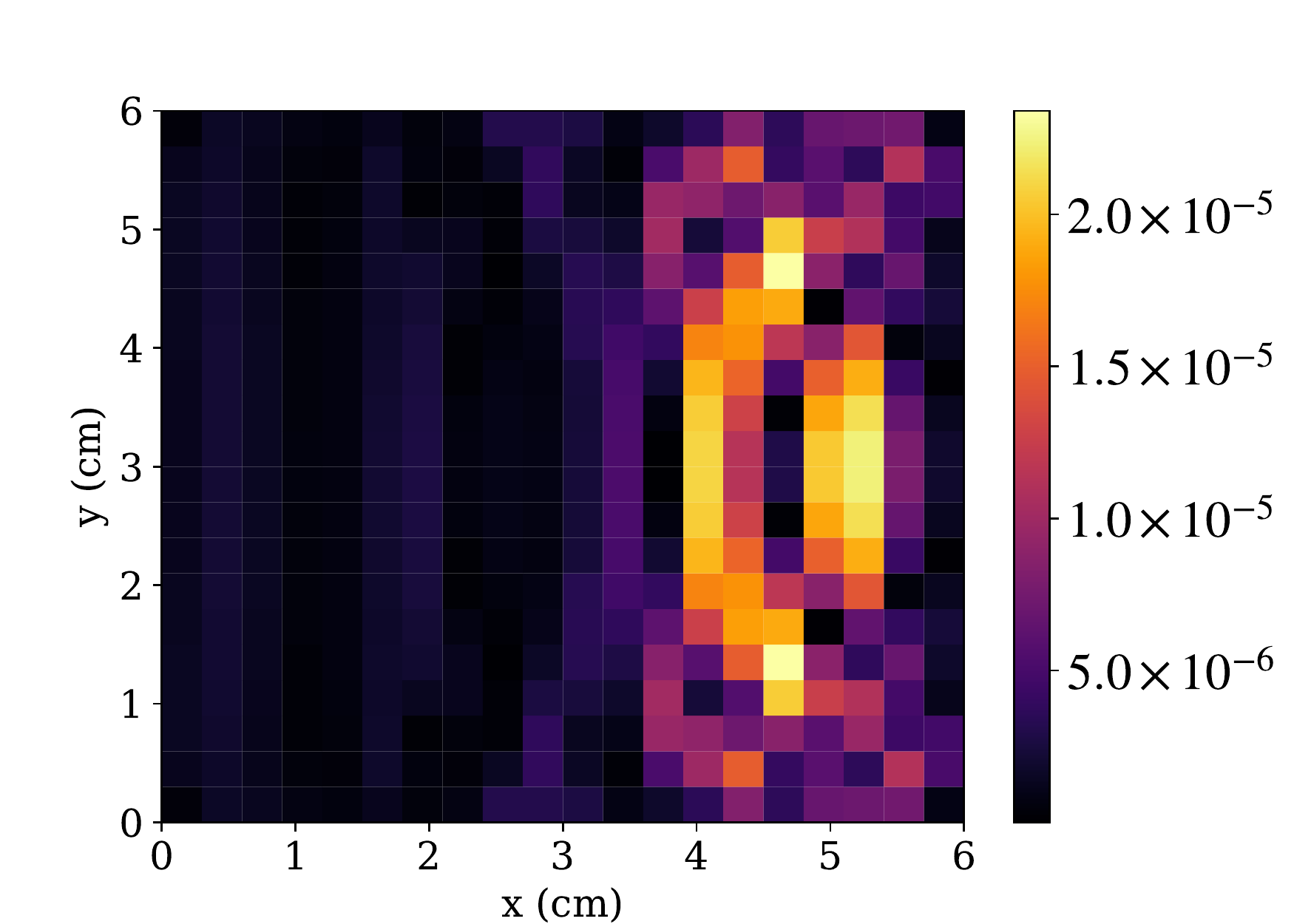}}&
		\raisebox{-.5\height}{\includegraphics[trim=1cm 0cm 0cm .5cm,clip,height=3.5cm]{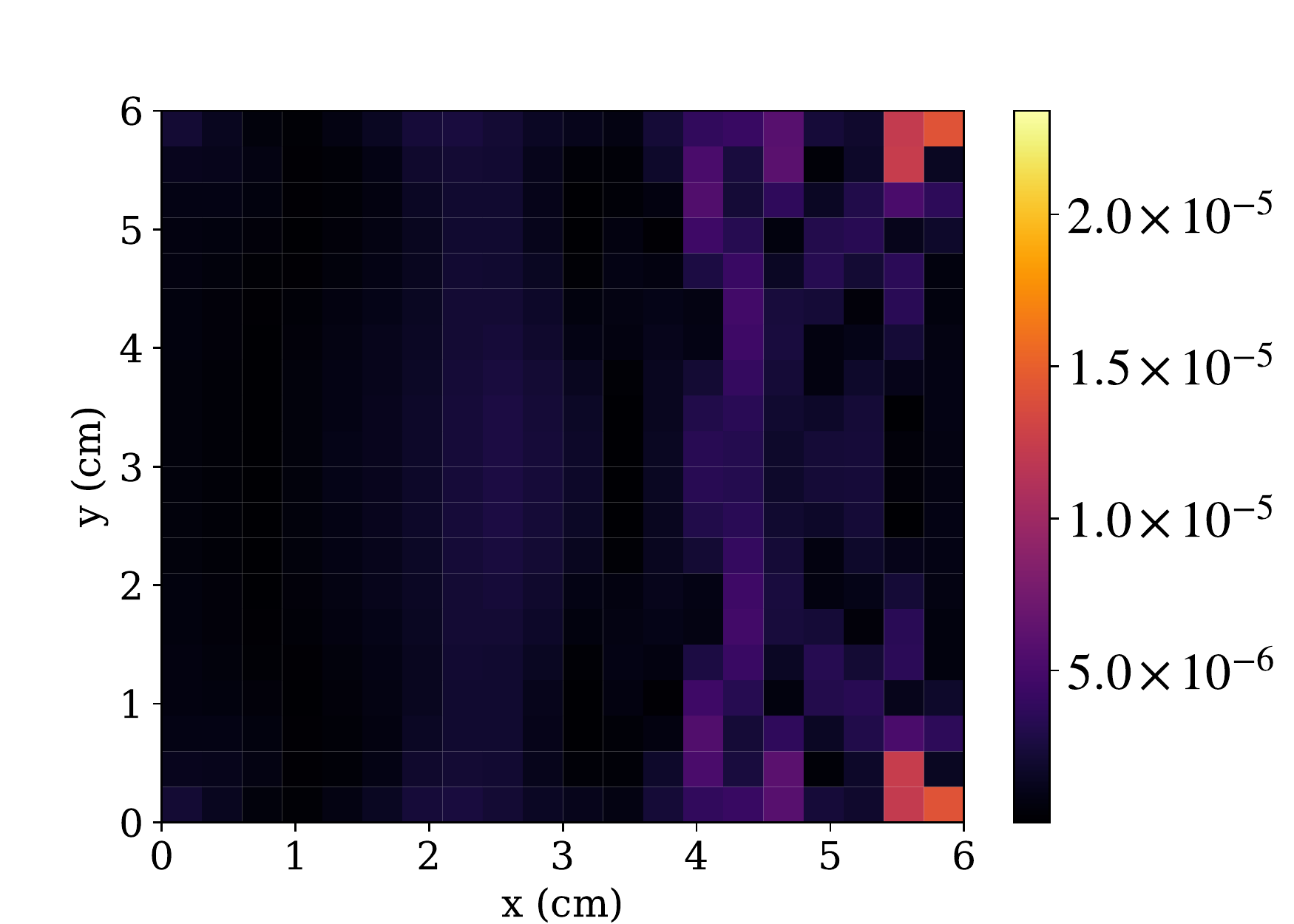}}
	\end{tabular} }
	\vspace*{.25cm}
	\caption{Cell-wise relative error in material temperature $(T)$ and total radiation energy density $(E)$ over the spatial domain at times t=1, 2, 3 ns for the DDET ROM equipped with the POD for $\xi_{\text{rel}}=10^{-2},10^{-4}$.}

	\label{fig:spatial_errs_pod}
\end{figure}

\begin{figure}[ht!]
	\centering
	\subfloat[Material Temperature]{
		\begin{tabular}{c|c|c|c}
			$\xi_\text{rel} $& t=1ns & t=2ns & t=3ns \\ \hline &&& \\[-.3cm]
			%%%
			$10^{-2}$ &
			\raisebox{-.5\height}{\includegraphics[trim=1cm 0cm 4cm .5cm,clip,height=3.5cm]{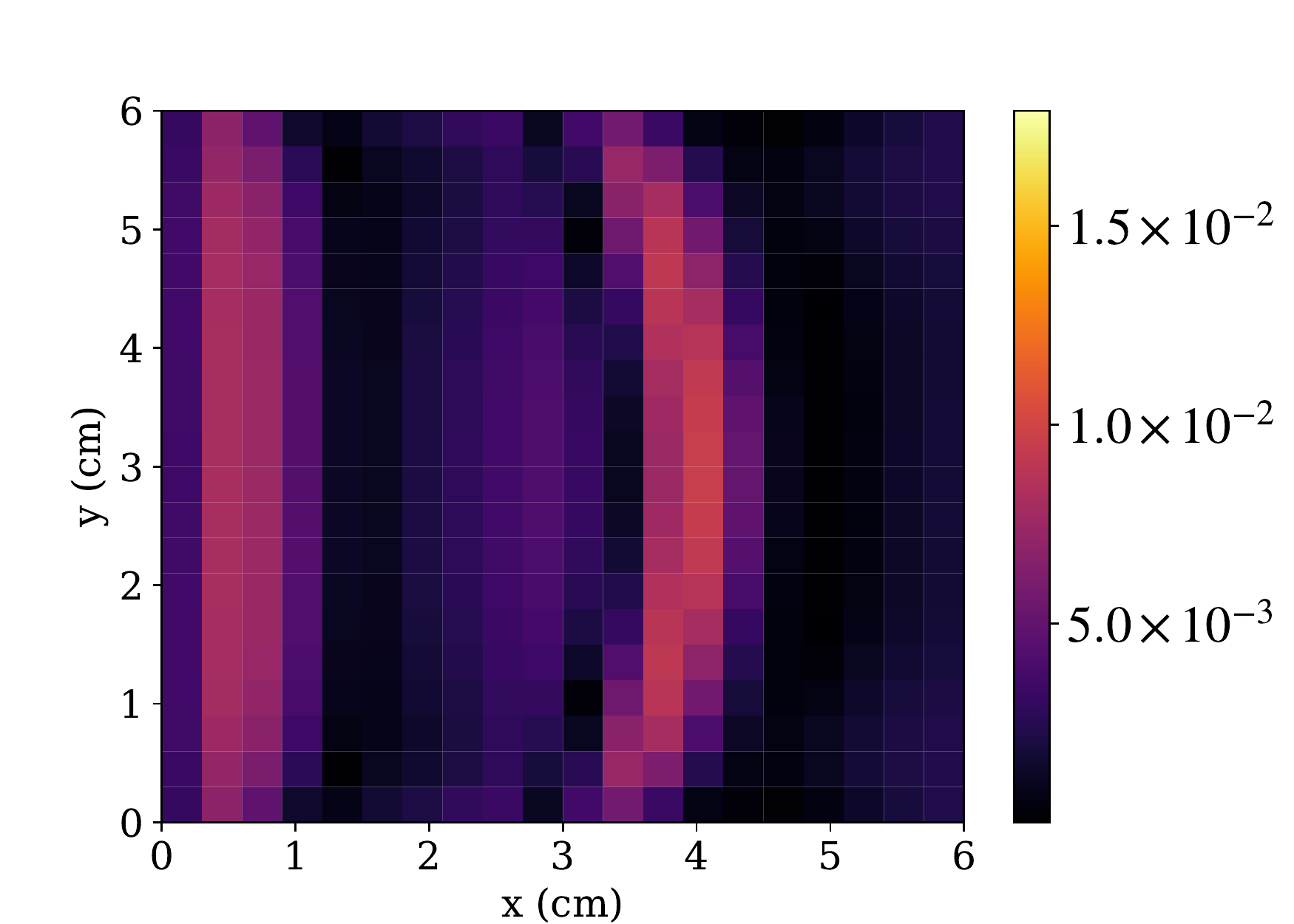}}&
			\raisebox{-.5\height}{\includegraphics[trim=1cm 0cm 4cm .5cm,clip,height=3.5cm]{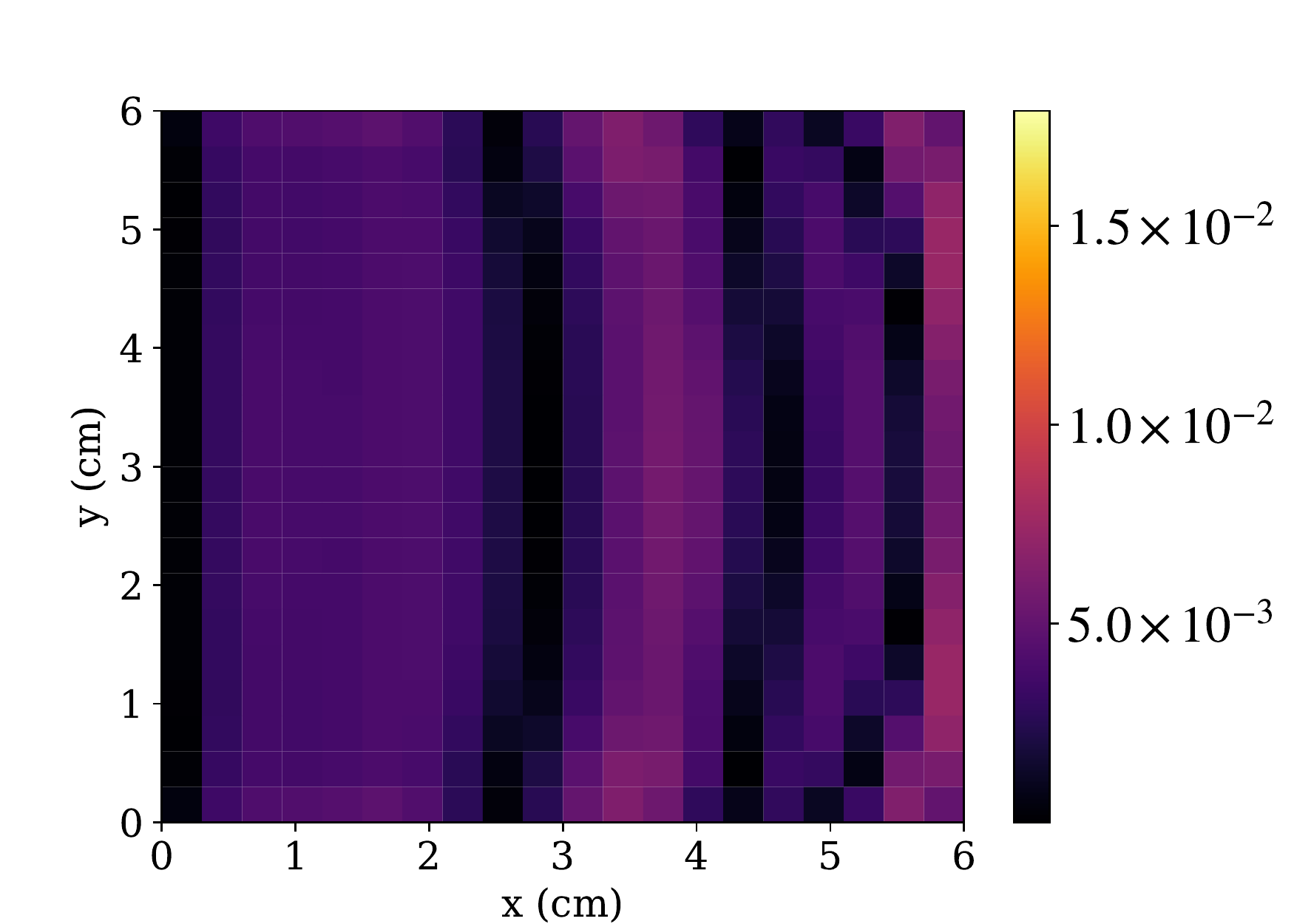}}&
			\raisebox{-.5\height}{\includegraphics[trim=1cm 0cm 0cm .5cm,clip,height=3.5cm]{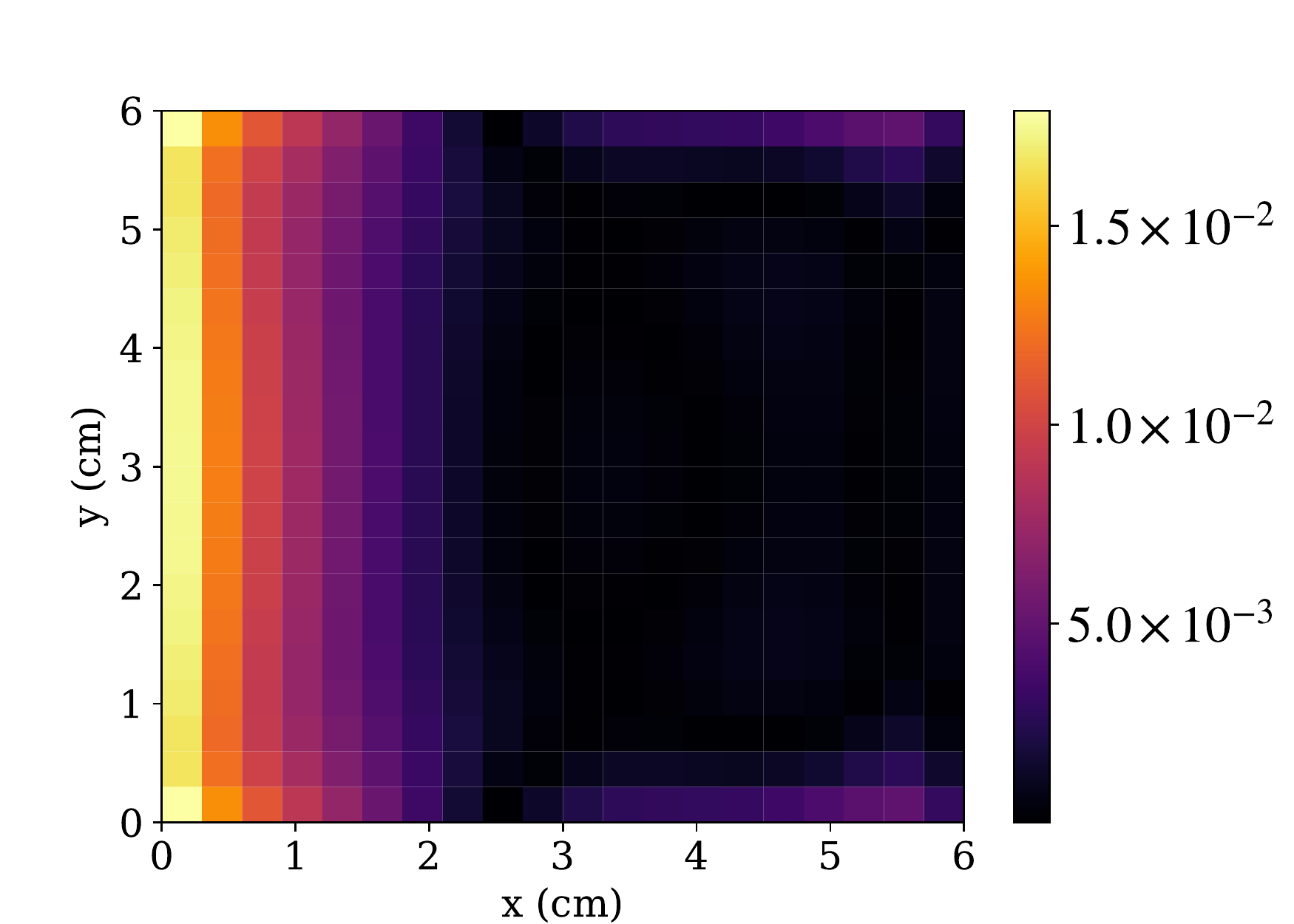}}\\[.1cm] \hline &&& \\[-.3cm]
			%%%
			$10^{-4}$ &
			\raisebox{-.5\height}{\includegraphics[trim=1cm 0cm 4cm .5cm,clip,height=3.5cm]{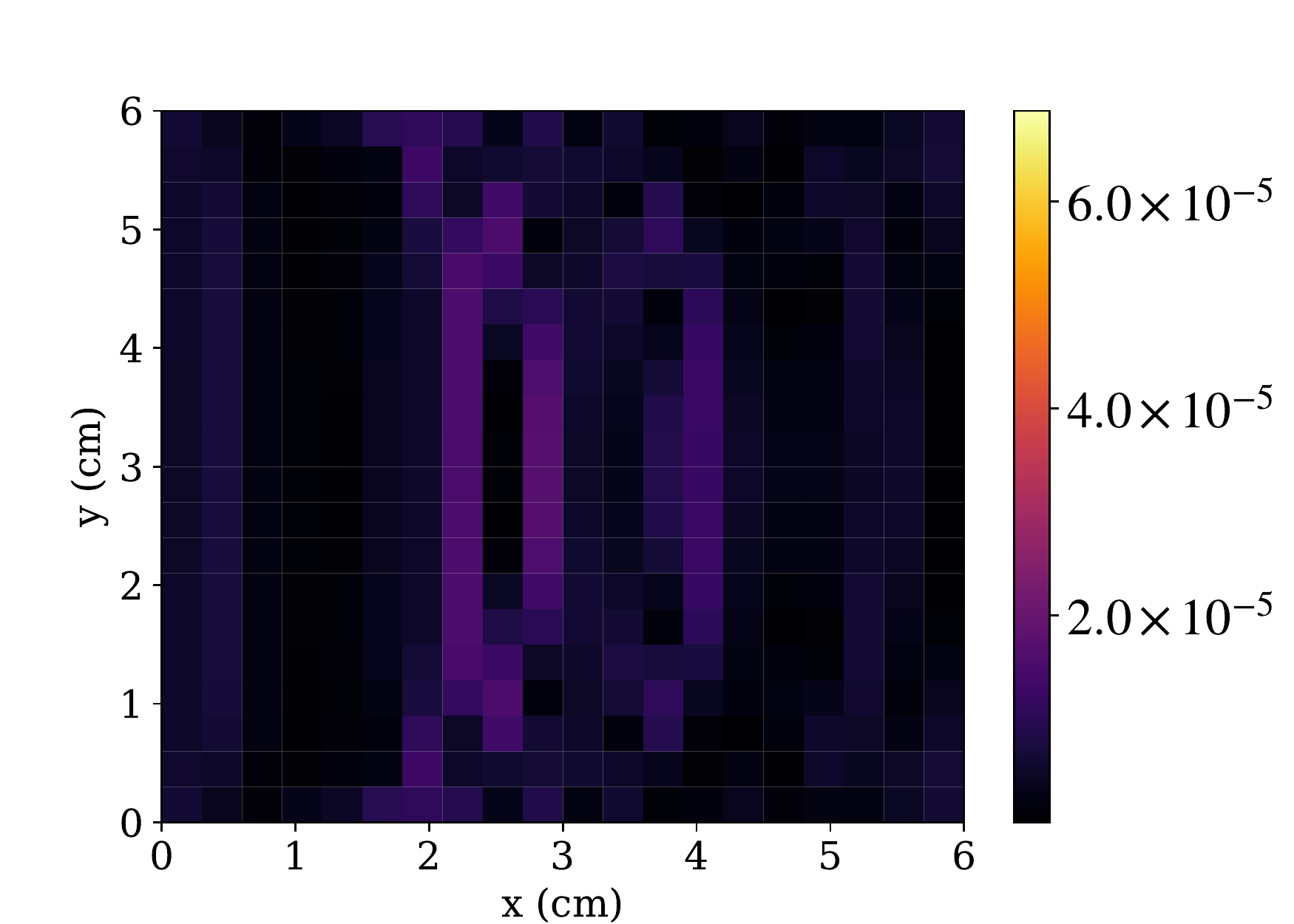}}&
			\raisebox{-.5\height}{\includegraphics[trim=1cm 0cm 4cm .5cm,clip,height=3.5cm]{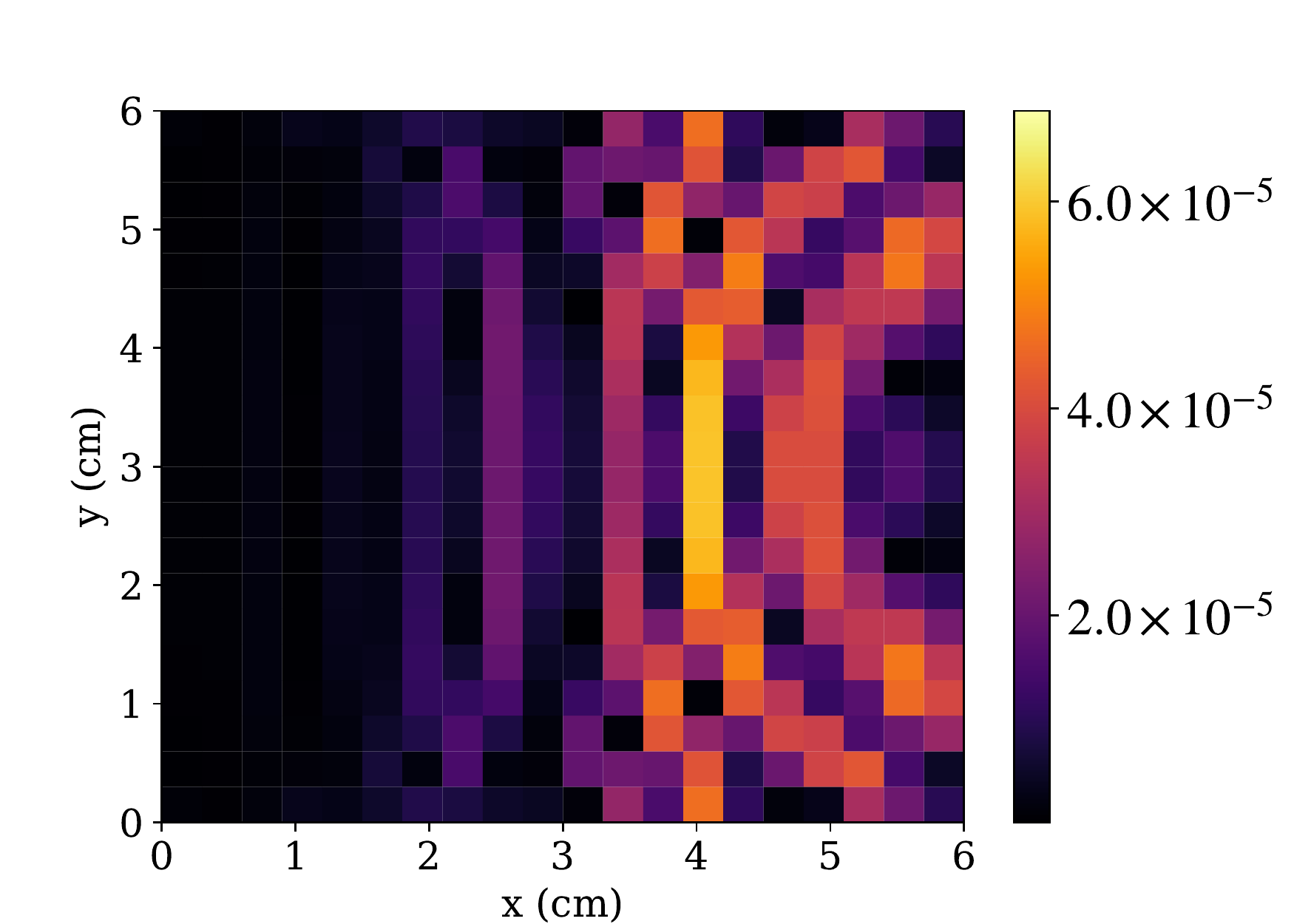}}&
			\raisebox{-.5\height}{\includegraphics[trim=1cm 0cm 0cm .5cm,clip,height=3.5cm]{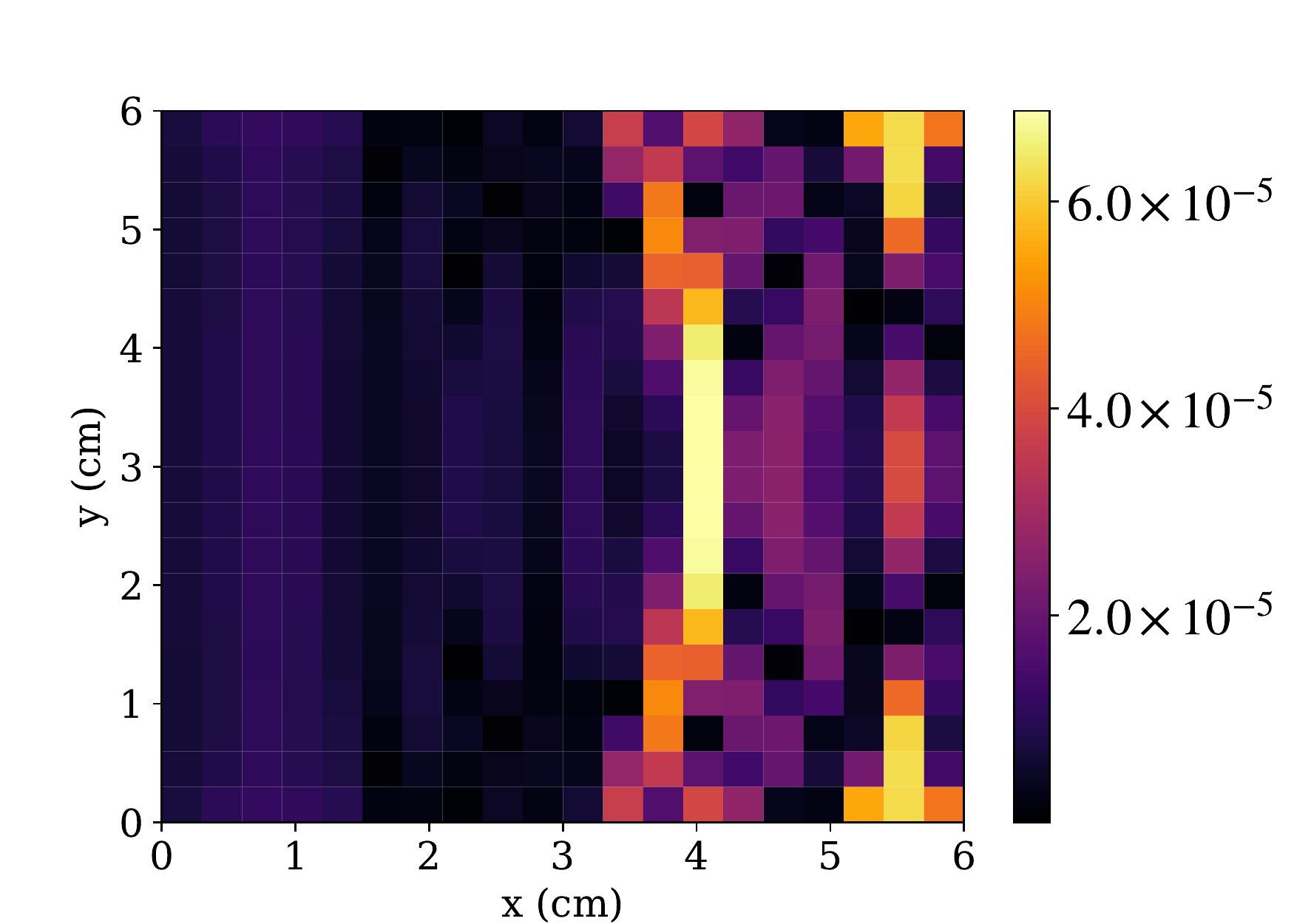}}
	\end{tabular} }
	\vspace*{.25cm}
	\subfloat[Radiation Energy Density]{
		\begin{tabular}{c|c|c|c}
			$\xi_\text{rel} $& t=1ns & t=2ns & t=3ns \\ \hline &&& \\[-.3cm]
			%%%
			%%%
			$10^{-2}$ &
			\raisebox{-.5\height}{\includegraphics[trim=1cm 0cm 4cm .5cm,clip,height=3.5cm]{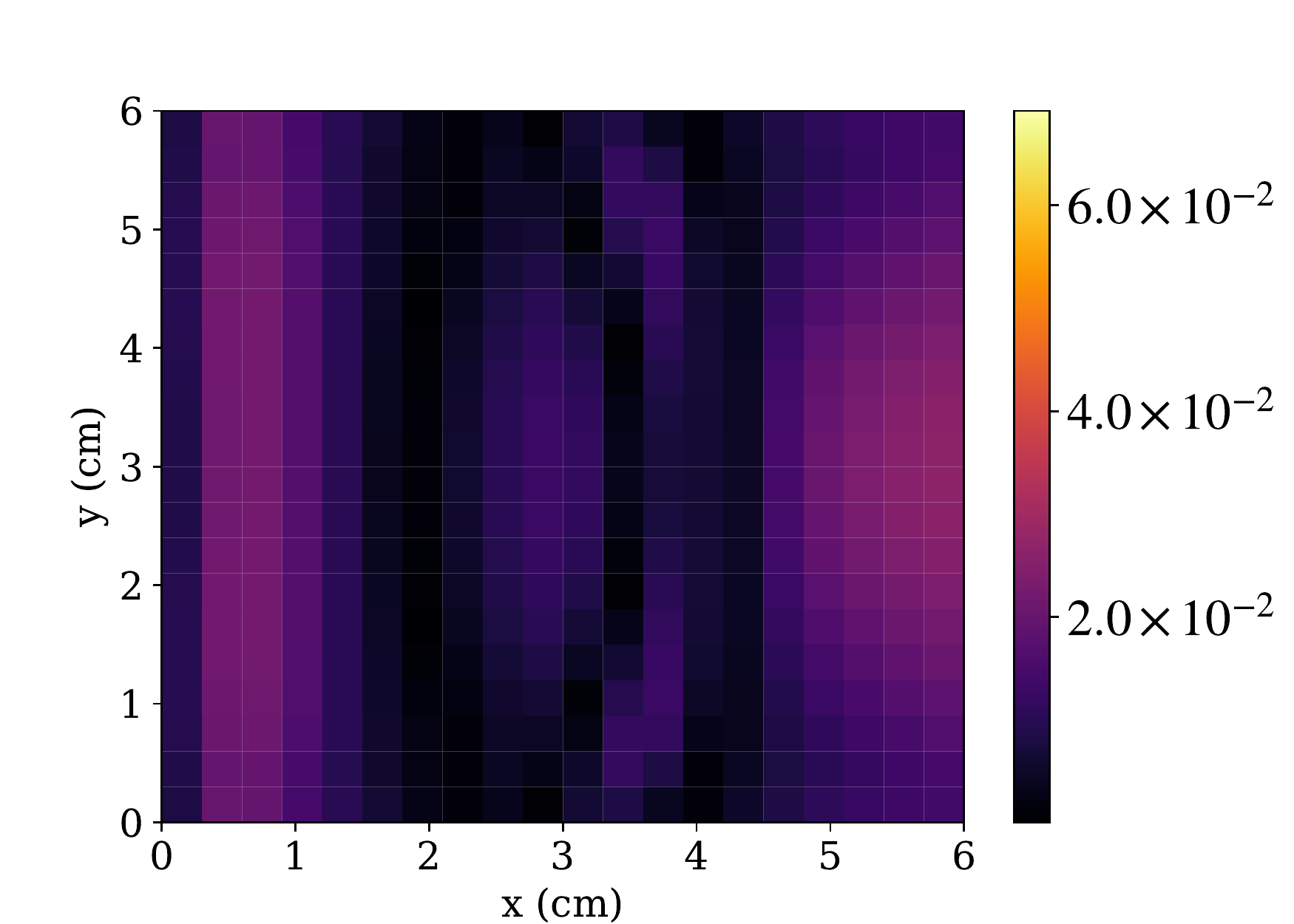}}&
			\raisebox{-.5\height}{\includegraphics[trim=1cm 0cm 4cm .5cm,clip,height=3.5cm]{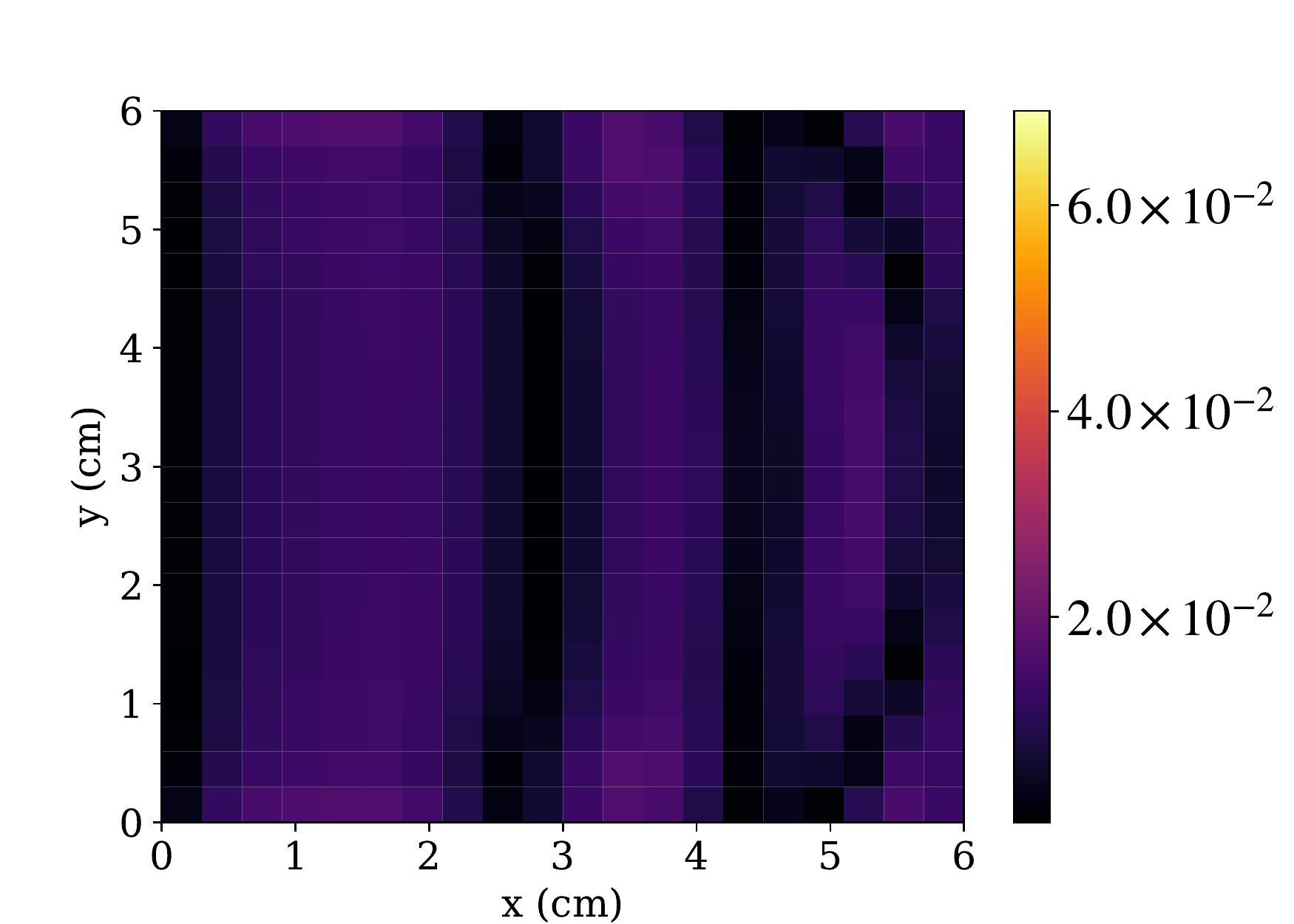}}&
			\raisebox{-.5\height}{\includegraphics[trim=1cm 0cm 0cm .5cm,clip,height=3.5cm]{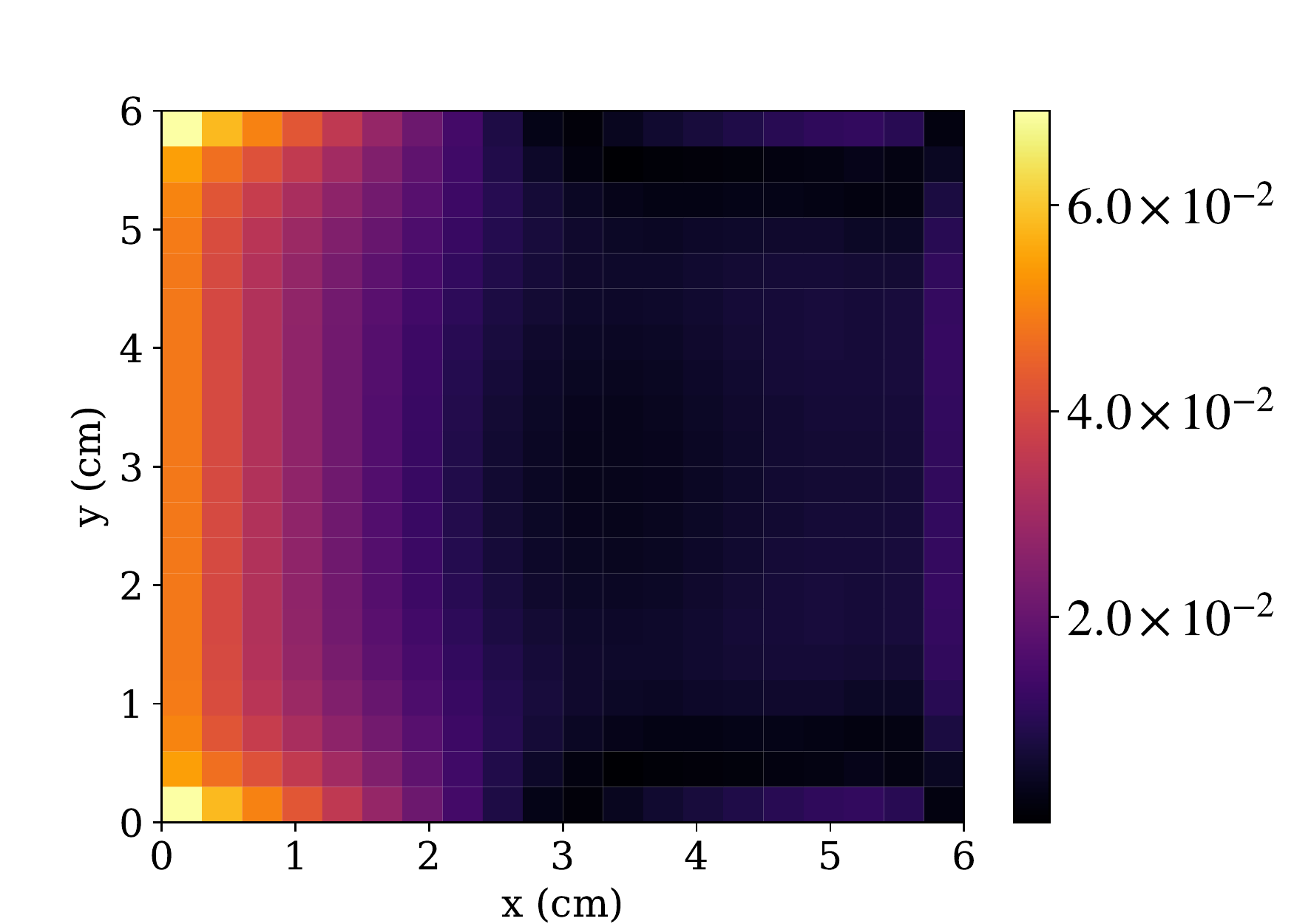}}\\[.1cm] \hline &&& \\[-.3cm]
			%%%
			$10^{-4}$ &
			\raisebox{-.5\height}{\includegraphics[trim=1cm 0cm 4cm .5cm,clip,height=3.5cm]{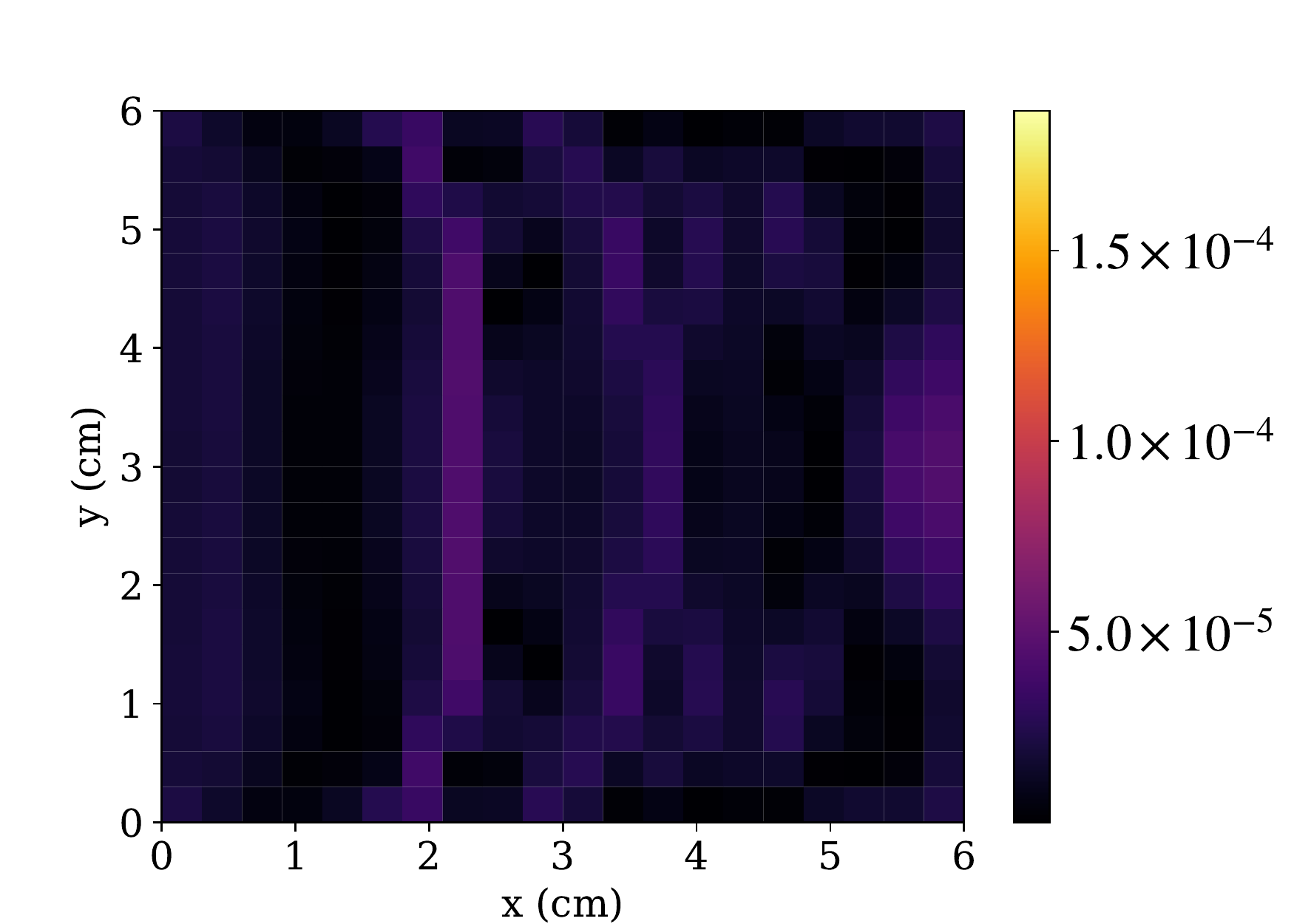}}&
			\raisebox{-.5\height}{\includegraphics[trim=1cm 0cm 4cm .5cm,clip,height=3.5cm]{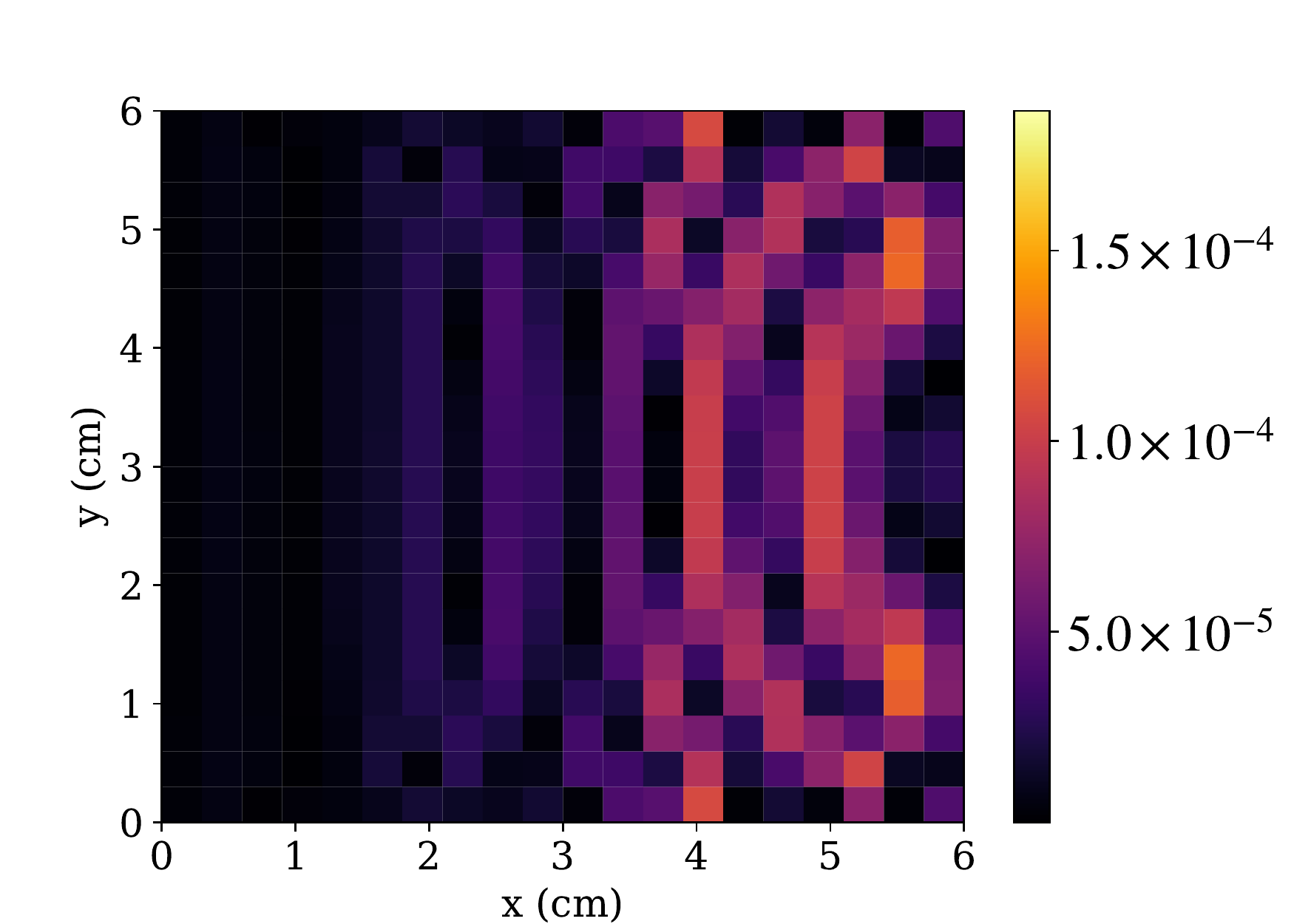}}&
			\raisebox{-.5\height}{\includegraphics[trim=1cm 0cm 0cm .5cm,clip,height=3.5cm]{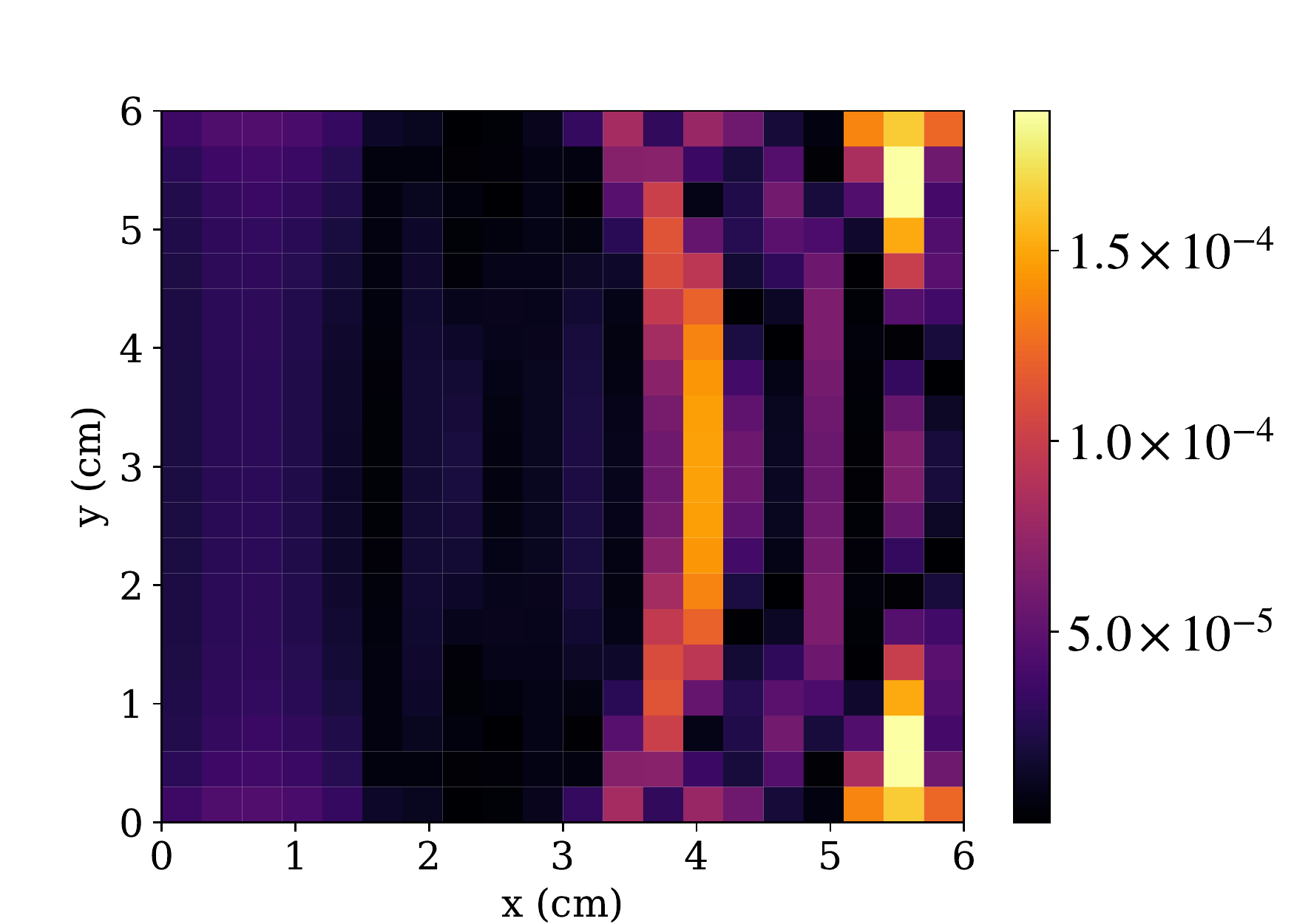}}
	\end{tabular} }
	\vspace*{.25cm}
	\caption{Cell-wise relative error in material temperature $(T)$ and total radiation energy density $(E)$ over the spatial domain at times  t=1, 2, 3 ns for the DDET ROM equipped with the DMD for $\xi_{\text{rel}}=10^{-2},10^{-4}$.}

	\label{fig:spatial_errs_dmd}
\end{figure}

\begin{figure}[ht!]
	\centering
	\subfloat[Material Temperature]{
		\begin{tabular}{c|c|c|c}
			$\xi_\text{rel} $& t=1ns & t=2ns & t=3ns \\ \hline &&& \\[-.3cm]
			%%%
			$10^{-2}$ &
			\raisebox{-.5\height}{\includegraphics[trim=1cm 0cm 4cm .5cm,clip,height=3.5cm]{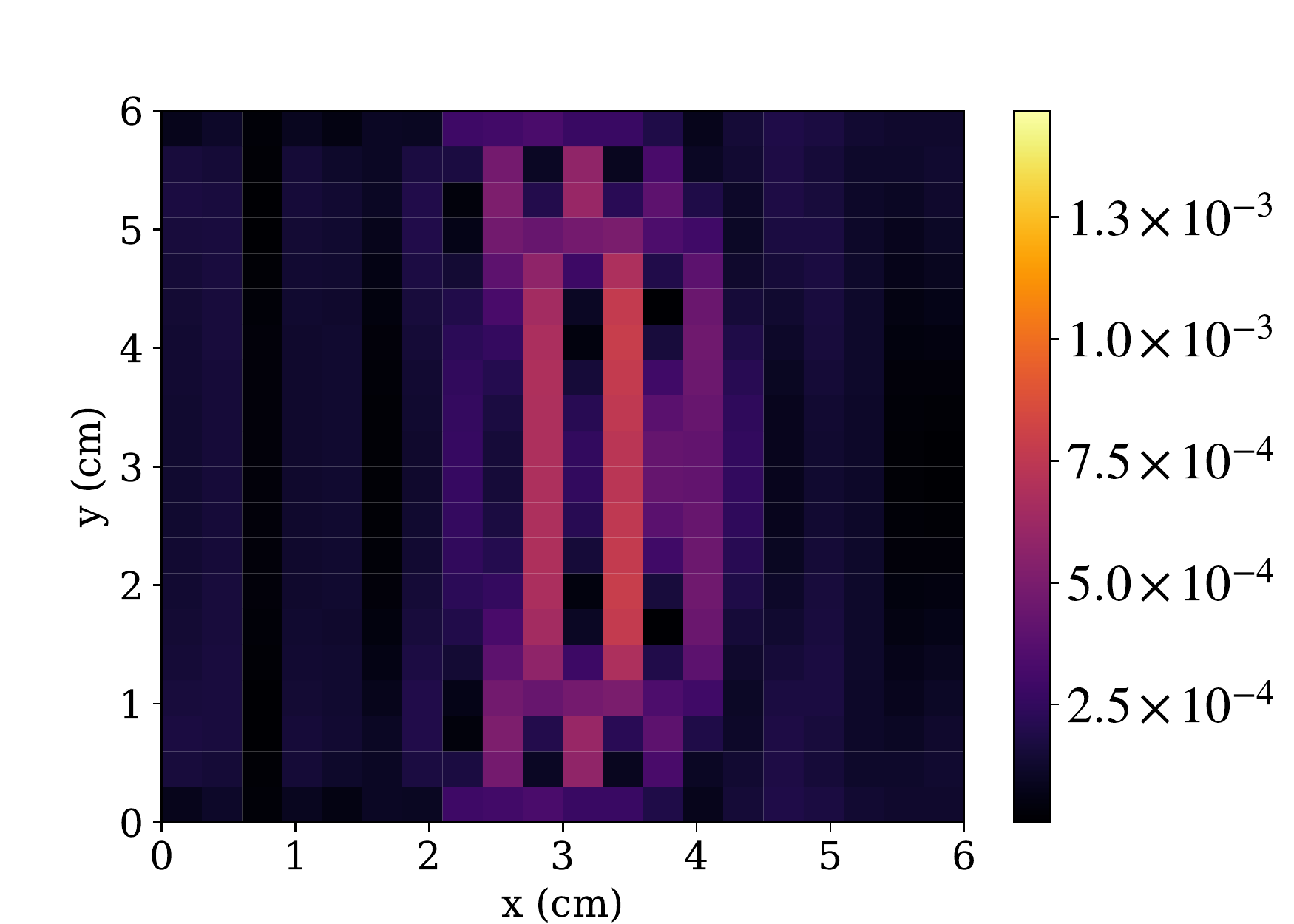}}&
			\raisebox{-.5\height}{\includegraphics[trim=1cm 0cm 4cm .5cm,clip,height=3.5cm]{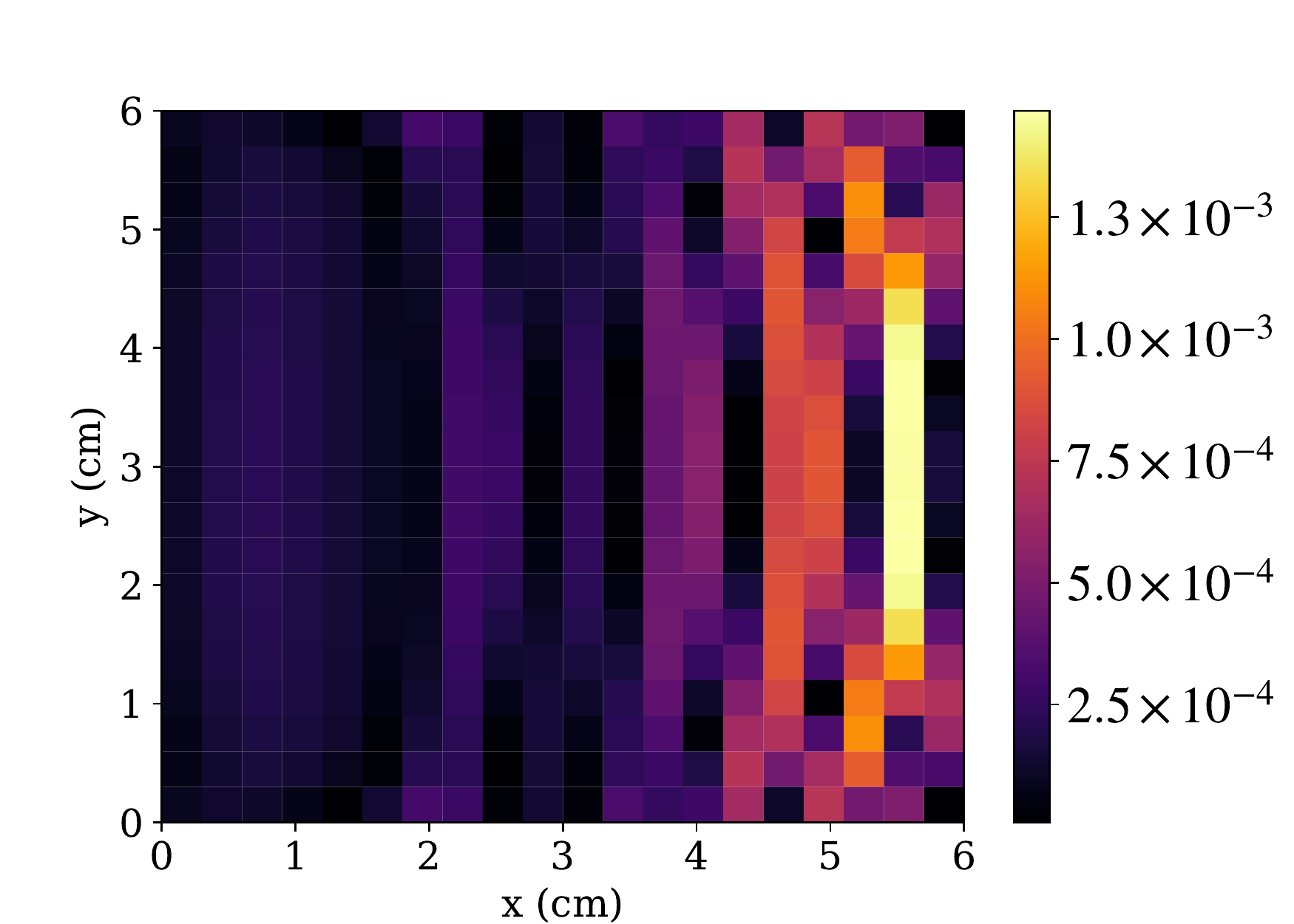}}&
			\raisebox{-.5\height}{\includegraphics[trim=1cm 0cm 0cm .5cm,clip,height=3.5cm]{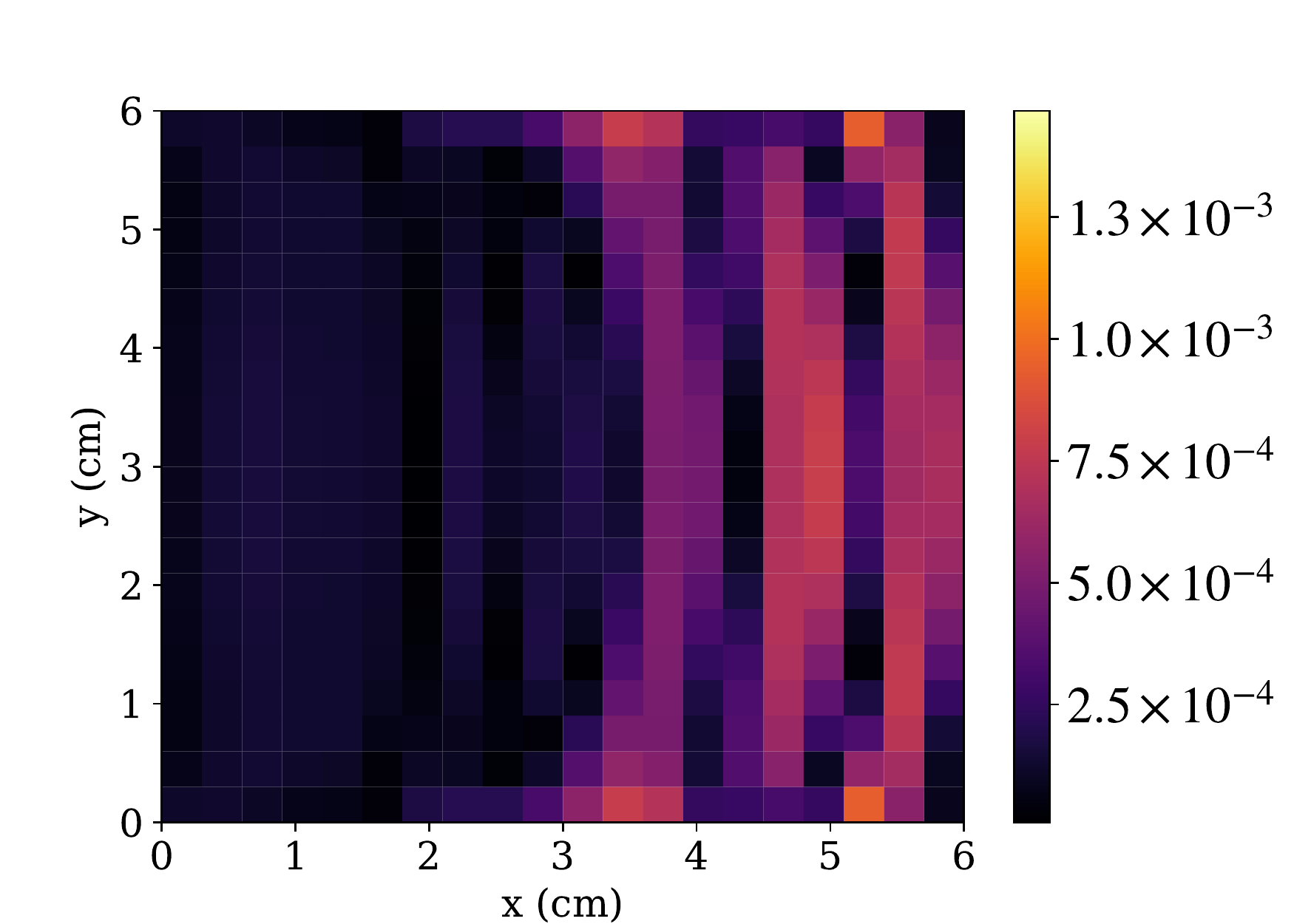}}\\[.1cm] \hline &&& \\[-.3cm]
			%%%
			$10^{-4}$ &
			\raisebox{-.5\height}{\includegraphics[trim=1cm 0cm 4cm .5cm,clip,height=3.5cm]{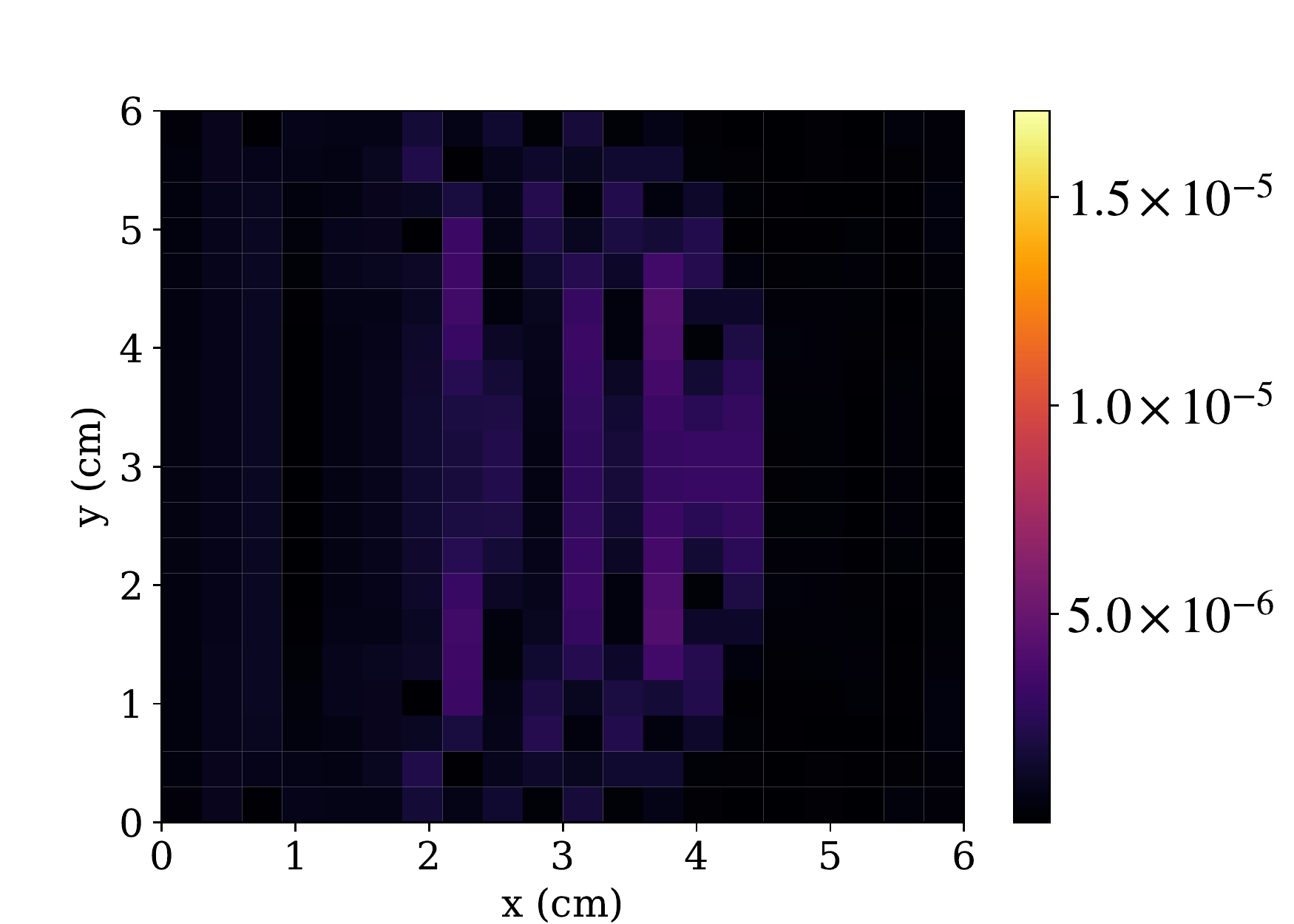}}&
			\raisebox{-.5\height}{\includegraphics[trim=1cm 0cm 4cm .5cm,clip,height=3.5cm]{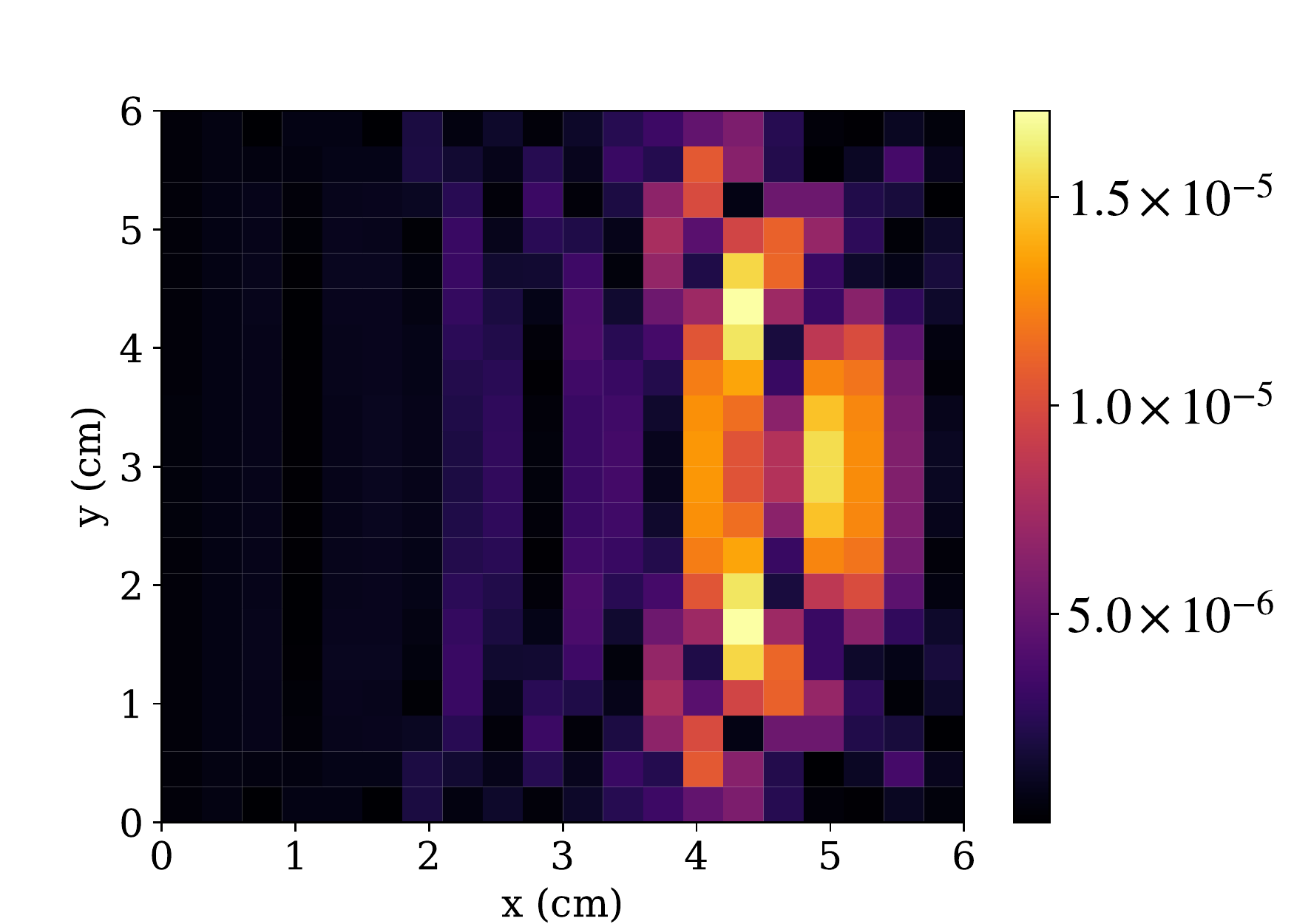}}&
			\raisebox{-.5\height}{\includegraphics[trim=1cm 0cm 0cm .5cm,clip,height=3.5cm]{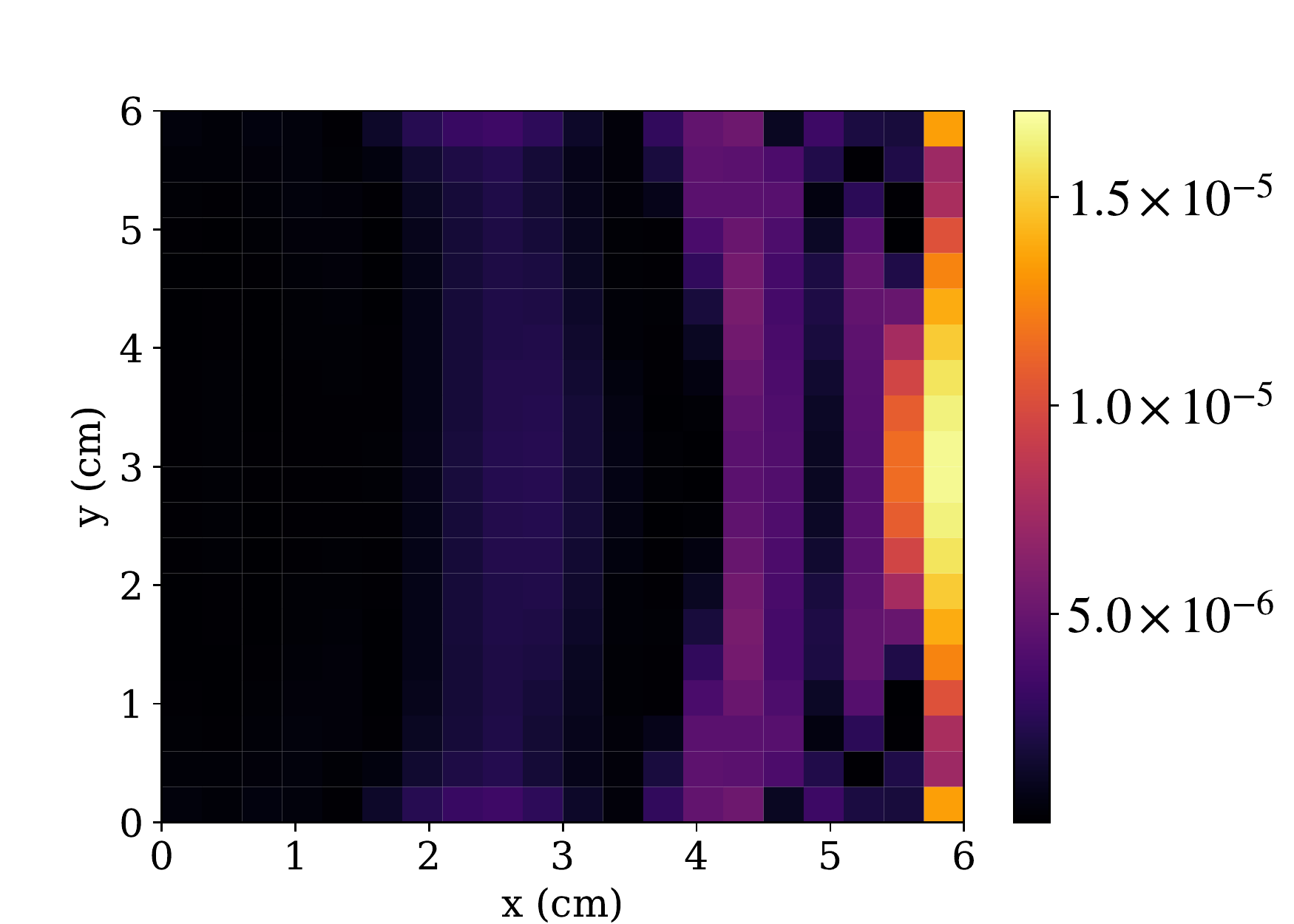}}
	\end{tabular} }
	\vspace*{.25cm}
	\subfloat[Radiation Energy Density]{
		\begin{tabular}{c|c|c|c}
			$\xi_\text{rel} $& t=1ns & t=2ns & t=3ns \\ \hline &&& \\[-.3cm]
			%%%
			%%%
			$10^{-2}$ &
			\raisebox{-.5\height}{\includegraphics[trim=1cm 0cm 4cm .5cm,clip,height=3.5cm]{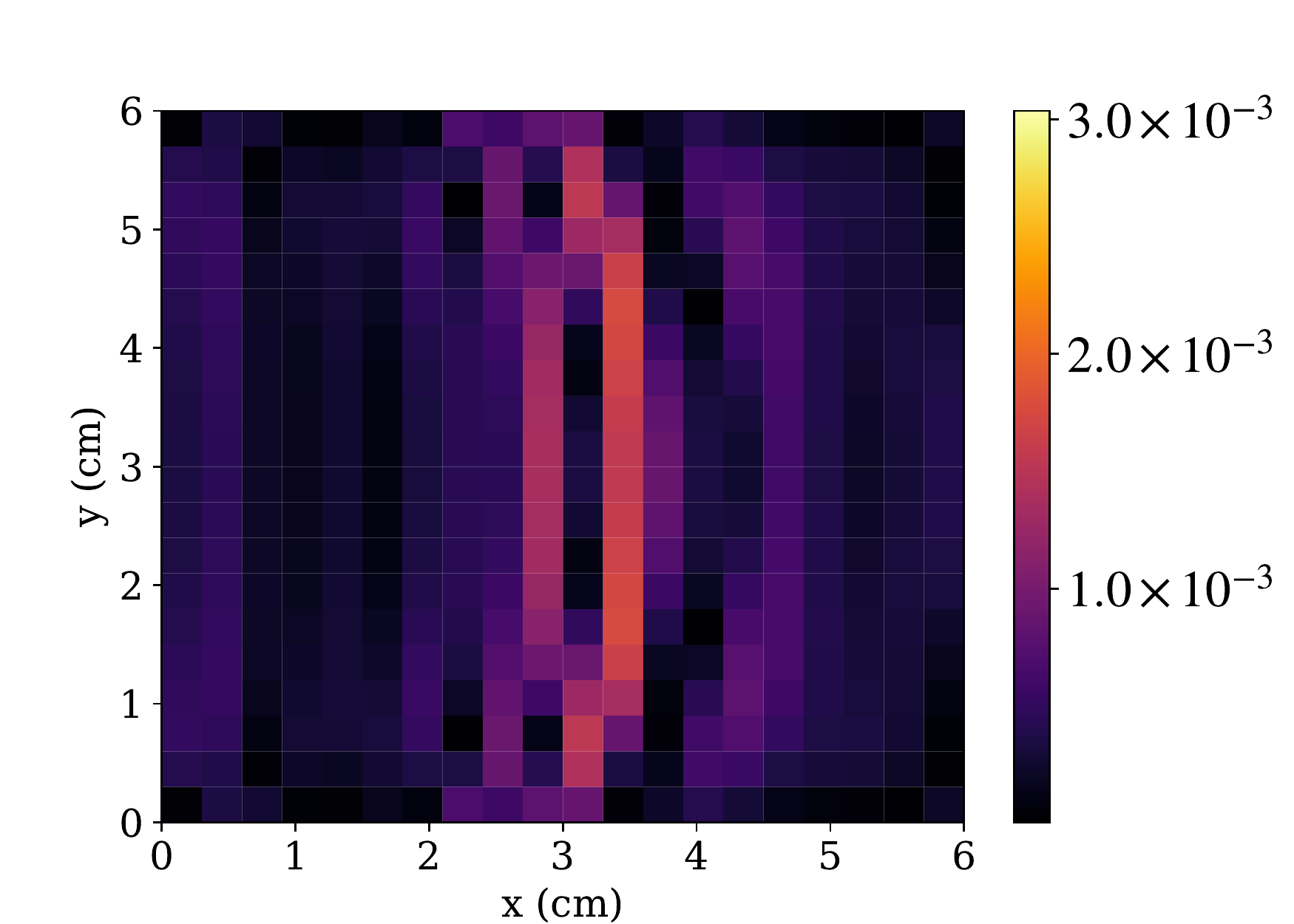}}&
			\raisebox{-.5\height}{\includegraphics[trim=1cm 0cm 4cm .5cm,clip,height=3.5cm]{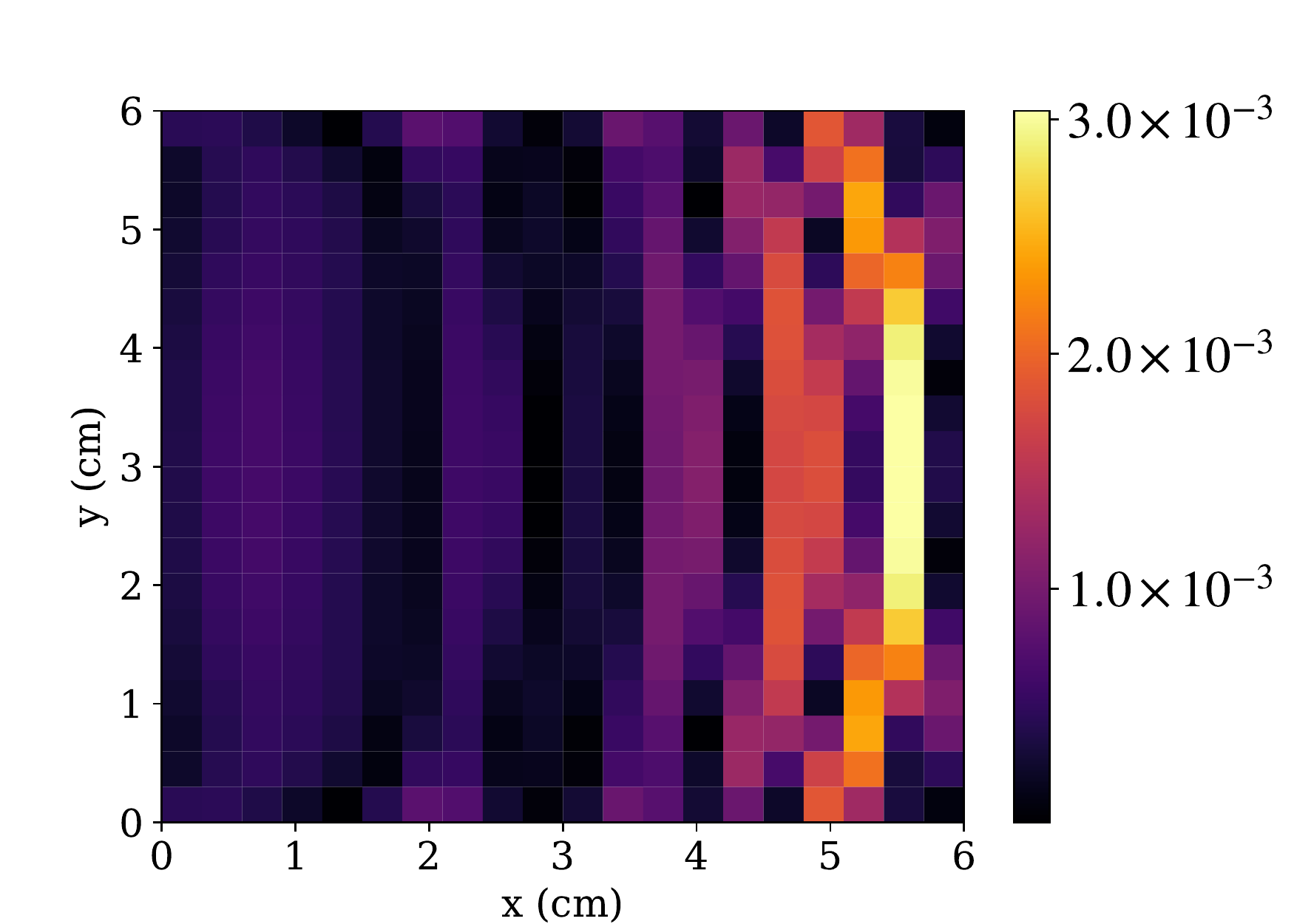}}&
			\raisebox{-.5\height}{\includegraphics[trim=1cm 0cm 0cm .5cm,clip,height=3.5cm]{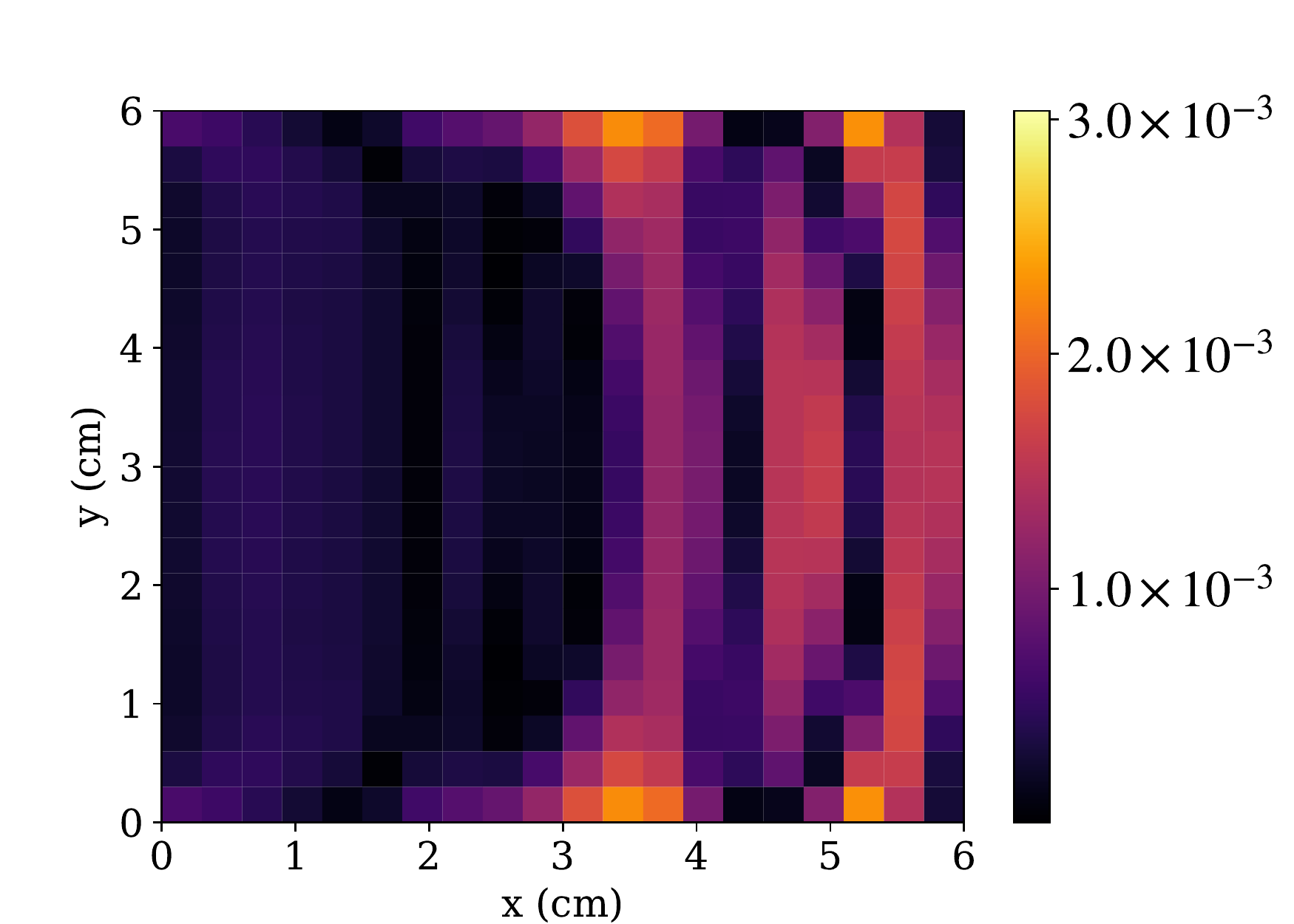}}\\[.1cm] \hline &&& \\[-.3cm]
			%%%
			$10^{-4}$ &
			\raisebox{-.5\height}{\includegraphics[trim=1cm 0cm 4cm .5cm,clip,height=3.5cm]{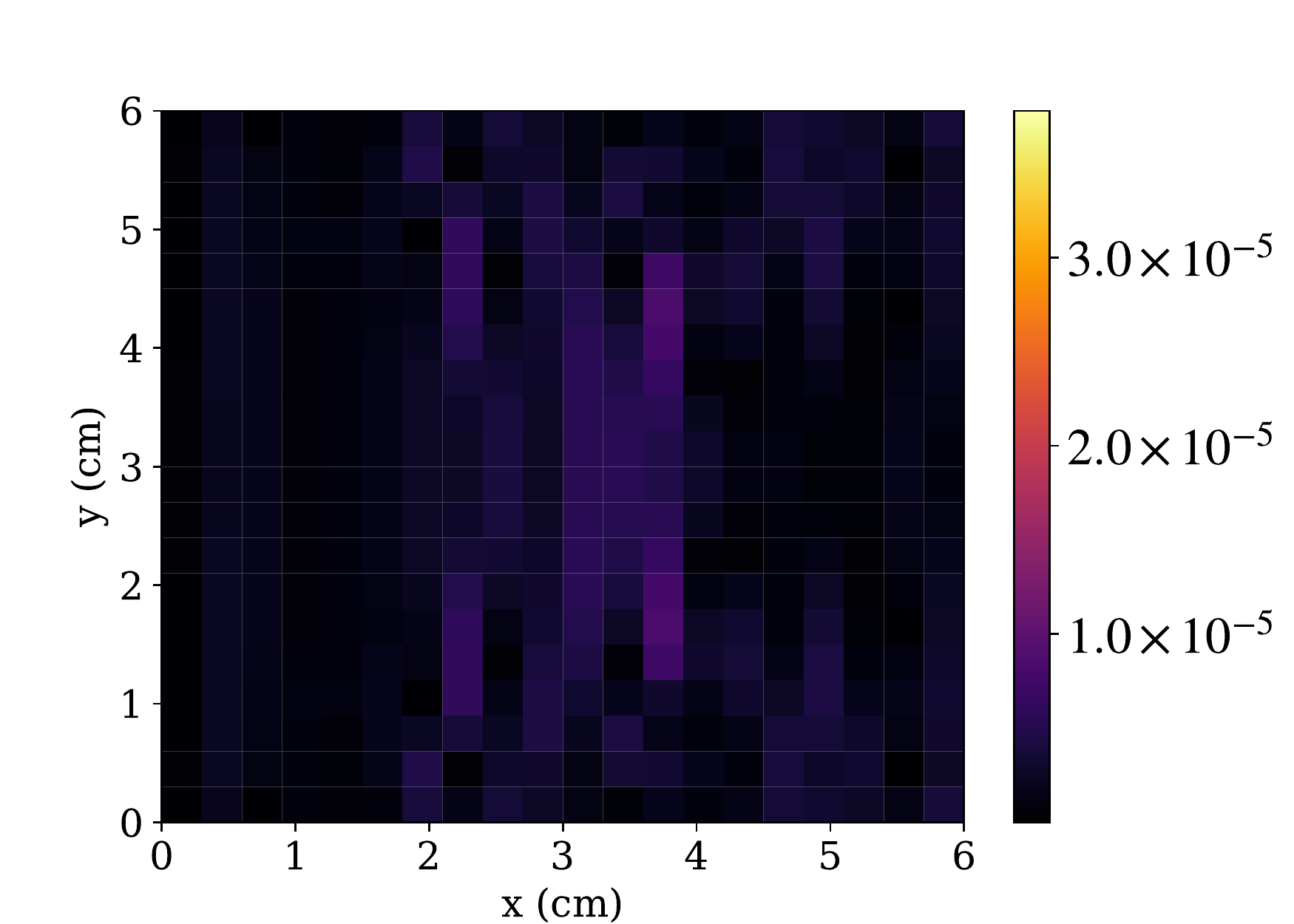}}&
			\raisebox{-.5\height}{\includegraphics[trim=1cm 0cm 4cm .5cm,clip,height=3.5cm]{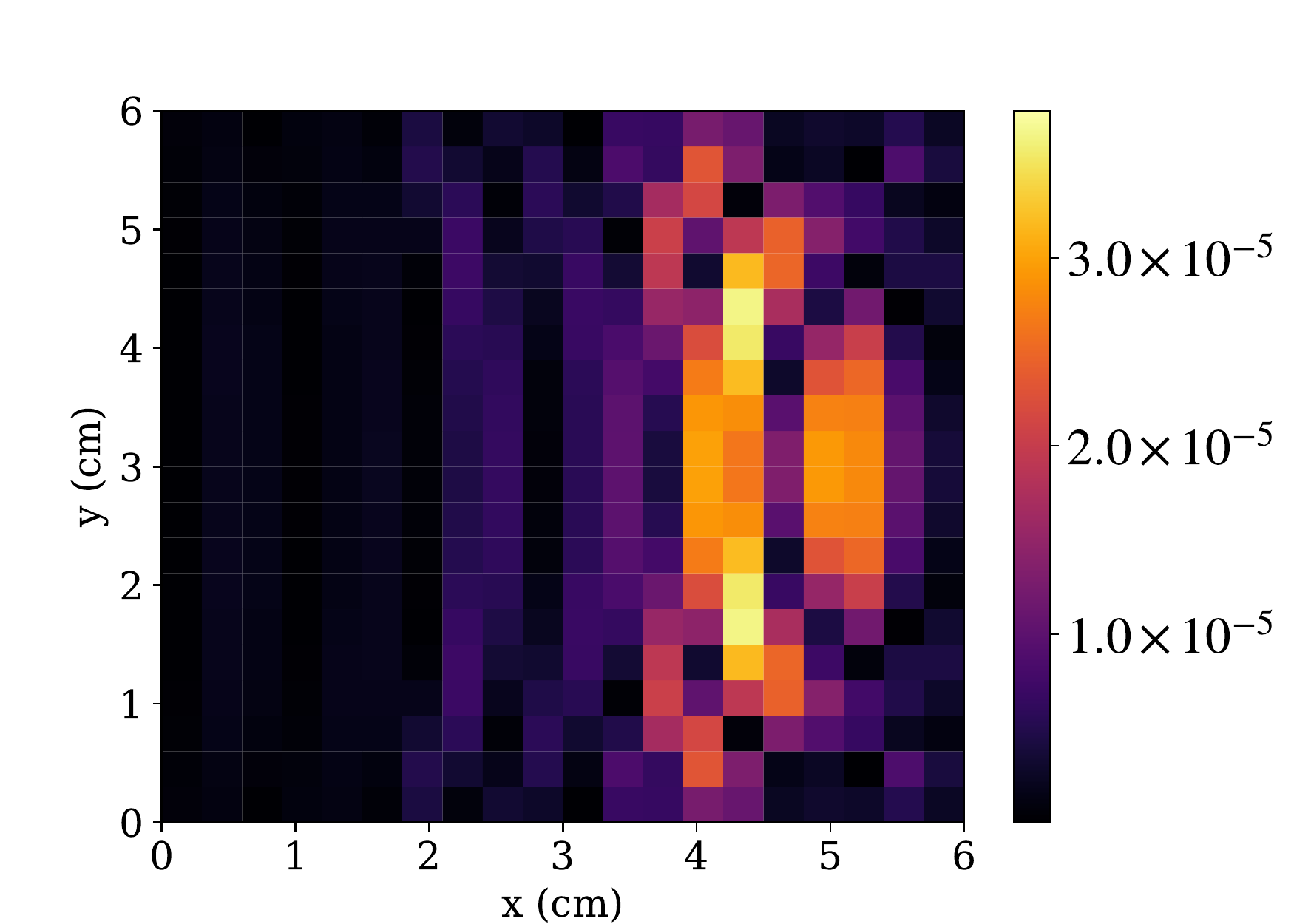}}&
			\raisebox{-.5\height}{\includegraphics[trim=1cm 0cm 0cm .5cm,clip,height=3.5cm]{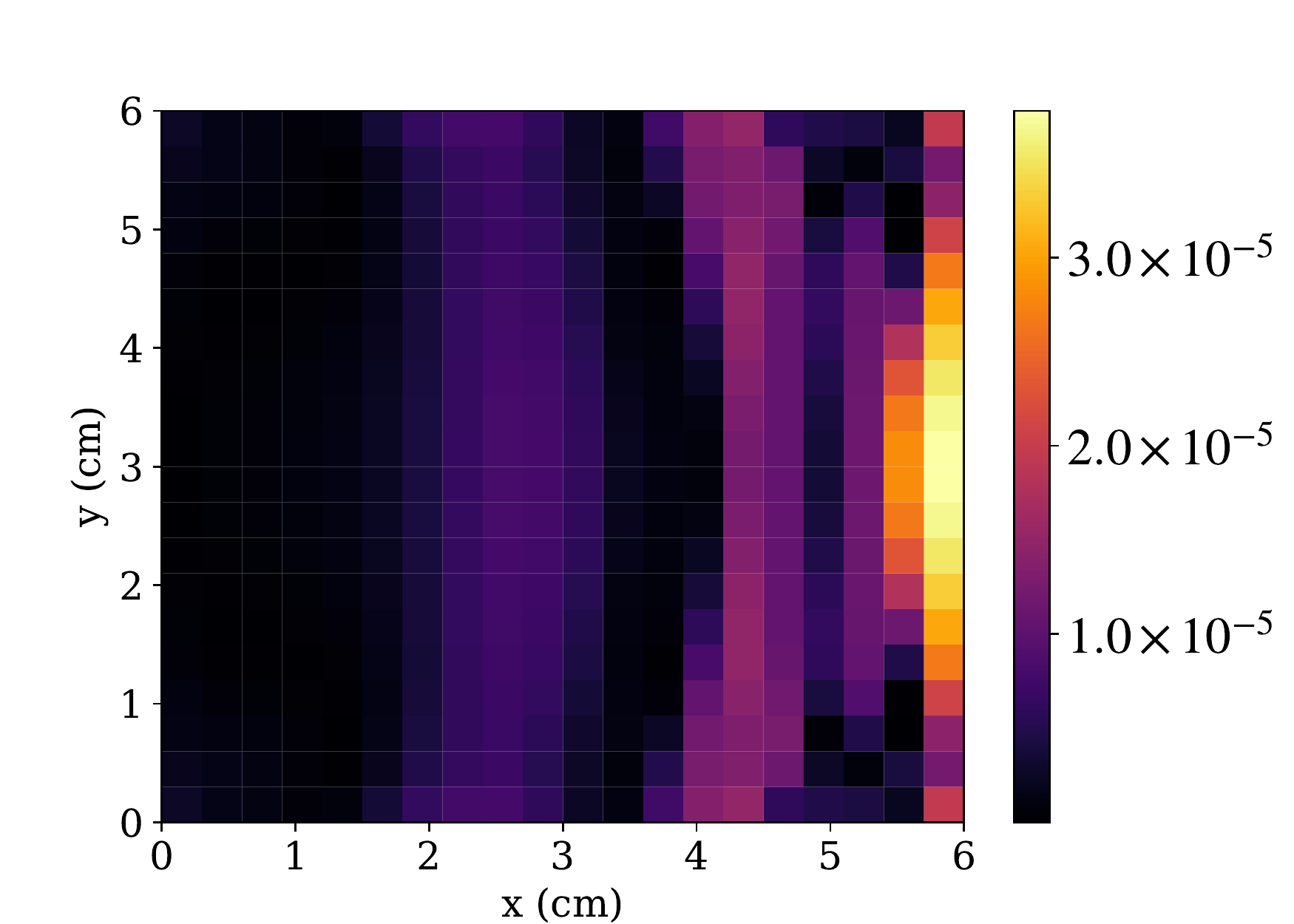}}
	\end{tabular} }
	\vspace*{.25cm}
	\caption{Cell-wise relative error in material temperature $(T)$ and total radiation energy density $(E)$ over the spatial domain at times t=1, 2, 3 ns for the DDET ROM equipped with the DMD-E for $\xi_{\text{rel}}=10^{-2},10^{-4}$.}

	\label{fig:spatial_errs_dmdb}
\end{figure}

In practice the ROMs with lowest-rank can become the most important as \linebreak undoubtedly the DDET ROM reaps the most computational benefits with larger $\xi_{\text{rel}}$ \linebreak(i.e. $\xi_{\text{rel}}=10^{-2},10^{-4}$).
In this way the properties and accuracy of the DDET ROMs with large $\xi_{\text{rel}}$ are important to understand.
The results shown so far indicate that the DDET ROM performs well with very low-rank representations of the Eddington tensor given by the POD, DMD and DMD-E.
The spatial errors have been studied the 2-norm.  We now analyse local behavior of errors over the spatial domain.
Figures \ref{fig:spatial_errs_pod}, \ref{fig:spatial_errs_dmd} and \ref{fig:spatial_errs_dmdb} show
cell-wise relative errors in $T$ and $E$ at selected instances.
Each of these figures takes the form of two tables that display the relative pointwise error in the DDET ROM across the spatial domain of the F-C test. The first (top) table shows errors in the material temperature $(T)$ and the second (bottom) shows errors in the total radiation energy density $(E)$. Each row corresponds to a different value of $\xi_\text{rel}$ and each column corresponds to the specific instant of time. These include $\xi_\text{rel}=10^{-2},10^{-4}$ and $t=1,2,3$ ns, respectively. Furthermore all plots contained on a single row use the same scale for their color distributions.
In order, Figures \ref{fig:spatial_errs_pod}, \ref{fig:spatial_errs_dmd} and \ref{fig:spatial_errs_dmdb} correspond to the errors in the DDET ROMs equipped with the POD, DMD and DMD-E. These figures clearly demonstrate that the spatial distribution of errors in those low-rank ROMs is relatively uniform. There are no sharp changes in the error about spatial position and each point has an error value residing in a close neighborhood to the relative 2-norm error for the corresponding ROM and time point shown in Figures \ref{fig:2nrm-errs_POD}, \ref{fig:2nrm-errs_DMD} and \ref{fig:2nrm-errs_DMDB}.


%=================================================================================
%
%=================================================================================
\clearpage
\subsection{Breakout Time}
\begin{figure}[ht!]
	\centering
	\subfloat[ $\bar{F}_R$]{\includegraphics[width=.33\textwidth]{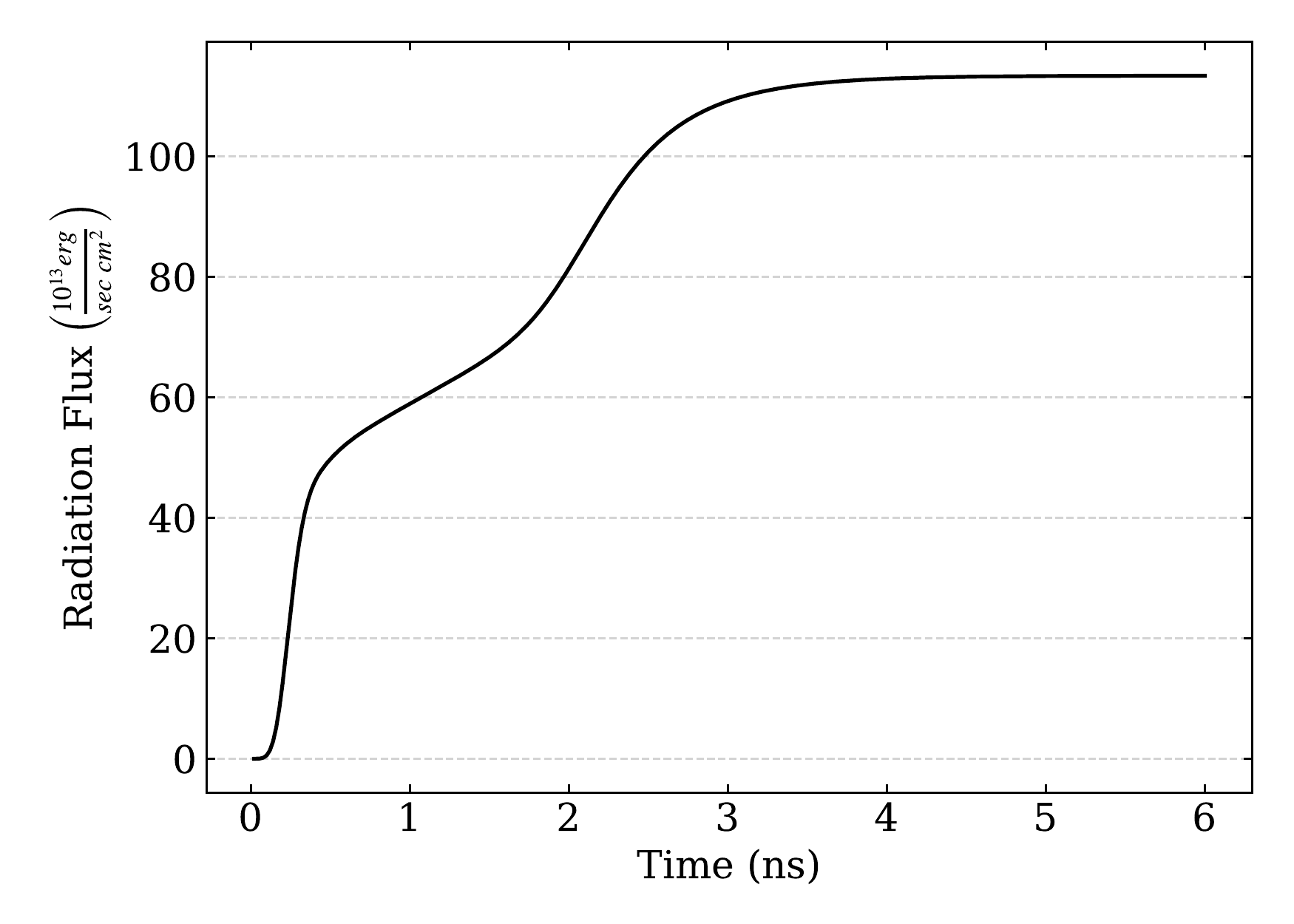}}
	\subfloat[$\bar{E}_R$]{\includegraphics[width=.33\textwidth]{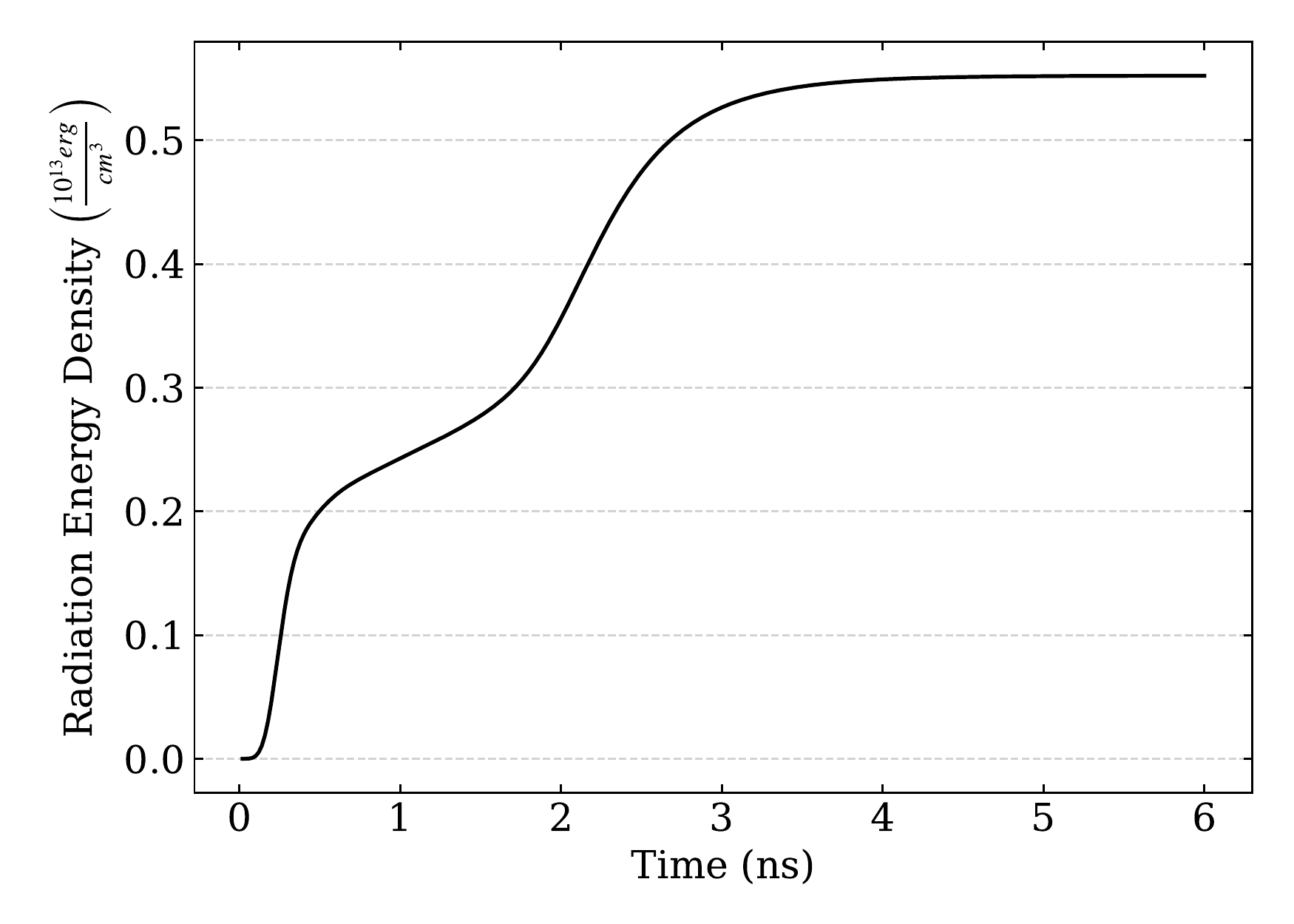}}
	\subfloat[$\bar{T}_R$]{\includegraphics[width=.33\textwidth]{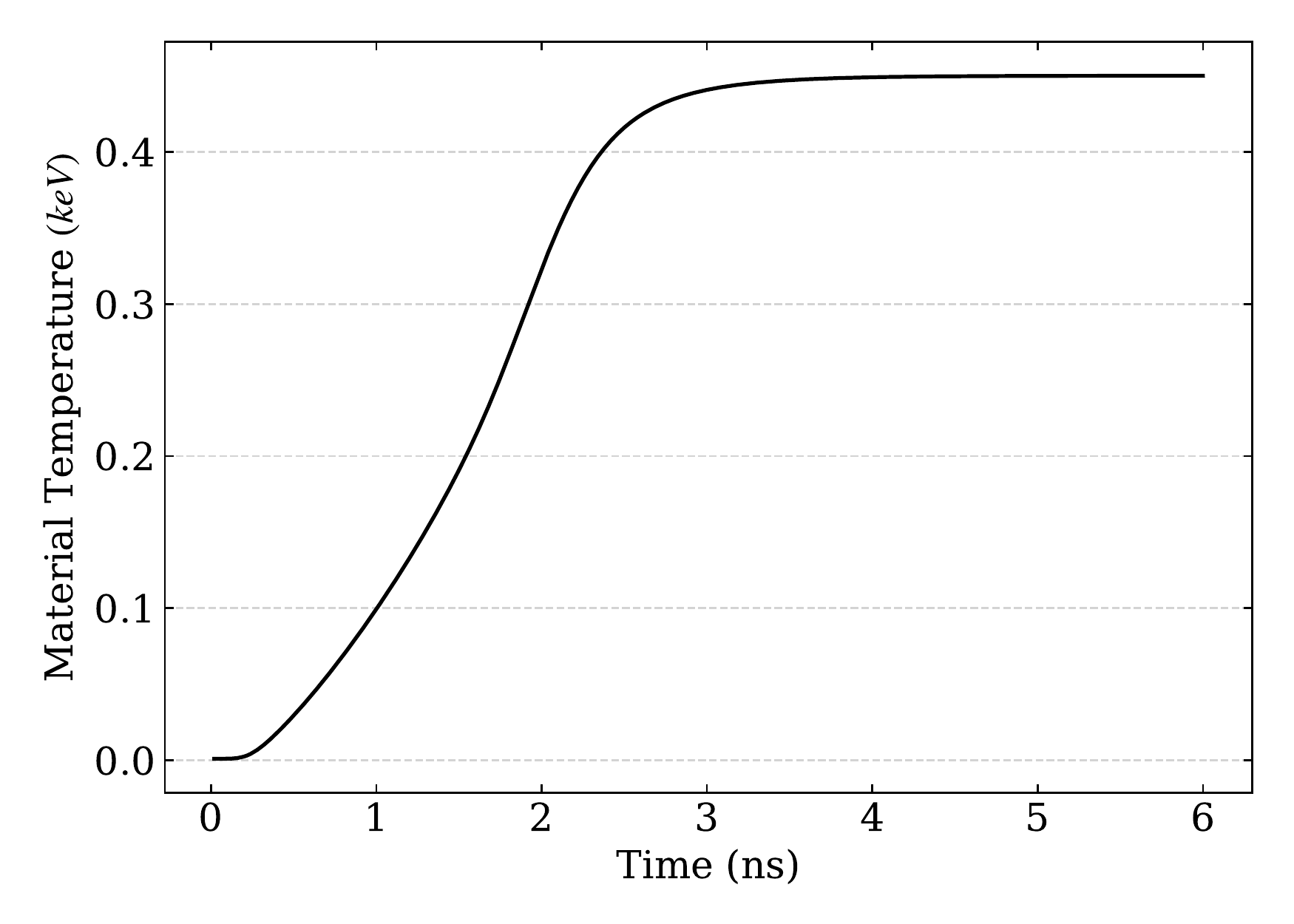}}
	\caption{Total radiation flux $(\bar{F}_R)$, total radiation energy density $(\bar E_R)$ and material temperature $(\bar T_R)$
		averaged over the right boundary of the spatial domain plotted vs time. Shown solutions are generated by the FOM.
		\label{fig:rbndvals_fom} }
\end{figure}

We now consider properties of the DDET class of ROMs in  capturing certain physics of TRT problems.
One metric of particular importance is {\it breakout time} of radiation that characterises how well the ROMs are able to reproduce the FOM radiation wavefront as
it propagates through the spatial domain \cite{fryer-2016,Moore-2015}.
The radiation wave produced in the F-C test travels from left to right and correspondingly the notion of breakout time is associated with radiation levels at the right boundary. Typically in the literature, breakout time is measured as the elapsed time until a certain level of radiative flux is detected \cite{fryer-2016,Moore-2015}.
Here we consider not only the radiation flux, but the energy density and material temperature at the right boundary of the F-C test as well.
We consider  the boundary-averages of these quantities, defined as follows:
\begin{subequations}
\begin{gather}
	\bar{F}_R = \frac{1}{L_R}\int_{0}^{L_R} \vec{e}_x\cdot\vec{F}(x_R,y)\ dy,\\[5pt] \bar{E}_R = \frac{1}{L_R}\int_{0}^{L_R} E(x_R,y)\ dy,\\[5pt] \bar{T}_R = \frac{1}{L_R}\int_{0}^{L_R} T(x_R,y)\ dy,
\end{gather}
\end{subequations}

\noindent where $L_R = x_R = 6\text{cm}$. The time evolution of $\bar{F}_R$, $\bar{E}_R$, and $\bar{T}_R$, calculated with the FOM is depicted in Figure \ref{fig:rbndvals_fom}. These figures show two sharp increases in $\bar{F}_R$ and $\bar{E}_R$ followed by plateaus, whereas $\bar{T}_R$  increases smoothly until reaching a final plateau. The initial plateaus for $\bar{F}_R$, $\bar{E}_R$ occur at roughly 0.5 ns and indicate when the high-energy radiation has penetrated the domain. The final plateaus for each $\bar{F}_R$, $\bar{E}_R$, and $\bar{T}_R$ occurs at about $2.5$ ns, indicating full penetration of the domain.

\begin{figure}[ht!]
	\subfloat[$\bar{F}_R$]{\includegraphics[width=.33\textwidth]{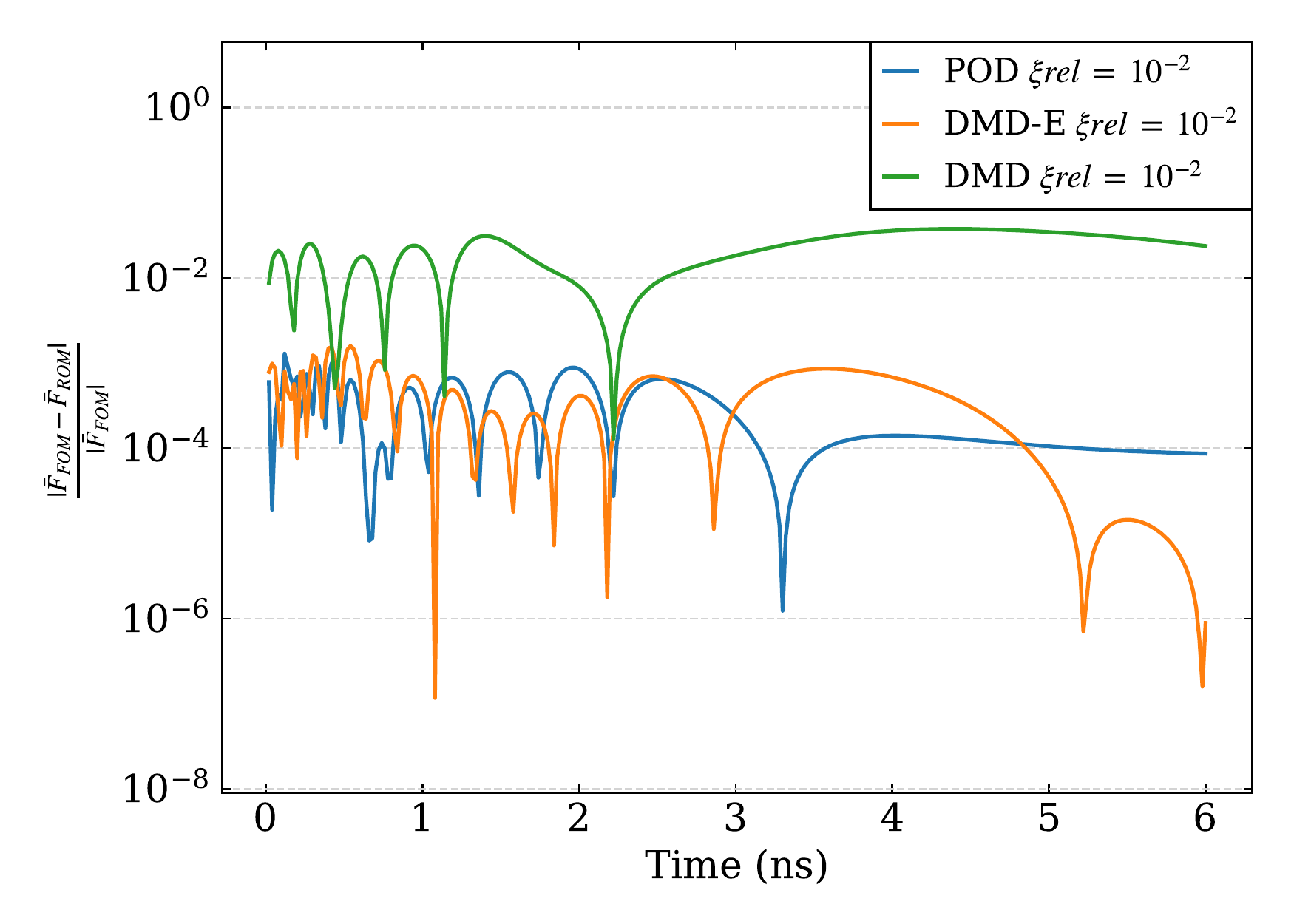}}
	\subfloat[$\bar E_R$]{\includegraphics[width=.33\textwidth]{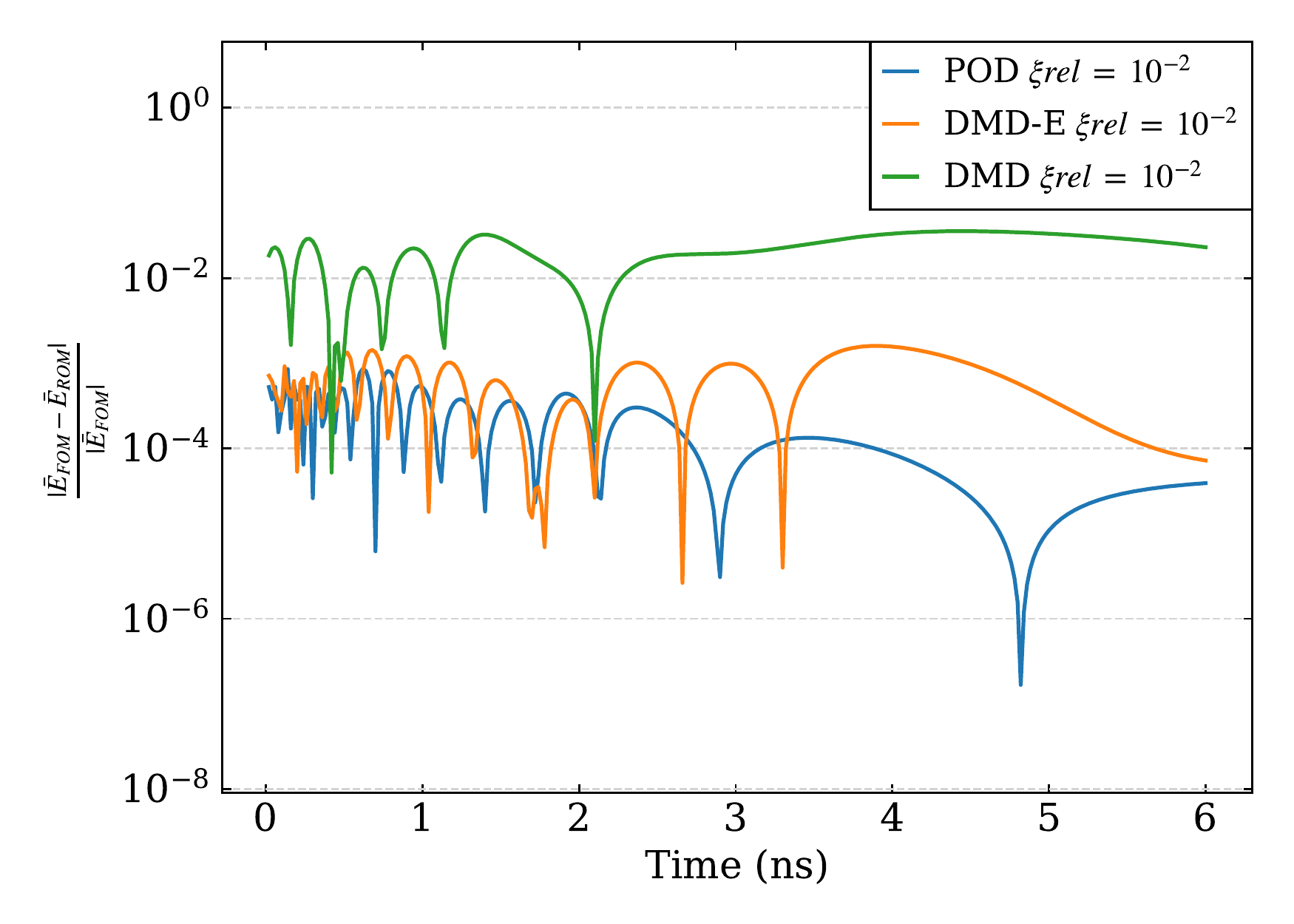}}
	\subfloat[$\bar T_R$]{\includegraphics[width=.33\textwidth]{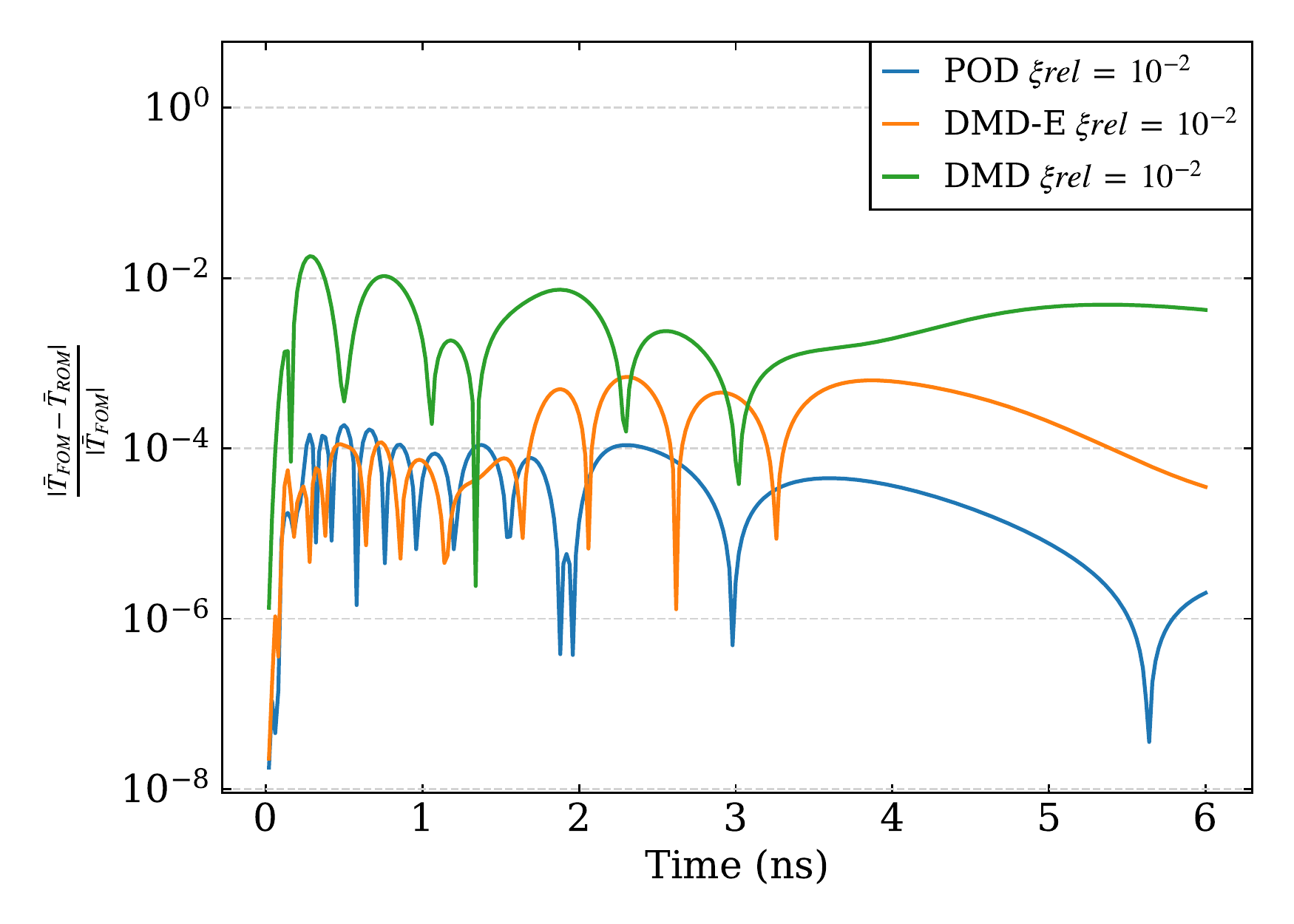}}
	\caption{Relative error  for the DDET ROMs with $\xi_{\text{rel}}=10^{-2}$
 for data located at and integrated over the right boundary of the domain.
		\label{fig:rbndvals_1e-2} }
	\end{figure}

\begin{figure}[ht!]
	\subfloat[$\bar{F}_R$]{\includegraphics[width=.33\textwidth]{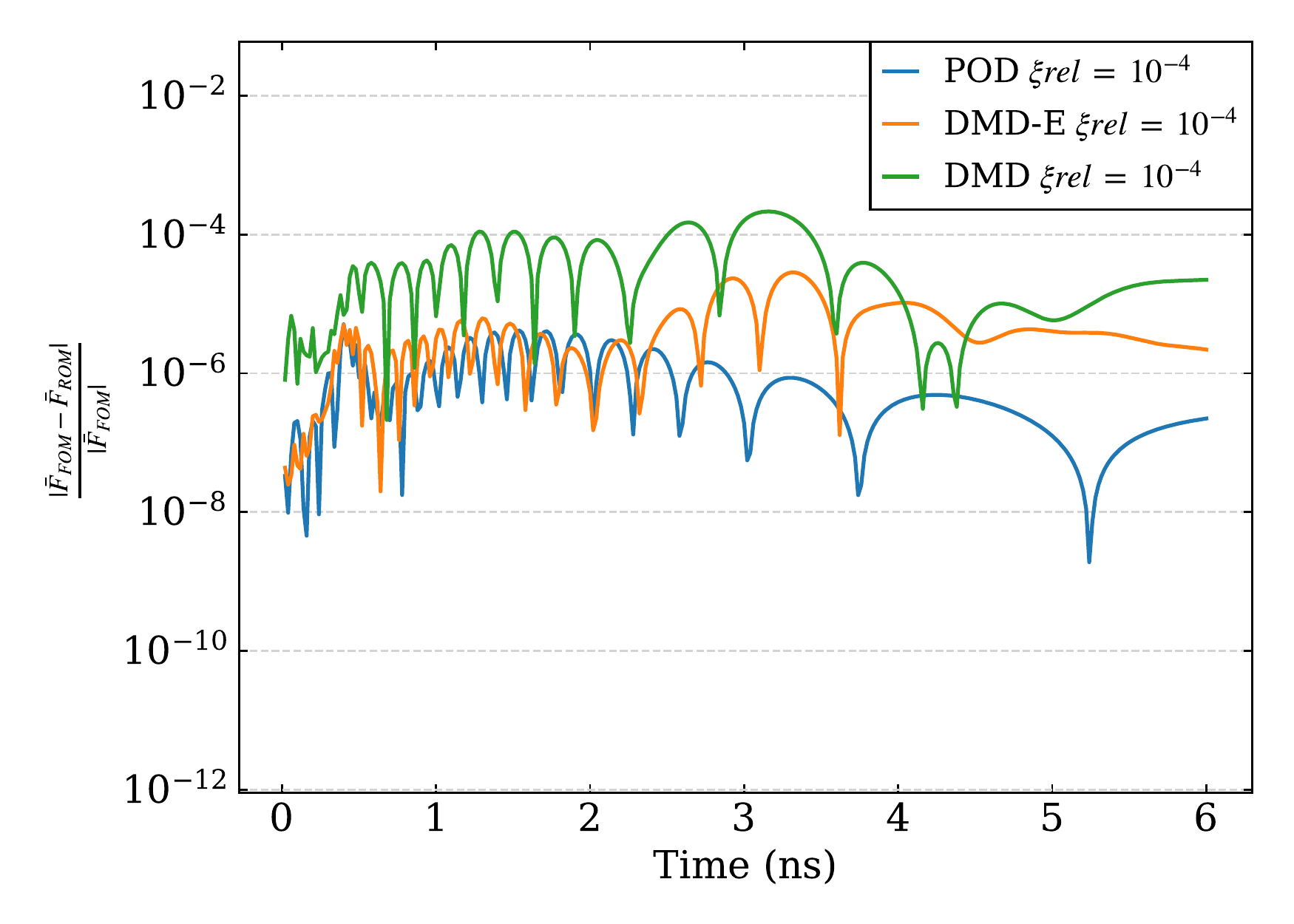}}
	\subfloat[$\bar E_R$]{\includegraphics[width=.33\textwidth]{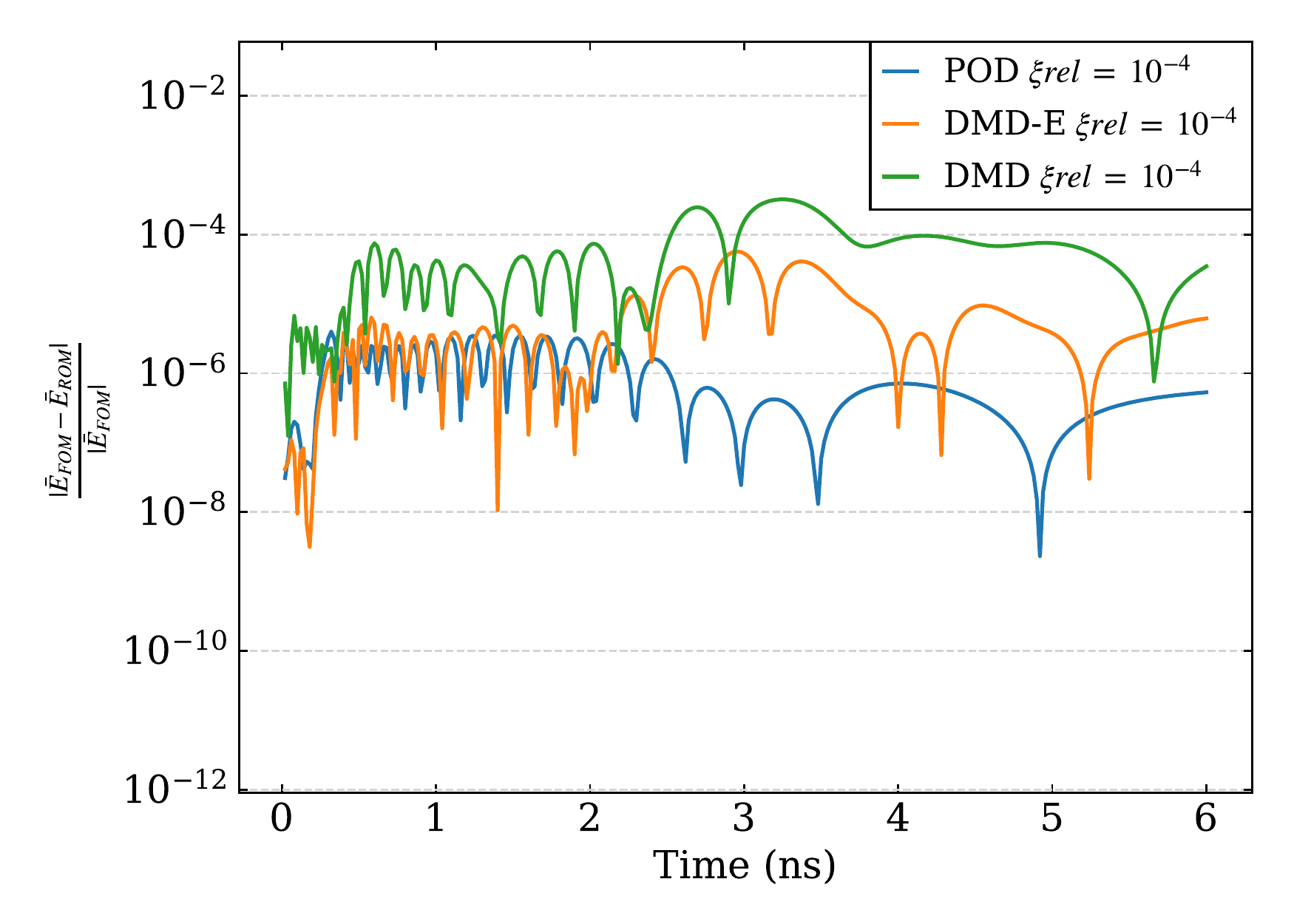}}
	\subfloat[$\bar T_R$]{\includegraphics[width=.33\textwidth]{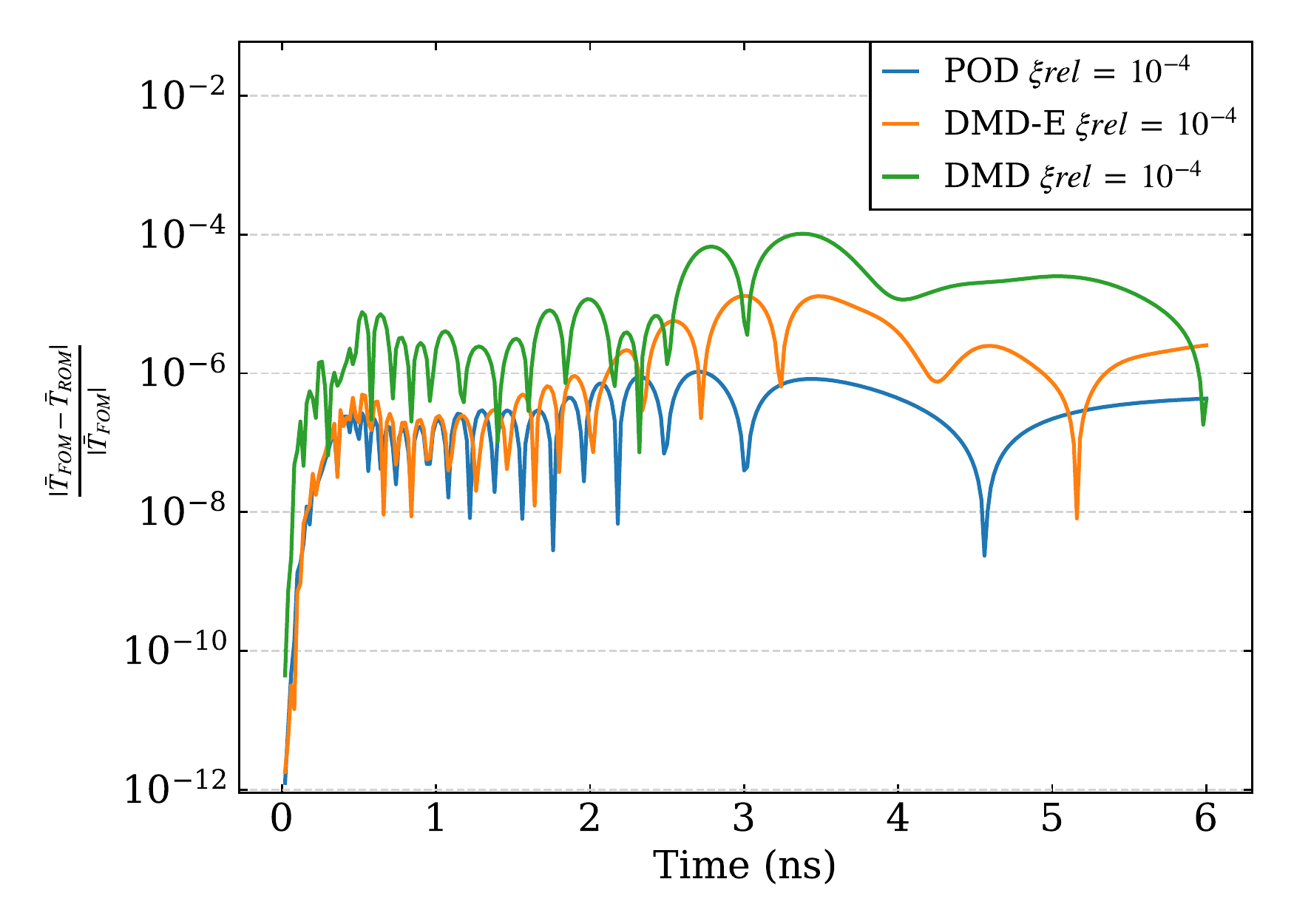}}
	\caption{Relative error for the DDET ROMs with $\xi_{\text{rel}}=10^{-4}$
for data located at and integrated over the right boundary of the domain.
		\label{fig:rbndvals_1e-4} }
\end{figure}

\newpage Figures \ref{fig:rbndvals_1e-2} and \ref{fig:rbndvals_1e-4} plot the relative error in each of these quantities produced by the DDET ROMs using the POD, DMD and DMD-E
with a very low rank corresponding to $\xi_{\text{rel}}=10^{-2}, 10^{-4}$, respectively.
In a similar manner to the results shown above, the DDET ROM with the DMD is observed to reproduce the FOM with the lowest accuracy and the ROMs with the POD and DMD-E achieve similar levels of error to one another. High accuracy is achieved for all considered ROMs with the highest errors on the order of $10^{-2}$. When using the POD or DMD-E even with $\xi_{\text{rel}}=10^{-2}$,
the relative error is about  $10^{-3}$.

These results show that the low-rank DDET ROMs generate good predictions of breakout times across multiple measured quantities.
In the vast majority of cases each of these ROMs yielded the same time step as the FOM for when either $\bar{F}_R$, $\bar{E}_R$, and $\bar{T}_R$ reached a certain arbitrary value.


%=================================================================================
%
%=================================================================================
\newpage
\section{Conclusion} \label{sec:conclusion}
In this paper, new  ROMs for TRT  problems  are presented.
The ROMs are  formulated by a  multilevel system of moment equations
consisting of (i) the multigroup LOQD equations for group radiation energy densities and fluxes
and (ii) the effective grey LOQD equations for total radiation energy density and flux coupled with the MEB equation.
This hierarchy of  equations is derived by means of  the nonlinear-projective  methodology.
The  exact closures of the radiation pressure tensor in the multigroup LOQD equations are formulated by the group Eddington (QD) tensor
 defined by the specific intensity.
The effective grey LOQD equations are closed by means of the spectrum averaged Eddington tensor and opacities.
The proposed ROMs  employ data-driven approximations for the group Eddington tensor   based on available
data to formulate approximate closures. Three data compression techniques are applied, namely, POD, DMD and DMD-E.
The low-rank approximation of the group Eddington tensor  is performed over the whole phase space and time interval.
The hierarchy of  moment equations include radiation energy and momentum balance equations.
These equations are discretized by conservative discretization schemes.
The  solution of the discretized moment equations satisfies the corresponding  conservation of laws.

The analysis of these ROMs has been performed on the classical Fleck-Cummings TRT test problem
with a radiation-driven   Marshak wave.
The  DDET ROMs have been proven effective in efficiently reducing dimensionality of TRT problems, and shown capable of producing a variety of levels of accuracy as the rank parameter is tuned.
Each ROM used the parameter $\xi_{\text{rel}}$ (Eq. \eqref{PODerr_rel}) to determine the rank of approximation. The cases with $\xi_{\text{rel}}=10^{-2},10^{-4}$ lead to the lowest-rank ROMs.
Under these conditions,  DDET ROMs with POD and DMD-E produced comparable levels of accuracy. For all other values of $\xi_{\text{rel}}$ the ROM with POD produced the lowest errors of all the ROMs. Errors in $T$ and $E$ were also observed to converge linearly with $\xi_{\text{rel}}$ for
the ROM with POD. The ROM with DMD consistently produced slightly higher errors compared to the ROM with POD and converged at a similar rate with $\xi_{\text{rel}}$ until $\xi_{\text{rel}}=10^{-10}$, after which numerical errors were seen to dominate. The ROM with DMD-E was shown to possess some numerical instability as $\xi_{\text{rel}}$ decreased and no steady convergence of errors was observed.
All ROMs at low rank were shown to produce errors in $T$ and $E$ that are relatively uniform across both space and time.

The POD ROM closely matches the DMD-E ROM for low-rank approximations, and is more accurate than either DMD ROM for high-rank approximations,  although  the DMD ROMs have the strength of being continuous in time. In order for the POD ROM to be used on time intervals other than those used to generate the training data, an interpolation scheme must be  applied. A possible technique to overcome the shortcomings of both DMD ROMs would be to combine them such that the DMD-E ROM is used for low-rank approximations and the DMD is used instead for high-rank approximations. This would produce a DMD-based DDET ROM with accuracy similar to the POD ROM while using low-rank that is able to converge as the rank is increased.

The broad class of TRT problems involves various parameters, for example, material opacities, incoming radiation fluxes, initial distribution of temperature etc.    The TRT solution depends differently on these parameters.
The next step in future research is to parameterize the ROMs.  In this way the DDET ROMs can be used for efficient parameter sampling for design calculations or experimental simulations. Other  paths for future research include the exploration of other data-based methods of approximating the Eddington tensor and the generation of enhanced databases that allow for lower-rank approximations. An example technique would be to leverage methods of symmetry-reduction \cite{rowley-mardsen-2000,reiss-2004} to improve basis generation given the wave-structure of our considered problems.

The proposed approach for development of ROMs can also be applied to a wide class of multiphysical high-energy density problems, such as radiative hydrodynamics problems. Regardless of the involved multiphysics equations, which can all be coupled to the effective grey low-order system, all that is required for the use of these ROMs is data on the BTE solution.

  \section{Acknowledgements}

The project or effort depicted is sponsored  by the Department of Defense, Defense Threat Reduction Agency, grant number HDTRA1-18-1-0042. The content of the information does not necessarily reflect the position or the policy of the federal government, and no official endorsement should be inferred.

%%Vancouver style references.
\bibliographystyle{model1-num-names}
\bibliography{DDET_rom-jmc-dya}

\end{document}